\documentclass{amsart}

\usepackage{latexsym}

\usepackage{amsfonts,amssymb,euscript,epsfig}

\usepackage[all]{xy}

\newtheorem{theo}{Theorem}[section]
\newtheorem{defi}[theo]{Definition}
\newtheorem{ex}[theo]{Example}

\newtheorem{lem}[theo]{Lemma} 
\newtheorem{prop}[theo]{Proposition}

\newtheorem{rem}[theo]{Remark}
\newtheorem{coro}[theo]{Corollary}

\newcommand{\cgot}{\ensuremath{\mathfrak{c}}}

\newcommand{\sugot}{\ensuremath{\mathfrak{su}}}

\newcommand{\kgot}{\ensuremath{\mathfrak{k}}}

\newcommand{\tgot}{\ensuremath{\mathfrak{t}}}
\newcommand{\Rgot}{\ensuremath{\mathfrak{R}}}
\newcommand{\Sgot}{\ensuremath{\mathfrak{S}}}

\newcommand{\Acal}{\ensuremath{\mathcal{A}}}
\newcommand{\Bcal}{\ensuremath{\mathcal{B}}}
\newcommand{\Ccal}{\ensuremath{\mathcal{C}}}
\newcommand{\Dcal}{\ensuremath{\mathcal{D}}}
\newcommand{\Ecal}{\ensuremath{\mathcal{E}}}

\newcommand{\Fcal}{\ensuremath{\mathcal{F}}}
\newcommand{\Hcal}{\ensuremath{\mathcal{H}}}

\newcommand{\Lcal}{\ensuremath{\mathcal{L}}}
\newcommand{\Qcal}{\ensuremath{\mathcal{Q}}}
\newcommand{\Mcal}{\ensuremath{\mathcal{M}}}
\newcommand{\Ncal}{\ensuremath{\mathcal{N}}}

\newcommand{\Zcal}{\ensuremath{\mathcal{Z}}}
\newcommand{\Ucal}{\ensuremath{\mathcal{U}}}
\newcommand{\Vcal}{\ensuremath{\mathcal{V}}}

\newcommand{\Zbb}{\ensuremath{\mathbb{Z}}}
\newcommand{\Cbb}{\ensuremath{\mathbb{C}}}

\newcommand{\Rbb}{\ensuremath{\mathbb{R}}}
\newcommand{\Nbb}{\ensuremath{\mathbb{N}}}

\newcommand{\Fbb}{\ensuremath{\mathbb{F}}}
\newcommand{\Sbb}{\ensuremath{\mathbb{S}}}

\newcommand{\fpol}{\ensuremath{\mathcal{S}^{*}}}
\newcommand{\hd}{\ensuremath{\rm DH}}
\newcommand{\conv}{\ensuremath{\rm conv}}
\newcommand{\cone}{\ensuremath{\rm Cone}}

\newcommand{\f}{\ensuremath{\mathcal{C}^{\infty}}}
\newcommand{\Hom}{\ensuremath{\rm Hom}}
\newcommand{\fgene}{\ensuremath{\mathcal{C}^{-\infty}}}
\newcommand{\croc}{\ensuremath{\hookrightarrow}}

\newcommand{\indice}{\ensuremath{\hbox{\rm Index}}}

\newcommand{\loc}{\ensuremath{\hbox{\rm P}}}
\newcommand{\p}{\ensuremath{\Delta}}
\newcommand{\mm}{\ensuremath{\hbox{\rm m}}}
\newcommand{\crpt}{\ensuremath{\hbox{\rm Cr}}}
\newcommand{\Vol}{\ensuremath{\hbox{\rm vol}}}
\newcommand{\T}{\ensuremath{\hbox{\bf T}}}
\newcommand{\tr}{\ensuremath{\hbox{\bf Tr}}}
\newcommand{\K}{\ensuremath{\hbox{\bf K}}}
\newcommand{\Eul}{\ensuremath{\hbox{\rm Eul}}}
\newcommand{\Char}{\ensuremath{\hbox{\rm Char}}}
\newcommand{\End}{\ensuremath{\hbox{\rm End}}}
\newcommand{\Thom}{\ensuremath{\hbox{\rm Thom}}}

\newcommand{\chern}{\ensuremath{\hbox{\rm Ch}}}
\newcommand{\todd}{\ensuremath{\hbox{\rm Td}}}

\newcommand{\stab}{\ensuremath{\hbox{\rm Stab}}}

\newcommand{\kir}{\ensuremath{\hbox{\rm \bf Kir}}}
\newcommand{\CW}{\ensuremath{\hbox{\rm \bf cv}}}

\def \spinc {{\rm Spin}^{c}}
\def \SU {{\rm SU}}

\def \Clif {{\rm Cl}}

\begin{document}
    
\title{Jump formulas in Hamiltonian geometry}

\author{Paul-Emile  PARADAN}

\maketitle

{\center
UMR 5582, Institut Fourier, B.P. 74, 38402, Saint-Martin-d'H\`eres 
cedex, France\\
e-mail: Paul-Emile.Paradan@ujf-grenoble.fr\\
}


{\small
\tableofcontents
}

\section{Introduction}

This paper is concerned with the Hamiltonian actions of a (compact)
torus $T$ on a symplectic manifold $(M,\Omega)$. We assume that the
moment map $\Phi: M\to\tgot^*$ is {\em proper} and that 
the generic stabilizer of $T$ on $M$ is {\em finite}. We are
interested here in two global invariants: 

- the Duistermaat-Heckman measure $\hd(M)$ which is the pushforward by 
$\Phi$ of the Liouville volume form,

- the Riemann-Roch characters $RR(M,L^{\otimes k}),k\geq 1$, which are 
virtual representations of $T$. Here $M$ is compact and the data 
$(M,\Omega,\Phi)$ is prequantized by a Kostant-Souriau line bundle
$L$. For every couple $(\mu,k)\in\Lambda^*\times\Zbb^{>0}$, we 
denote by $\mm(\mu,k)\in\Zbb$ the multiplicity of the weight $\mu$ in 
$RR(M,L^{\otimes k})$.

\medskip

One can associate to a connected component $\cgot$ of regular values
of $\Phi$ the local invariants:

- the Duistermaat-Heckman polynomial $\hd_\cgot:\tgot^*\to\Rbb$ 
which coincides with $\hd(M)$ on $\cgot$,

- the periodic polynomial $\mm_\cgot:\Lambda^*\times\Zbb\to\Zbb$ which 
coincides with the map $\mm:\Lambda^*\times\Zbb^{>0}\to\Zbb$ on
the cone of $\tgot^*\times\Rbb$ generated by $\cgot\times\{1\}$.

\medskip

The main results of this paper concern the differences
$\hd_{\cgot_+}-\hd_{\cgot_-}$ and $\mm_{\cgot_+}-\mm_{\cgot_-}$ when
$\cgot_\pm$ are two {\em adjacent} connected components of 
regular values of $\Phi$. Let us introduce some notations. 
We denote by $\p\subset\tgot^*$ the hyperplane that separates 
$\cgot_\pm$, and by $T_\p\subset T$ the subtorus of dimension $1$ 
that has for Lie algebra the one dimensional subspace $\tgot_\p$ 
which is orthogonal to the direction of $\p$. We make the choice 
of a decomposition  $T=T_\p\times T/T_\p$, 
where $T/T_\p$ denotes a subtorus de $T$. At the level of Lie
algebras, we have then $\tgot=\tgot_\p\oplus(\tgot/\tgot_\p)$ and 
$\tgot^*=\tgot_\p^*\oplus(\tgot/\tgot_\p)^*$:
 hence $\xi + (\tgot/\tgot_\p)^*=\p$ for any $\xi\in\p$.

Let $\cgot'\subset\p$ be
the relative interior of $\overline{\cgot_+}\cap\overline{\cgot_-}$ 
in $\p$. 

\medskip

{\em In order to give a clear idea of our results we suppose in the 
introduction that $\Phi^{-1}(\xi)\cap M^{T_\p}$ is {\em connected}
when $\xi\in\cgot'$. It means that there exists only one connected 
component $Z\subset M^{T_\p}$ such that $\cgot'\subset\Phi(Z)$.  
We denote by $N_Z$ the normal bundle of $Z$ in $M$.
}

\medskip

Let $\Omega^\p_\xi$ be the induced symplectic form on the reduced space 
$\Mcal^\p_\xi:=$ \break $(\Phi^{-1}(\xi)\cap M^{T_\p})/(T/T_\p)$. Let 
$\omega^\p_\xi\in\Hcal^2(\Mcal^\p_\xi)\otimes \tgot/\tgot_\p$ be the
curvature of the $T/T_\p$-principal bundle $\Phi^{-1}(\xi)\cap
M^{T_\p}\to\Mcal^\p_\xi$. Let $\beta\in\tgot_\p$ be the 
primitive vector of the lattice $\ker(\exp:\tgot\to T)$  
which is orthogonal to $\p$ and is pointing out $\cgot_-$. 
We prove in Section \ref{DH-mesures} the following 

\medskip

{\bf Theorem A.}\quad {\em Let $2d$ be the dimension of $\Mcal^\p_\xi$, and
  $2r$ be the rank of the bundle $N_Z\to Z$.  
We have the following equality of polynomials: for $a=a'+ a''\in
(\tgot/\tgot_\p)^*\oplus\tgot^*_\p= \tgot^*$ , we have 
$$
({\hd_{\cgot_+}-\hd_{\cgot_-}})(a+\xi)=
\frac{\vert S_Z\vert^{-1}}{\det^{1/2}_Z(\frac{-\Lcal_\beta}{2\pi})}
\int_{\Mcal^\p_\xi}
e^{\Omega^\p_\xi+\langle a',\omega^\p_\xi\rangle}\hbox{{\bf P}}(a'')
$$
where $\hbox{{\bf P}}:\tgot^*_\p\to\Hcal(\Mcal^\p_\xi)$ is the
polynomial mapping defined by 
$$
{\bf P}(a'')=\sum_{k=0}^{d}
\frac{\alpha_k}{(r-1+k)!}\langle a'',\beta\rangle^{r-1+k}
$$
Here the $\alpha_k\in\Hcal^{2k}(\Mcal^\p_\xi)$ are characteristic
classes and $\alpha_0=1$. In the first equation 
$\det^{1/2}_Z(\frac{-\Lcal_\beta}{2\pi})$ is the
Pfaffian of the infinitesimal action of $\frac{-\beta}{2\pi}$ on 
the fibers of $N_Z$, and $\vert S_Z\vert^{-1}$ 
is the cardinal of the generic stabilizer of $T/T_\p$ on $Z$.
}

\bigskip

Theorem  {\bf A.} generalizes previous results of 
Guillemin-Lerman-Sternberg \cite{GLS96} and Brion-Procesi 
\cite{Brion-Procesi}. In Section \ref{subsec:Jump-dh} we give the 
precise definition of the characteristic classes $\alpha_k$.

\medskip

Suppose now that $M$ is compact and is prequantized by a
Kostant-Souriau line bundle. The hyperplane $\p$ is defined by the
equation $\frac{\langle\xi,\beta\rangle}{2\pi}=r_\p,\ \xi\in\tgot^*$, 
for some $r_\p\in\Zbb$. The bundle $N_Z$ decomposes as the sum 
of two polarized sub-bundles
$N_Z^{\pm,\beta}$. Let $s_Z^\pm$ be the absolute value of
the trace of $\frac{1}{2\pi}\Lcal_\beta$ on $N_Z^{\pm,\beta}$. 
In (\ref{eq:bundle-xi-mu}) we define a familly of orbifold 
vector bundles 
$$
\Sbb^k_{Z,\mu}\to\Mcal^\p_\xi,\quad\quad
(\mu,k)\ \in\ \Lambda^*\times\Zbb,
$$ 
with the fundamental property that $\Sbb^k_{Z,\mu}=0$ when 
$$
- s_Z^- < \frac{\langle\mu,\beta\rangle}{2\pi}- k r_\p < s_Z^+.
$$
We prove in Section  \ref{subsec:jump-formula} the following 

\medskip

{\bf Theorem B.}\quad {\em For all $(\mu,k)\in\Lambda^*\times\Zbb$ 
we have
$$
\mm_{\cgot_+}(\mu,k)-\mm_{\cgot_-}(\mu,k)=RR(\Mcal^\p_\xi,\Sbb^k_{Z,\mu}).   
$$
In particular $\mm_{\cgot_+}(\mu,k)=\mm_{\cgot_-}(\mu,k)$ if
$$
 - s_Z^- < \frac{\langle\mu,\beta\rangle}{2\pi}- k r_\p < s_Z^+. 
$$
The integer $s_Z^+ + s_Z^-$ is larger than half of the codimension of 
$Z$ in $M$. The previous inequalities are optimal, i.e. there exists
$(\mu,k)$ such that 
$\frac{\langle\mu,\beta\rangle}{2\pi}-k r_\p=\pm s_Z^\pm$ and 
$\mm_{\cgot_+}(\mu,k)\neq\mm_{\cgot_-}(\mu,k)$
}

\medskip

In Section \ref{sec:cas-coadjoint} we apply Theorem {\bf B} 
to the particular cases where $M$ is a integral coadjoint orbit of a 
compact Lie group $G$. In Section \ref{subsec:su-n}, we study 
more precisely the case $G=\SU(n)$: here our result precises some of the
results of Billey-Guillemin-Rassart \cite{BGR03}.

\medskip

In Section \ref{sec:cas-C-d}, we study the case where $M$ is a complex 
vector space. In this situation, we interpret Theorem {\bf B} 
as a combinatorial formula between vector partition
functions. We recover also a Theorem of Szenes-Vergne \cite{Szenes-Vergne02}.

\bigskip

{\bf Acknowledgments.}\ I am very grateful to Mich\`ele Vergne for
bringing me the reference \cite{BGR03} to my attention, and for
explaining me her work with Andr\'as Szenes \cite{Szenes-Vergne02}.

\medskip

\begin{center}
    {\bf Notations}
\end{center}

{\em Throughout the paper $T$ will denote a compact, connected abelian Lie group, 
and $\tgot$ its Lie algebra. The integral lattice 
$\Lambda\subset\tgot$ is defined as the kernel of $\exp:\tgot\to T$, 
and the real weight lattice $\Lambda^* \subset\tgot^*$ is defined by : 
$\Lambda^*:=\hom(\Lambda,2\pi\Zbb)$. Every $\mu\in\Lambda^*$ defines a 
1-dimensional $T$-representation, denoted by $\Cbb_{\mu}$, where 
$t=\exp X$ acts by $t^{\mu}:=e^{i\langle\mu,X\rangle}$. 
We denote by $R(T)$ the ring of characters of finite-dimensional 
$T$-representations. We denote by $R^{-\infty}(T)$  the 
set of generalized characters of $T$. An element 
$\chi\in R^{-\infty}(T)$ is of the form 
$\chi=\sum_{\mu\in\Lambda^{*}}a_{\mu}\, \Cbb_{\mu}$, 
where $\mu\mapsto a_{\mu}, \Lambda^{*}\to\Zbb$ has at most 
polynomial growth. 

The symplectic manifolds are oriented by their Liouville volume forms. 
If $(Z,o_Z)$ is an oriented  submanifold of an
oriented manifold $(M,o_M)$, we take on the fibers of the normal
bundle $N$ of $Z$ in $M$, the orientation $o_N$ satisfying
$o_M=o_Z\cdot o_N$. }


\section{Duistermaat-Heckman measures}\label{DH-mesures}

Let $(M,\Omega)$ be a symplectic manifold of dimension $2n$ equipped 
with an Hamiltonian action of a torus $T$, with Lie algebra
$\tgot$. The moment map $\Phi: M\to\tgot^*$ satisfies the relations
$\Omega(X_M,-)+d\langle\Phi,X\rangle=0$, $X\in\tgot$. We assume in this section 
that $\Phi$ is {\em proper}, and that the generic stabiliser $\Gamma_M$ of $T$ on $M$ 
is {\em finite}.

The Duistermaat-Heckman measure $\hd(M)$ is defined as the pushforward by $\Phi$ 
of the Liouville volume form $\frac{\Omega^n}{n !}$ on $M$.  
For every $f\in\f(\tgot^*)$ with compact
support one has $\int_{\tgot^*}{\hd(M)}(a)f(a)=\int_M f(\Phi)\frac{\Omega^n}{n !}$. 
In other terms ${\hd(M)}(a)= \int_M \delta(a-\Phi)\frac{\Omega^n}{n !}$. We can define
$\hd(M)$ in terms of equivariant forms as follows. Let $\Acal(M)$ 
be the space of differential forms on $M$ with complex
coefficients. We denote by $\Acal^{-\infty}_{temp}(\tgot,M)$ the space 
of tempered generalized functions 
over $\tgot$ with values in $\Acal(M)$, and by $\Mcal^{-\infty}_{temp}(\tgot^*,M)$ 
the space of tempered distributions over $\tgot^*$ with values 
in $\Acal(M)$. Let $\Fcal:\Acal^{-\infty}_{temp}(\tgot,M)\to 
\Mcal^{-\infty}_{temp}(\tgot^*,M)$ be the Fourier transform 
normalized by the condition that $\Fcal(X\mapsto e^{i\langle \xi,X\rangle})$ is 
equal to the Dirac distribution $a\mapsto\delta(a-\xi)$.

Let $\Omega_\tgot(X)=\Omega -\langle\Phi,X\rangle$ be the equivariant symplectic form.
We have then $\Fcal(e^{-i\Omega_\tgot})=e^{-i\Omega}\delta(a-\Phi)$ and so 
  \begin{equation}
    \label{eq:dh-equi}
    {\hd(M)}=(i)^n \int_M \Fcal(e^{-i\Omega_\tgot}).
  \end{equation}

\subsection{Equivariant cohomology and localization}
We first recall the Cartan model of equivariant cohomology  
with polynomial coefficients and the extension to generalized  coefficients 
defined by Kumar and Vergne \cite{Kumar-Vergne}. We give after a brief account 
to the method of localization developped in \cite{pep1,pep2},

Let $M$ be a manifold provided with an action of a compact connected
Lie group $K$ with Lie algebra $\kgot$. Let $d:\Acal(M)\to\Acal(M)$ be 
the exterior differentiation. 
Let $\Acal_{c}(M)$ be the sub-algebra of compactly supported
differential forms. If $\xi$ is a vectors field on $M$ we denote by 
$c(\xi):\Acal(M)\to\Acal(M)$ the contraction by $\xi$. The action of 
$K$ on $M$ gives a morphism $X\to X_{M}$ from $\kgot$ to the Lie algebra
of vectors fields on $M$.

We consider the space of $K$-equivariant maps
$\kgot\to\Acal(M),\ X\mapsto\eta(X)$, equipped with the derivation 
$(D\eta)(X):= (d-c(X_{M}))(\eta(X)),\ X\in\kgot$.
Since $D^2=0$, one can define the cohomology space $\ker D/{\rm Im} D$. 
The Cartan model \cite{B-G-V,Guill-Stern00} 
considers {\em  polynomial}  maps and the
associated cohomology is denoted $\Hcal_{K}^{*}(M)$. 
Kumar and Vergne \cite{Kumar-Vergne} 
studied the cohomology spaces $\Hcal_{K}^{\pm\infty}(M)$ obtained by taking
$\Ccal^{\pm\infty}$ maps. Recall the 
construction $\Hcal_{K}^{-\infty}(M)$.

The space $\Ccal^{-\infty}(\kgot,\Acal(M))$ of generalized functions on
$\kgot$ with values in the space $\Acal(M)$ is, by definition,
 the space ${\Hom}({\hbox{\sl m}}_c(\kgot),\Acal(M))$ of continuous $\Cbb$-linear maps
from the space ${\hbox{\sl m}}_c(\kgot)$ of smooth compactly supported densities on 
$\kgot$ to the space $\Acal(M)$, both endowed with the $\f$-topologies. We define 
$\Acal_K^{-\infty}(M):=\Ccal^{-\infty}(\kgot,\Acal(M))^K$ 
as the space of $K$-equivariant $\Ccal^{-\infty}$-maps from $\kgot$ 
to $\Acal(M)$. 
The differential $D$ defined on $\Ccal^{\infty}(\kgot,\Acal(M))$ 
admits a natural extension to $\Ccal^{-\infty}(\kgot,\Acal(M))$ and $D^2=0$ on
$\Acal_K^{-\infty}(M)$ \cite{Kumar-Vergne}.  
The cohomology associated to $(\Acal_K^{-\infty}(M),D)$ is
called the $K$-equivariant cohomology with generalized coefficients 
and is denoted by 
$\Hcal_K^{-\infty}(M)$. The subspace $\Acal_{K,c}^{-\infty}(M):=
\Ccal^{-\infty}(\kgot,\Acal_{c}(M))^K$ is stable under the differential $D$, and we
denote by $\Hcal_{K,c}^{-\infty}(M)$ the associated cohomology. When $M$ is 
oriented, the integration over $M$ gives rise to a map $\int_M:
\Hcal_{K,c}^{-\infty}(M)\to\Ccal^{-\infty}(\kgot)^K$.

\medskip

{\em Localization procedure.} Let $\lambda$ be a $K$-invariant
$1$-form on $M$ and let 
\begin{equation}
  \label{eq:phi-lambda}
  \Phi_\lambda:M\to\kgot^*
\end{equation}
be the $K$-equivariant map defined by $\langle\Phi_\lambda(m),X\rangle=
\lambda(X_M)_m$ : then 
$D\lambda(X)=d\lambda-\langle\Phi_\lambda,X\rangle$.  The localization 
procedure developped in \cite{pep1,pep2} is based on 
the existence of an inverse $[D\lambda]^{-1}$ of the $K$-equivariant
form $D\lambda$. It is an equivariantly closed element of 
$\Acal_K^{-\infty}(M- \Phi_\lambda^{-1}(0))$
defined by the integral
\begin{equation}
  \label{eq:lambda-inverse}
  [D\lambda]^{-1}(X)=i\int_0^{\infty}e^{-i\, t\, D\lambda(X)}dt.
\end{equation}
An open subset $\Ucal\subset M$ is called {\em adapted to} $\lambda$ 
if $\Ucal$ is $K$-invariant and if $(\partial\Ucal)
\cap\Phi_\lambda^{-1}(0)=\emptyset$. In \cite{pep2},
we associate to an open subset $\Ucal$ adapted to $\lambda$, the following 
equivariantly closed form with generalized  coefficients 
\begin{equation}
  \label{eq:def-P}
\loc_{\lambda}^\Ucal=\chi^\Ucal + d\chi^\Ucal[D\lambda]^{-1}\lambda\ .
\end{equation}
Here $\chi^\Ucal\in\f(M)$ is a $K$-invariant function supported in 
$\Ucal$ which is equal to $1$ in a neighborhood of $\Ucal\cap\Phi_\lambda^{-1}(0)$.
The cohomology class defined by $\loc_{\lambda}^\Ucal$ in 
$\Hcal_K^{-\infty}(M)$ does  not depend of $\chi^\Ucal$ (in particular 
$\loc_{\lambda}^\Ucal=0$ in $\Hcal_K^{-\infty}(M)$ if 
$\Ucal\cap\Phi_\lambda^{-1}(0)=\emptyset$). If $\Ucal\cap
\Phi_\lambda^{-1}(0)$ is compact, we take  $\chi^\Ucal$ with compact
support, then $\loc_{\lambda}^\Ucal$ defines a cohomology class in 
$\Hcal_{K,c}^{-\infty}(M)$.

\subsection{Localization of $\hd(M)$}\label{sec:hd}

We come back to the situation of a Hamiltonian action of a torus 
$T$ on a symplectic manifold $(M,\omega)$. We keep the same
notations and hypothesys of the introduction. We need two 
auxilliary data : a $T$-invariant Riemannian metric on $M$ denoted  
$(-,-)_{_M}$, and a scalar product $(-,-)$ on $\tgot^*$ 
which induces an identification $\tgot^*\simeq\tgot$.  
  
Let $\Hcal$ be the Hamiltonian vectors field of the function 
$\frac{-1}{2}\|\Phi\|^2: M\to \Rbb$ : for every $m\in M$ we have 
$\Hcal_m=(\Phi(m))_M\vert_m$. Then for every $\xi\in\tgot^*$,  
the Hamiltonian vectors field of $\frac{-1}{2}\|\Phi-\xi\|^2$ is 
$\Hcal-\xi_M$. For every $\xi\in\tgot^*$, 
we consider the following  $T$-invariant $1$-form
\begin{equation}
  \label{eq:lambda-xi}
  \lambda_\xi=(\Hcal-\xi_M,-)_{M}
\end{equation}
and the corresponding map $\Phi_{\lambda_\xi}:M\to\tgot^*$ 
(see (\ref{eq:phi-lambda})). Here $\Phi_{\lambda_\xi}^{-1}(0)$ 
coincides with the subset $\crpt(\|\Phi-\xi\|^2)\subset M$ of 
critical points of the function $\|\Phi-\xi\|^2$, and \break 
$m\in \crpt(\|\Phi-\xi\|^2)$ if and only if $(\Phi(m)-\xi)_M$ 
vanishes at $m$ \cite{pep1,pep2}. 
\begin{defi}
  \label{def:p-xi}
Let $\loc_\xi\ \in\ \Hcal_{T,c}^{-\infty}(M)$ be the cohomology 
class defined by $\loc_{\lambda_\xi}^\Ucal$, where $\Ucal$ is 
a $T$-invariant relatively compact 
neighborhood of $\Phi^{-1}(\xi)$ such that \break 
$\overline{\Ucal}\cap \crpt(\|\Phi-\xi\|^2)=\Phi^{-1}(\xi)$.
\end{defi}

The cohomology class $\loc_\xi$ will be used to localized the 
Duitermaat-Heckman measure. For every $\xi\in\tgot^*$, we define 
the distribution $\hd_\xi(M)$ by 
\begin{equation}
  \label{eq:dh-xi}
  {\hd_\xi(M)}=(i)^n \Fcal\left(\int_M \loc_\xi e^{-i\Omega_\tgot}\right).
\end{equation}
Here we can put the Fourier transform outside the integral because 
$\loc_\xi$ is compactly 
supported on $M$. For any $\xi\in\tgot^*$ 
let $r_\xi>0$ be the smallest non-zero critical value of the function 
$\|\Phi-\xi\|^2$. As a particular case of 
Proposition 3.8 in \cite{pep2}, we have

\begin{prop}
  \label{prop:dh-loc}
Let $\xi$ be any point in $\tgot^*$. The following equality of distributions 
on $\tgot^*$ 
$$
\hd(M)=\hd_\xi(M)
$$
holds in the open ball $B(\xi,r_\xi)\subset \tgot^*$.
\end{prop}

We will now use the last Proposition, first to recover the classical result of 
Duistermaat and Heckman \cite{D-H} concerning the polynomial behaviour of 
$\hd(M)$ on the open subset of regular values of $\Phi$. After we determine 
the difference taken by $\hd(M)$ between two adjacent regions of regular values.

\subsection{Polynomial behaviour}

We recall now the computation of the cohomology class $\loc_\xi$ when $\xi$ is a 
regular value of $\Phi$, that is given in \cite{pep1}[Section 6] for the 
torus case (and in \cite{pep2} [Section 3.1] for the case of Hamiltonian 
action of a compact Lie group). First recall the following 
basic result which shows that $\xi\mapsto\hd_\xi(M)$ is 
locally constant on the open subset of regular values of $\Phi$.
\begin{lem}[\cite{cras04}]\label{lem:P-xi-0}
  If $\xi$ and $\xi'$ belong to the same connected component of
  regular values of $\Phi$, we have $\loc_\xi=\loc_{\xi'}$ in $\Hcal_{T,c}^{-\infty}(M)$.
\end{lem}

Associated to a regular value of $\xi$, we have the $T$-principal bundle 
$\Phi^{-1}(\xi)\to \Mcal_\xi:=\Phi^{-1}(\xi)/T$ with curvature form 
$\omega_\xi\in\Hcal^2(\Mcal_\xi)\otimes\tgot$. The orbifold $\Mcal_\xi$ 
carries a canonical symplectic $2$-form $\Omega_\xi$. We denote 
$\kir_{\xi}:\Hcal_T^{\infty}(M)\to\Hcal^*(\Mcal_\xi)$ the Kirwan morphism. 
For any $\psi\in\f(\tgot)$ and $\eta\in\Hcal_T^{\infty}(M)$ we have 
$\kir_{\xi}(\eta\psi)=\kir_{\xi}(\eta)\psi(\omega_\xi)$, where the 
characteristic class 
$\psi(\omega_\xi)$ is the value of the 
differential operator $e^{\omega_\xi(\frac{\partial}{\partial X}|_0)}$ 
against $\psi$. After \cite{pep2}[Prop. 3.11],  the integral 
$\int_\tgot\int_{M}\loc_\xi(X)\eta(X)\psi(X)dX$ is equal to 
\begin{equation}
  \label{eq:loc-xi-psi}
\frac{(-2i\pi)^{dimT}{\rm vol}(T,dX)}{|\Gamma_M|}\int_{\Mcal_\xi}
\kir_{\xi}(\eta)\psi(\omega_\xi)  
\end{equation}
for every equivariant class $\eta\in\Hcal_{T}^{\infty}(M)$. Here 
$\Vol(T,dX)$ is the volume of $T$ for the Haar mesure compatible with $dX$, and
$\vert \Gamma_M\vert$ is the cardinal of $\Gamma_M$ (Note that 
the generic stabilizer of $T$ on $\Phi^{-1}(\xi)$ is $\Gamma_M$). In other words, 
for every $\eta\in\Hcal_{T}^{\infty}(M)$  we have the following equality 
of generalized functions on $\tgot^*$ supported at $0$
\begin{equation}
  \label{eq:formule-p}
\int_{M}\loc_\xi(X)\eta(X)=\frac{(-2i\pi)^{dimT}}{|\Gamma_M|}
\int_{\Mcal_\xi}\kir_{\xi}(\eta)
e^{\omega_\xi(\frac{\partial}{\partial X}|_0)}\Vol(T,-).
\end{equation}

For $\eta=e^{-i\Omega_\tgot}$ we have $\kir_\xi(\eta)=
e^{-i(\Omega_\xi-\langle\xi,\omega_\xi\rangle)}$, and a 
small computation shows that 
\begin{equation}
  \label{eq:fourier-omega}
  \Fcal\left(e^{\omega_\xi(\frac{\partial}{\partial X}|_0)}\Vol(T,-)\right)(a)=
e^{-i\langle a ,\omega_\xi\rangle}  \frac{da}{(2\pi)^{dimT}}, \quad a\in\tgot^*.  
\end{equation}
where $da$ is the Lebesgue measure on $\tgot^*$ normalized by the condition:  
$\Vol(T,dX)=1$ for the Lebesgue measure $dX$ on $\tgot$ which is dual to $da$.

Finally (\ref{eq:dh-xi}), (\ref{eq:formule-p}) and (\ref{eq:fourier-omega}) 
give  
\begin{eqnarray}
  \label{eq:dh-xi-calcul}
   {\hd_\xi(M)}(a)&=& \frac{(i)^{p}}{|\Gamma_M|}\int_{\Mcal_\xi}
e^{-i(\Omega_\xi+\langle a-\xi,\omega_\xi\rangle)} \ da \nonumber\\
&=&\frac{1}{|\Gamma_M|}\int_{\Mcal_\xi}
\frac{(\Omega_\xi+\langle a-\xi,\omega_\xi\rangle)^p}{p!}\ da,
\end{eqnarray}
where $2p= {\rm dim}\Mcal_\xi$. 

\begin{defi}\label{def:DH-c}
  For any connected component $\cgot$ of regular values of $\Phi$ we denote
$\hd_{\cgot}$ the polynomial function $a\mapsto\frac{1}{|\Gamma_M|}\int_{\Mcal_\xi}
\frac{(\Omega_\xi+\langle a-\xi,\omega_\xi\rangle)^p}{p!}$, where $\xi$ is any point 
of $\cgot$. Hence, if $\xi\in\cgot$ we have the equality 
$\hd_\xi(M)(a)=\hd_{\cgot}(a) da,\ a\in\tgot^*$, 
where $da$ is the Lebesgue measure on $\tgot^*$ normalized by the condition:  
$\Vol(T,dX)=1$ for the Lebesgue measure $dX$ on $\tgot$ which is dual to $da$.
\end{defi}

With the help of Proposition \ref{prop:dh-loc} we recover 
the classical result of Duistermaat and Heckman that says that the measure 
$\hd(M)$ is locally 
polynomial\footnote{It is a polynomial times a Lebesgue measure on 
$\tgot^*$.} on the open subset of regular values of $\Phi$, and 
it's value at a regular element $\xi$ is equal to the 
symplectic volume of the reduce space $\Mcal_\xi$. More 
precisely we have shown that for a connected component 
$\cgot$  of regular values of $\Phi$ we have
\begin{equation}
  \label{eq:f-C}
{\hd(M)}(a)={\hd_{\cgot}}(a) d a,\qquad a\in\cgot. 
\end{equation}

\subsection{Jump formulas}\label{subsec:Jump-dh} 
Consider now two connected regions $\cgot_\pm$ 
of regular values of $\Phi$ separated by an hyperplane
$\p\subset\tgot^*$. In this section 
we compute the polynomial $\hd_{\cgot_+}-\hd_{\cgot_-}$. It generalizes 
previous results of Guillemin-Lerman-Sternberg \cite{GLS96} and Brion-Procesi 
\cite{Brion-Procesi}.

Let $\xi_+,\xi_-$ be  respectively two elements of $\cgot_+$ and $\cgot_-$. 
We know from (\ref{prop:dh-loc}), (\ref{eq:dh-xi-calcul}) and 
Definition (\ref{def:DH-c}) that
\begin{equation}
  \label{eq:f-difference}
 {(\hd_{\cgot_+}}-{\hd_{\cgot_-})}(a)da=
(i)^n \Fcal\left(\int_M (\loc_{\xi_+}-\loc_{\xi_-}) 
e^{-i\Omega_\tgot}\right)(a),\quad a\in\tgot^*. 
\end{equation}
We recall now the computation of the cohomology class $\loc_{\xi_+}-\loc_{\xi_-}\in
\Hcal^{-\infty}_{T,c}(M)$ \cite{cras04}.

Let $T_{\p}\subset T$ be the subtorus of dimension $1$, 
with Lie algebra $\tgot_{\p}:=
\{X\in\tgot\vert\ \langle \xi-\xi',X\rangle=0,\ \forall \xi,\xi'\in\p\}$. 
We make the choice of a decomposition  $T=T_\p\times T/T_\p$, 
where $T/T_\p$ denotes a subtorus de $T$. At the level of Lie
algebras, we have then $\tgot=\tgot_\p\oplus(\tgot/\tgot_\p)$ and 
$\tgot^*=\tgot_\p^*\oplus(\tgot/\tgot_\p)^*$:
 we see that $\xi + (\tgot/\tgot_\p)^*=\p$ for any $\xi\in\p$.

Let $\cgot'$ be the relative interior of 
$\overline{\cgot_+}\cap\overline{\cgot_-}$ in $\p$. It is not 
difficult to see that for any $m\in\Phi^{-1}(\cgot')$ the stabilizer
$\tgot_m\subset\tgot$ is either equal to $\tgot_\p$ or reduced to $\{0\}$.

\begin{defi}\label{def:M-delta-quantum}
 We denote $M_{\cgot'}$ the union of 
the connected component $Z$  of the fixed point set 
$M^{T_{\p}}$ for which we have $\cgot'\subset\Phi(Z)$.
 Let $M_{\cgot'}^o$ be the $T$-invariant open subset of 
$M_{\cgot'}$ where $T/T_\p$ acts locally freely.
\end{defi}

The symplectic manifold $M_{\cgot'}$ carries a Hamiltonian action of 
$T/T_{\p}$ with moment map $\Phi_{\cgot'}:M_{\cgot'}\to\p$ equal to the
restriction of $\Phi$ on $M_{\cgot'}$. Let $\xi$ be a point in $\cgot'$. 
We remark before that for all $m\in\Phi^{-1}(\xi)$ the 
stabilizer $\tgot_m$ is either equal to $\tgot_\p$ or to $\{0\}$: in particular
$\xi$ is a regular value of $\Phi_{\cgot'}$, i.e. $\Phi^{-1}(\xi)\in
M_{\cgot'}^o$. Following Definition \ref{def:p-xi} 
we associate to $\xi$ the cohomology 
class $\loc_\xi^\p\in\Hcal^{-\infty}_{T/T_\p,c}(M_{\cgot'}^o)$. 

Let $\Hcal^*(M_{\cgot'}^o)^{bas}$ 
be the sub-algebra of $\Hcal^*(M_{\cgot'}^o)$ formed by the $T$-basic elements.
Since the $T_\p$-action on $M_{\cgot'}^o$ is trivial we have a canonical
product operation 
\begin{equation}
  \label{eq:produit}
  \Hcal^{-\infty}_{T/T_\p,c}(M_{\cgot'}^o)
\times\fgene(\tgot_\p, \Hcal^*(M_{\cgot'}^o)^{bas})
\longrightarrow \Hcal^{-\infty}_{T,c}(M_{\cgot'}^o).
\end{equation}
 
\begin{prop}[\cite{cras04}]
  \label{prop:p-diff}
There exists a generalized function supported at $0$,
$\delta^\p\in\fgene(\tgot_\p, \Hcal^*(M_{\cgot'}^o)^{bas})$, such that
$$
 \loc_{\xi^+}-\loc_{\xi^-}=
(i_\p)_*\left(\loc_\xi^\p\delta^\p\right)\qquad {\rm in} 
\quad \Hcal^{-\infty}_{T,c}(M).
$$
Here $(i_\p)_*:\Hcal^{-\infty}_{T,c}(M_{\cgot'}^o)\to\Hcal^{-\infty}_{T,c}(M)$ 
is the direct image map relative to the inclusion 
$i_\p: M_{\cgot'}^o\croc M$.
\end{prop}

We will now give the precise definition of $\delta^\p$. 
The decomposition $T=T_\p\times T/T_\p$ and the trivial action of $T_\p$
on $M_{\cgot'}^o$ determine a canonical isomorphism
$$
j_\p:\Hcal^*_{T}(M_{\cgot'}^o)\stackrel{\sim}{\longrightarrow}
\fpol(\tgot_\p)\otimes\Hcal^*_{T/T_{\p}}(M_{\cgot'}^o),
$$ 
where $\fpol(\tgot_\p)$ is the algebra of complex polynomial 
functions on $\tgot_\p$. Since the $T/T_\p$-action on $M_{\cgot'}^o$ 
is locally free, we have the Chern-Weil isomorphism  
$$
\CW_\p:\Hcal^*_{T/T_{\p}}(M_{\cgot'}^o)
\stackrel{\sim}{\longrightarrow}\Hcal^*(M_{\cgot'}^o)^{bas}.
$$ 

Let $N_{\p}$ be the $T$-equivariant normal bundle of $M^{T_\p}$ in
$M$, and  let $\Eul{}(N_{\p})\in\Hcal_{T}^*(M^{T_\p})$ 
be the $T$-equivariant Euler class  of $N_\p$. Now we consider  the
restriction of $\Eul{}(N_{\p})$ on the open subset $M_{\cgot'}^o\subset
M^{T_\p}$, that we look through the isomorphism $\CW_\p\circ j_\p$ 
as an element of $\fpol(\tgot_\p)\otimes\Hcal^*(M_{\cgot'}^o)^{bas}$ 
(for simplicity we keep the same notations $\Eul{}(N_{\p})$ 
for this element). Following \cite{pep1}, we define inverses 
$\Eul{}^{-1}_{\pm}(N_{\p})\in\fgene(\tgot_\p,\Hcal^*(M_{\cgot'}^o)^{bas})$ 
by 
\begin{equation}
  \label{eq:eul-1}
  \Eul{}^{-1}_{\pm\beta}(N_{\p})(X) 
=\lim_{s\to +\infty}\frac{1}{\Eul{}(N_{\p})(X\pm is\beta)}.
\end{equation}
Here $\beta\in\tgot_\p-\{0\}$ is choosen such that 
$\langle\xi^+ - \xi^-,\beta\rangle>0$. 

\begin{defi}
The generalized function $\delta^\p\in\fgene(\tgot_\p,\Hcal^*(M_{\cgot'}^o)^{bas})$ 
is defined by
\begin{equation}
  \label{eq:def-delta}
  \delta^\p:=\Eul{}^{-1}_{\beta}(N_{\p})-\Eul{}^{-1}_{-\beta}(N_{\p}).
\end{equation}
Since the polynomial $\Eul{}(N_{\p})$ is  invertible in a smooth manner
on $\tgot_\p-\{0\}$ the generalized function $\delta^\p$ is supported at $0$.
If we restrict $\delta^\p$ to the submanifold $\Phi_{\cgot'}^{-1}(\xi)$ we get 
the generalized function
$\delta^\p_\xi\in\fgene(\tgot_\p,\Hcal^*(\Mcal^\p_\xi))$, 
where $\Mcal^\p_\xi$ denotes the reduced space 
$\Phi_{\cgot'}^{-1}(\xi)/(T/T_\p)=(M^{T_\p}\cap\Phi^{-1}(\xi))/(T/T_\p)$.
\end{defi}

Now we are able to compute the RHS of (\ref{eq:f-difference}). Let 
$\omega^\p_\xi\in\Hcal^2(\Mcal^\p_\xi)\otimes\tgot/\tgot_\p$ be the curvature of the 
$T/T_\p$-principal bundle  $\Phi_{\cgot'}^{-1}(\xi)\to\Mcal^\p_\xi$. Let 
$\vert S^\p_\xi\vert$ be locally constant function on $\Phi_{\cgot'}^{-1}(\xi)$  
which is equal to the cardinal of the generic stabilizer of $T/T_\p$.  
From (\ref{eq:formule-p}) and Proposition \ref{prop:p-diff} we have
\begin{eqnarray}
  \label{eq:diff-1}
 \lefteqn{\int_M (\loc_{\xi_+}-\loc_{\xi_-})(X)e^{-i\Omega_\tgot(X)}}\nonumber\\
& & =\int_{M_{\cgot'}^o}\loc_\xi^\p(X')\delta^\p(X'')e^{-i\Omega_\tgot(X'+X'')}\nonumber\\
& & = \frac{(-2i\pi)^{dim T-1}}{|S^\p_\xi|}\int_{\Mcal^\p_\xi}
e^{\omega^\p_\xi(\frac{\partial}{\partial X'}|_0)}\Vol(T/T_\p,-)
\kir^\p_{\xi}(e^{-i\Omega_\tgot})(X'')
\delta^\p_\xi(X'')
\end{eqnarray}
In the last equation the notations are the following:

- $X=X'+X''$ with $X'\in\tgot/\tgot_\p$ and $X''\in\tgot_\p$, 

- the Kirwan map $\kir^\p_{\xi}: \Hcal^{\infty}_{T}(M)\to 
\f(\tgot_\p,\Hcal^*(\Mcal^\p_\xi))$ is the composition of the 
restriction $\Hcal^{\infty}_{T}(M)\to\Hcal^{\infty}_{T}(\Phi_{\cgot'}^{-1}(\xi))$
with the Chern-Weil isomorphism \break $\Hcal^{\infty}_{T}(\Phi_{\cgot'}^{-1}(\xi))
\stackrel{\sim}{\to}\f(\tgot_\p,\Hcal^*(\Mcal^\p_\xi))$.

A direct computation gives that 
$\kir^\p_{\xi}(\Omega_\tgot)(X'')=\Omega_\xi^\p-\langle\xi,\omega^\p_\xi+X''\rangle$
where $\Omega_\xi^\p$ is the induced symplectic form on the reduced space 
$\Mcal^\p_\xi$. If we take the fourier transform 
in (\ref{eq:diff-1}) we get
\begin{eqnarray}
  \label{eq:diff-2}
\lefteqn{{(\hd_{\cgot_+}-\hd_{\cgot_-})}(a)da}\nonumber\\
&= &  \frac{(i)^{n+1-dim T }}{|S^\p_\xi|}\left(\int_{\Mcal^\p_\xi}
e^{-i(\Omega^\p_\xi+\langle a',\omega^\p_\xi\rangle)}\, da'\, 
\Fcal_{\tgot_\p}(\delta^\p_\xi)(a'')\right)(a-\xi),\nonumber\\
&= & \sum_{Z} \frac{(i)^{n+1-dim T }}{|S^Z_\xi|}\left(\int_{\Zcal_\xi}
e^{-i(\Omega^Z_\xi+\langle a',\omega^Z_\xi\rangle)}\, da'\, 
\Fcal_{\tgot_\p}(\delta^Z_\xi)(a'')\right)(a-\xi)
\end{eqnarray}
where $a=a'+a''$ with $a'\in(\tgot/\tgot_\p)^*$ and
$a''\in(\tgot_\p)^*$. In (\ref{eq:diff-2}), we write
$\int_{\Mcal^\p_\xi}=\sum_{Z}\int_{\Zcal^\p_\xi}$ where the sum is 
taken over the connected components $Z$ of
$M_{\cgot'}^o$, and $\Zcal_\xi=(Z\cap\Phi^{-1}(\xi))/(T/T_\p)$. The $2$-forms
$\Omega^\p_\xi,\omega^\p_\xi$, the generic stabiliser $S^\p_\xi$, 
the vector bundle $N^\p$, the generalized function $\delta^\p_\xi$ 
restrict to each component $Z$: we denote them respectively 
$\Omega^Z_\xi,\omega^Z_\xi$, $S^Z_\xi$, $N_Z$, 
$\delta^Z_\xi$.

\medskip

We recall now the computation of the Fourier tranform of the inverses  
$\Eul{}^{-1}_{\pm\beta}(N_Z)$ $=\Eul{}^{-1}_{\pm\beta}(N_\p)_{\vert Z}$ 
that is given in \cite{pep1}[Proposition 4.8.]. We consider a 
$T$-invariant scalar product on the fibers of the bundle 
$N_\p$. Let $R\in\Acal^2(M_{\cgot'}^o,{\rm so}(N_\p))^{bas}$ be 
the curvature of a $T$-invariant and $T/T_\p$-horizontal 
Euclidean connexion on $N_\p$: we denote by 
$R^Z\in\Acal^2(Z,{\rm so}(N_Z))^{bas}$ the restriction of 
$R$ to a component $Z$. The curvature commutes with the 
infinitesimal action $\Lcal_X$ of $X\in\tgot_\p$, and 
with the complex structure $J_\beta=\Lcal_\beta(-\Lcal_\beta^2)^{1/2}$ 
on $N_\p$ defined by $\beta\in\tgot_\p$.

We denote by $S^\bullet$ the symmetric algebra of the complex 
vector bundle  $(N_\p,J_\beta)$. We keep the same notation 
for the restriction of $S^\bullet$ on the submanifolds $Z$, 
$\Phi_{\cgot'}^{-1}(\xi)$, and for the induced orbifold 
vector bundle on the reduced spaces $\Zcal_\xi$ and 
$\Mcal^\p_\xi$. For each $k\in\Nbb$, we denote by 
$\tr_{S^k}$ the trace operator defined
on the complex endomorphism of $S^k$. For a complex 
endomorphism $A$ of $N_\p$, we denote 
by $A^{\otimes k}$ the induced endomorphism on $S^k$. 
For any $X\in\tgot_\p$, the complex endomorphism $\Lcal_X^{-1}R^Z$ 
is symmetric. Hence the trace 
$\tr_{S^k}((\Lcal_X^{-1}R^Z)^{\otimes k})$ is a basic {\em real} differential
form of degree $2k$ on $Z$ which does not depend of the choice of 
complex structures ($J_\beta$ or $J_{-\beta}$). 

\begin{prop}[\cite{pep1}]
For a smooth function $f$ on $\tgot_\p^*$ with compact support we have
$\int_{\tgot_\p^*}\Fcal_{\tgot_\p}(\Eul{}^{-1}_{\beta}(N_Z))(a'')f(a'')=\int_0^\infty
\hbox{\rm P}_Z(t)f(t\beta^*)dt$ 
where $\hbox{P}_Z$ is the polynomial on $\Rbb$ defined by:
\begin{equation}
  \label{eq:p-beta}
  \hbox{\rm P}_Z(t)=
\frac{(2\pi i)^{r_Z}}{\det^{1/2}_Z(\Lcal_\beta)}\left(\frac{t^{r_Z-1}}{(r_Z-1)!}+
\sum_{k=1}^{dim(Z)/2} \,(i)^k \,\tr_{S^k}((\Lcal_\beta^{-1}R^Z)^{\otimes k})
\frac{t^{r_Z-1+k}}{(r_Z-1+k)!}
\right).
\end{equation}
Here $\det^{1/2}_Z(\Lcal_\beta)$ is the Pfaffian of $\Lcal_\beta$ on
$N_Z$, and $r_Z=\hbox{\rm rk}_\Cbb(N_Z)$. For *
$\Fcal_{\tgot_\p}(\Eul{}^{-1}_{-\beta}(N_Z))$ 
we have 
\begin{eqnarray*}
  \int_{(\tgot_\p)^*}\Fcal_{\tgot_\p}(\Eul{}^{-1}_{-\beta}(N_Z))(a'')f(a'')
&=& \int_0^\infty -\hbox{\rm P}_Z(-t)f(-t\beta^*)dt\\
&=& -\int_{-\infty}^0 \hbox{\rm P}_Z(t)f(t\beta^*)dt.
\end{eqnarray*}
Hence the distribution $\Fcal_{\tgot_\p}(\delta^Z)$ is equal to 
$\hbox{\rm P}_Z(\beta)d\beta$. 
\end{prop}

\medskip

Let $R^Z_\xi$ be the restriction of the curvature $R^Z$
to the submanifold $Z\cap\Phi^{-1}(\xi)$. 
Since $R^Z$ is $T/T_\p$-basic, $\tr_{S^k}((\Lcal_\beta^{-1}R^Z_\xi)^{\otimes k})$ 
can be seen as a real differential form of degree $2k$ on the orbifold 
$\Zcal_\xi=(Z\cap\Phi^{-1}(\xi))/(T/T_\p)$. 

Each connected component $Z$ of $M_{\cgot'}^o$ is a $T/T_\p$ Hamiltonian
manifold: we take for moment map $\Phi_Z: Z\to (\tgot/\tgot_\p)^*$ the
restriction of $\Phi-\xi$ to $Z$. Hence $0$ is a regular value of
$\Phi_Z$. Let $\hd_0(Z)$ be the polynomial function on 
$(\tgot/\tgot_\p)^*=\{a\in\tgot^*\,\vert\,\langle\beta,a\rangle=0\}$ such that 
${\hd(Z)}(a')={\hd_0(Z)}(a')da'$ near $0$. Finally (\ref{eq:diff-2})
together with the last proposition give the following 

\begin{theo}\label{DH-theo-principal}{}

{\bf (a).}We have $\hd_{\cgot_+}-\hd_{\cgot_-}= \sum_{Z\subset
  M_{\cgot'}}\hbox{\bf D}_Z$ where
\begin{equation}
  \label{eq:hd-plus-moins}
 \hbox{\bf D}_Z(a)=\frac{(-2\pi)^{r_Z}}{\det^{1/2}_Z(\Lcal_\beta)}
\left(\sum_{k=0}^{d_Z} \hbox{\rm Q}_{Z,k}\frac{\beta^{r_Z-1+k}}{(r_Z-1+k)!} 
\right)(a-\xi).
\end{equation}
The $\hbox{\rm Q}_{Z,k}$ are the polynomials of degree $d_Z-k$ 
on $(\tgot/\tgot_\p)^*$ defined by
\begin{equation}
  \label{eq:pol-q-k}
  \hbox{\rm Q}_{Z,k}(a')=\frac{(-1)^k}{|S^Z_\xi|}\int_{\Zcal_\xi}
\frac{(\Omega^Z_\xi+\langle a',\omega^Z_\xi\rangle)^{d_Z-k}}{(d_Z-k)!}
\tr_{S^k}((\Lcal_\beta^{-1}R^Z_\xi)^{\otimes k}).
\end{equation}
Here $2d_Z=\dim\Zcal_\xi$ and $2r_Z=\dim M-\dim Z$. 
The polynomial $\hbox{\rm Q}_{Z,0}$
correspond to the Duistermaat-Heckmann polynomial $\hd_0(Z)$. In
particular we see that the polynomial 
$\hd_{\cgot_+}-\hd_{\cgot_-}$ is divisible by the factor
$a\mapsto \langle\beta,\xi-a\rangle^{r-1}$ where $r=\inf_Z r_Z$. 
If $\p\cap\Phi(M)$ is not a facet of the polytope $\Phi(M)$ we have 
$r_Z\geq 2$ for all connected component $Z$ of $M_{\cgot`}$, hence 
$r\geq 2$.

{\bf (b).} Suppose now that $\cgot_-$ is a connected component 
of regular values of $\Phi$ bording a facet $\Phi(M)\cap\p$ 
of the polytope $\Phi(M)$. Here $Z=\Phi^{-1}(\p)$ is a connected 
component of the fixed point set $M^{T_\p}$. In this situation we have 
$\hd_{\cgot_-}= -\hbox{\bf D}_Z$ where the polynomial ${\bf D}_Z$
is defined by (\ref{eq:hd-plus-moins}).

In (\ref{eq:hd-plus-moins}) and (\ref{eq:pol-q-k}) the vector 
$\beta\in\tgot$ is normalized by the following conditions: 
$\beta$ is a primitive vector of the lattice $\ker(\exp:\tgot\to T)$, 
orthogonal to the hyperplane containing 
$\overline{\cgot_+}\cap\overline{\cgot_-}$, and pointing out $\cgot_-$.
\end{theo}

\section{Quantum version of Duistermaat-Heckman measures}
\label{sec:quantum-version}

We suppose here that the Hamiltonian $T$-manifold $(M,\omega,\Phi)$ 
is prequantized by a $T$-equivariant Hermitian line bundle $L$ over
$M$, which is equipped with an Hermitian 
connection $\nabla$ satisfying the Kostant formula
\begin{equation}
  \label{eq:kostant}
  \Lcal(X)-\nabla_{X_M}=i\langle\Phi,X\rangle,\quad X\in\tgot.
\end{equation}
The former equation implies that the first Chern class of $L$ is equal
to $\frac{\Omega}{2\pi}$. In this section we suppose that $M$ is 
{\em compact} and we still assume that the generic 
stabiliser $\Gamma_M$ of $T$ on $M$ is {\em finite}. The quantization of 
$(M,\Omega)$ is defined by the Riemann-Roch character $RR(M,L)\in
R(T)$ which is compute with a $T$-equivariant almost complex stucture 
on $M$ {\em compatible} with $\Omega$. For $k\geq 1$, 
we consider the tensor product $L^{\otimes k}$. Its Riemann-Roch character 
$RR(M,L^{\otimes k})$ decomposes as 
\begin{equation}
  \label{eq:mult-k}
  RR(M,L^{\otimes k})=\sum_{\mu\in\Lambda^*}\mm(\mu,k)\,\Cbb_\mu.
\end{equation}

Let us recall the well-known properties of the map $\mm:\Lambda^*\times
\Zbb^{\geq 0}\to \Zbb$. When $\frac{\mu}{k}$ is a regular value of $\Phi$, the 
"Quantization commutes with Reduction Theorem'' 
\cite{Meinrenken98,Meinrenken-Sjamaar} tell us that 
\begin{equation}
  \label{eq:Q-R}
  \mm(\mu,k)=RR(\Mcal_{\frac{\mu}{k}},\Lcal_\mu^k)
\end{equation}
where 
$\Lcal_\mu^k=(L^{\otimes  k}\vert_{\Phi^{-1}(\frac{\mu}{k})}\otimes\Cbb_{-\mu})/T$ 
is an orbifold line bundle over the symplectic orbifold 
$\Mcal_{\frac{\mu}{k}}=\Phi^{-1}(\frac{\mu}{k})/T$.
In particular if $\frac{\mu}{k}$ does not belong to $\Phi(M)$ 
we have $\mm(\mu,k)=0$. When $\frac{\mu}{k}\in\Phi(M)$ 
is not necessarilly a regular value of $\Phi$, 
one procceed by shift desingularization. If $\xi\in\Phi(M)$ is 
a regular value of $\Phi$ close enough to 
$\frac{\mu}{k}$ then (\ref{eq:Q-R}) becomes 
\begin{equation}
  \label{eq:Q-R-general}
  \mm(\mu,k)=RR(\Mcal_{\xi},\Lcal_{\xi,\mu}^k)
\end{equation}
where $\Lcal_{\xi,\mu}^k=(L^{\otimes k}\vert_{\Phi^{-1}(\xi)}\otimes\Cbb_{-\mu})/T$ 
(for a proof see \cite{Meinrenken-Sjamaar,pep4}). 

\begin{defi}
  A function $f:\Xi\to\Zbb$ defined over a lattice $\Xi\simeq\Zbb^r$ is 
called {\em periodic polynomial} if 
$$
f(x)=\sum_{i=1}^p e^{i\frac{\langle\alpha_j,x\rangle}{N}}P_j(x),\quad x\in\Xi,
$$
where $\alpha_1,\cdots,\alpha_p\in\Xi^*$, $N\geq 1$, and the functions 
$P_1,\cdots,P_p$ are polynomials with complex coefficients.
\end{defi}

\begin{rem}
Let $\Ccal$ a cone with {\em non-empty} interior in the real vector space 
$\Xi\otimes_\Zbb\Rbb$. Any {\em periodic-polynomial} function $f:\Xi\to\Zbb$ is 
completely determined by its restriction on $\Ccal\cap\Xi$.
\end{rem}

Let $\cgot\subset\tgot^*$ be a connected component of regular values of $\Phi$. 
In \cite{Meinrenken-Sjamaar} Meinrenken an Sjamaar proved that there exits a 
periodic polynomial function $\mm_{\cgot}:\Lambda^*\times\Zbb\to \Zbb$  
such that $\mm_{\cgot}(\mu,k)=\mm(\mu,k)$ for every $(\mu,k)$ in the cone
\begin{equation}
  \label{eq:cone-c}
  {\cone}(\cgot)=\{(\xi,s)\in\tgot^*\times\Rbb^{>0}\,\vert\, \xi\in s\cdot\cgot\}.
\end{equation}

\bigskip

Consider now two adjacent connected regions $\cgot_\pm$ of regular values of 
$\Phi$ separated by an hyperplane $\p\subset\tgot^*$. When $\p$ does
not contain a facet of the polytope $\Phi(M)$, Meinrenken an Sjamaar proved 
also that 
\begin{equation}
  \label{eq:MS-egalite}
  \mm_{\cgot_+}(\mu,k)=\mm_{\cgot_-}(\mu,k)=\mm(\mu,k)
\end{equation}
for every $(\mu,k)\in \overline{\cone(\cgot_+)}\cap\overline{\cone(\cgot_-)}
=\cone(\overline{\cgot_+}\cap\overline{\cgot_-})\subset\cone(\p)$. 
Our main objective is to prove that (\ref{eq:MS-egalite}) extends  
to a ``strip'' containing $\cone(\p)$. Let $\beta\in\Lambda$ be 
the primitive orthogonal vector to the hyperplane 
$\p\subset\tgot^*$ which is pointing out of $\cgot_-$. Then 
\begin{equation}
  \label{eq:delta-n}
  \p=\{\xi\in\tgot^*\,\vert\, 
\frac{\langle\mu,\beta\rangle}{2\pi} = r_\p\}
\end{equation}
for some $r_\p\in\Zbb$, ${\cone}(\p)=
\{(\xi,s)\in\tgot^*\times\Rbb^{\geq 0}\,\vert\,
\frac{\langle\xi,\beta\rangle}{2\pi} - s r_\p=0 \}$ and 
 $\cgot_- \subset$ \break $\{\xi\in\tgot^*\,\vert\,
\frac{\langle\mu,\beta\rangle}{2\pi} < r_\p \}$.

Let $T_{\p}$ be the subtorus of $T$ generated by $\beta$. 
Let $N_\p$ be the normal vector bundle of $M^{T_\p}$ in
$M$. The almost complex structure on  
$M$ induces a complex structure $J$ on the fibers of $N_\p$. 
We have a decomposition 
$N_\p=\sum_s N_\p^s$ where 
$N_\p^s=\{v\in N_\p\,\vert\, \Lcal_\beta v= s\, Jv\ \}$. 
We write  $N_\p=N_\p^{+,\beta}\oplus N_\p^{-,\beta}$ where 
\begin{equation}
  \label{eq:N-beta-plus}
 N_\p^{\pm,\beta}=\sum_{\pm s >0}N_\p^s.
\end{equation}

\begin{defi}
For every connected component $Z\subset M^{T_\p}$ we define $s^\pm_Z\in\Nbb$
respectively as the absolute value of the trace of $\frac{1}{2\pi}\Lcal_\beta$ on 
$N_\p^{\pm,\beta}\vert_Z$. Note that $s^+_Z + s^-_Z$ is larger than half of 
the codimension of $Z$ in $M$.
\end{defi}

We prove in Section \ref{subsec:jump-formula} the following 
 
\begin{theo} \label{theo:m-pm-section}
We have $\mm_{\cgot_+}(\mu,k)=\mm_{\cgot_-}(\mu,k)$ for all 
$(\mu,k)\in\Lambda^*\times\Zbb$ such that 
\begin{equation}\label{eq:encadrement}
  - s^- < \frac{\langle\mu,\beta\rangle}{2\pi}- k\, r_\p < s^+ \ .
\end{equation}
The number $s^-,s^+\in\Nbb$ are defined as follows. We take $s^\pm=\inf_Z s^\pm_Z$ 
where the minimum is taken over the connected components $Z$ of
$M^{T_\p}$ for which $\overline{\cgot_+}\cap\overline{\cgot_-}\subset\Phi(Z)$. 
\end{theo}

Similar results were obtained by Billey-Guillemin-Rassart \cite{BGR03}
in the case where $M$ is a coadjoint orbit of $\SU(n)$, and by 
Szenes-Vergne \cite{Szenes-Vergne02} in the case where $M$ is a
complex vector space. See Sections \ref{subsec:su-n} and
\ref{sec:cas-C-d} where we study these 
two particular cases in details. In Proposition
\ref{prop:inegalite-optimale}, we give also a 
criterium which says when the inequalities in (\ref{eq:encadrement}) 
are optimal. This criterium is fullfilled when there is 
{\em only one} component $Z$ of $M^{T_\p}$ such 
that $\overline{\cgot_+}\cap\overline{\cgot_-}\subset\Phi(Z)$. 
Then (\ref{eq:encadrement}) is {\em optimal} and 
$s^+ + s^-$ is larger than half of the codimension of $Z$ in $M$. 

The following easy Lemma (see Lemma 7.3. of \cite{pep4}) gives some
basic informations about the integer $s_Z^\pm$.

\begin{lem}\label{lem:pep4}
  Let $(M,\Omega,\Phi)$ be a compact Hamiltonian $T$-manifold equipped
  with a $T$-invariant almost complex structure compatible with
  $\Omega$. Consider a non-zero vector $\gamma\in\tgot$
and let $Z$ be a connected component of the fixed point set $M^\gamma$. Let 
$N$ be the normal vector of $Z$ in $M$ and let $N^{-,\gamma}$ be the 
negative polarized normal bundle (see (\ref{eq:N-beta-plus})). 
Then $N^{-,\gamma}=0$ if and only if the function 
$\langle\Phi,\gamma\rangle: M\to\Rbb$ takes its maximal value on $Z$.
\end{lem}

This Lemma insures that $s^\pm\geq 1$ in Theorem
\ref{theo:m-pm-section} when $\p\cap\Phi(M)$ is not a facet of the 
polytope $\Phi(M)$.

Consider the situation where  $\p\cap\Phi(M)$ is a facet 
of the polytope $\Phi(M)$ so that $\cgot_+\cap\Phi(M)=\emptyset$: 
hence $\mm_{\cgot_+}=0$. If we apply Lemma \ref{lem:pep4} to 
$\gamma=\beta$, one gets $N^{-,\beta}=0$ and so $s^-=0$. 
In this situation we get 

\begin{coro}Let $\cgot_-$ be a connected component of regular values
  of $\Phi$ bording a facet $\Phi(M)\cap\p$ of the polytope
  $\Phi(M)$. Let $\beta\in\Lambda$ be the unimodular orthogonal vector 
  to the hyperplane $\p\subset\tgot^*$ which is pointing out of
  $\cgot_-$. Here $Z=\Phi^{-1}(\p)$ is a connected 
component of the fixed point set $M^{T_\p}$. We have 
$\mm_{\cgot_-}(\mu,k)=0$ for all $(\mu,k)\in\Lambda^*\times\Zbb$ such that 
\begin{equation}
  \label{eq:encadrement-bord}
 0< \frac{\langle\mu,\beta\rangle}{2\pi}- k r_\p < s^+_Z .
 \end{equation} 
Here $s^+_Z\in \Nbb$ is the value of the trace of 
$\frac{1}{2\pi}\Lcal_\beta$ on the normal bundle of $Z$ in $M$, 
and then is larger
than half of the codimension of $Z$ in $M$. Moreover the 
inequalities (\ref{eq:encadrement-bord}) are optimal.
\end{coro}

We first review some of the results of \cite{pep4}.

\subsection{Elliptic and transversally elliptic symbols}

We work in the setting of a compact manifold $M$ equipped with a smooth  
action of a torus $T$.

Let $p:\T M\to M$ be the projection, and let $(-,-)_M$ be a 
$T$-invariant Riemannian metric.
If $E^{0},E^{1}$ are $T$-equivariant vector bundles over $M$, a 
$T$-equivariant morphism $\sigma \in \Gamma(\T M,\hom(p^{*}E^{0},p^{*}E^{1}))$  
is called a {\em symbol}. The subset of all $(m,v)\in \T M$ where 
$\sigma(m,v): E^{0}_{m}\to E^{1}_{m}$ is not invertible is called 
the {\em characteristic set} of $\sigma$, and is denoted by $\Char(\sigma)$. 

Let $\T_{T}M$ be the following subset of $\T M$ :
$$
   \T_{T}M\ = \left\{(m,v)\in \T M,\ (v,X_{M}(m))_{_{M}}=0 \quad {\rm for\ all}\ 
   X\in\kgot \right\} .
$$

A symbol $\sigma$ is {\em elliptic} if $\sigma$ is 
invertible outside a compact subset of $\T M$ ($\Char(\sigma)$ is
compact), and is {\em transversally elliptic} if the restriction of $\sigma$ 
to $\T_{T}M$ is invertible outside a compact subset  of $\T_{T}M$ 
($\Char(\sigma)\cap \T_{T}M$ is compact). An elliptic symbol $\sigma$ defines 
an element in the equivariant $K$-theory of $\T M$ with compact 
support, which is denoted by $\K_{T}(\T M)$, and the 
index of $\sigma$ is a virtual finite dimensional representation of $T$
\cite{Atiyah-Segal68,Atiyah-Singer-1,Atiyah-Singer-2,Atiyah-Singer-3}.

A {\em transversally elliptic} symbol $\sigma$ defines an element of 
$\K_{T}(\T_{T}M)$, and the index of $\sigma$ is defined as a trace class virtual 
representation of $T$ (see \cite{Atiyah74} for the analytic index and 
\cite{B-V.inventiones.96.1,B-V.inventiones.96.2} for the cohomological one). 
Remark that any elliptic symbol of $\T M$ is transversally elliptic, hence 
we have a restriction map $\K_{T}(\T M)\to \K_{T}(\T_{T}M)$, and 
a commutative diagram
\begin{equation}\label{indice.generalise}
\xymatrix{
\K_{T}(\T M)\ar[r]\ar[d]_{\indice_{M}^T} & 
\K_{T}(\T_{T}M)\ar[d]^{\indice_{M}^T}\\
R(T)\ar[r] &  R^{-\infty}(T)\ .
   }
\end{equation} 

\medskip

Using the {\em excision property}, one can easily show that the index map 
$\indice_{\Ucal}^{T}:\K_{T}(\T_{T}\Ucal)\to R^{-\infty}(T)$
is still defined when $\Ucal$ is a $T$-invariant relatively compact 
open subset of a $T$-manifold (see \cite{pep4}[section 3.1]).


\subsection{Localization of the Riemann-Roch character}\label{ssection.thom}

We suppose now that the compact $T$-manifold $M$ is equipped with a 
$T$-invariant almost complex structure $J$. Let us recall the 
definitions of the Thom symbol $\Thom(M,J)$ and of the 
Riemann-Roch character \cite{pep4}.

Consider a $T$-invariant Riemannian metric $q$ on $M$ such that 
$J$ is orthogonal relatively to $q$, and let 
$h$ be the  Hermitian structure on  $\T M$ defined by : 
$h(v,w)=q(v,w) -\imath q(Jv,w)$ for 
$v,w\in \T M$. The symbol 
$$
\Thom_{_{T}}(M,J)\in 
\Gamma\left(M,\hom(p^{*}(\wedge_{\Cbb}^{even} \T M),\,p^{*}
(\wedge_{\Cbb}^{odd} \T M))\right)
$$  
at $(m,v)\in \T M$ is equal to the Clifford map
\begin{equation}\label{eq.thom.complex}
 \Clif_{m}(v)\ :\ \wedge_{\Cbb}^{even} \T_m M
\longrightarrow \wedge_{\Cbb}^{odd} \T_m M,
\end{equation}
where $\Clif_{m}(v).w= v\wedge w - c_{h}(v).w$ for $w\in 
\wedge_{\Cbb}^{\bullet} \T_{x}M$. Here $c_{h}(v):\wedge_{\Cbb}^{\bullet} 
\T_{m}M\to\wedge^{\bullet -1} \T_{m}M$ denotes the 
contraction map relative to $h$. Since the map $\Clif_{m}(v)$ is invertible for all
$v\neq 0$, the symbol $\Thom_{_{T}}(M,J)$ is 
elliptic. 
 
The Riemann-Roch character $RR(M,-):\K_{T}(M)\to R(T)$ is 
defined by the following relation
\begin{equation}\label{def.RR}
RR(M,E)=\indice^{T}_{M}\left(\Thom_{_{T}}(M,J)\otimes p^{*}E\right)\ .
\end{equation}
The important point is that for any $T$-vector bundle $E$, 
$\Thom_{_{T}}(M,J)\otimes p^{*}E$ corresponds to the {\em principal symbol} of the  
twisted $\spinc$ Dirac operator $\Dcal_{E}$ \cite{Duistermaat96}, hence 
$RR(M,E)\in R(T)$ is also defined as the (analytical) index of the elliptic 
operator $\Dcal_E$.

\medskip

Consider now the case of a {\em compact} Hamiltonian $T$-manifold 
$(M,\omega,\Phi)$. Here $J$ is a $T$-invariant almost comlex structure
compatible with $\Omega$: $(v,w)\mapsto\Omega(v,Jw)$ defines a Riemannian 
metric on $M$. Like in Section \ref{sec:hd}, we make the choice of 
a scalar product $(-,-)$ on $\tgot^*$ 
(which induces an identification $\tgot^*\simeq\tgot$) 
and we consider for any $\xi\in\tgot^*$ 
the function $\frac{-1}{2}\parallel\Phi-\xi\parallel^2: M\to \Rbb$ and its 
Hamiltonian vectors field $\Hcal-\xi_M$.

\begin{defi}\label{def.thom.loc}
For any $\xi\in\tgot^*$ and any $T$-invariant open subset $\Ucal\subset M$ 
we define the symbol $\Thom_{\xi}(\Ucal)$ by the relation
$$
\Thom_{\xi}(\Ucal)(m,v):=\Thom(M,J)(m,v-(\Hcal-\xi_M)(m))\quad (m,v)\in\T \Ucal
$$
\end{defi}

The characteristic set of $\Thom_{\xi}(\Ucal)$ corresponds to
$\{(m,v)\in \T \Ucal,\ v=(\Hcal-\xi_M)(m)\}$, the graph of the vector 
field $\Hcal-\xi_M$ over $\Ucal$. Since $\Hcal-\xi_M$ belongs to the set of 
tangent vectors to the $T$-orbits, we have
\begin{eqnarray*}
\Char\left(\Thom_{\xi}(\Ucal)\right)\cap \T_{T}\Ucal 
&=&\{(m,0)\in \T \Ucal\,\vert\, (\Hcal-\xi_M)(m)=0 \}\\
&\cong& \{m\in \Ucal,\ d\parallel\Phi-\xi\parallel^2_m=0 \} \ .
\end{eqnarray*}
Therefore the symbol $\Thom_{\xi}(\Ucal)$ is transversally elliptic 
if and only if  
\begin{equation}
  \label{eq:U-xi}
  \crpt(\parallel\Phi-\xi\parallel^2)\cup\partial\Ucal=\emptyset.
\end{equation}
Here $\crpt(\parallel\Phi-\xi\parallel^2)$ denotes the set of critical points 
of the function $\parallel\Phi-\xi\parallel^2$. When (\ref{eq:U-xi}) holds
we say that the couple $(\Ucal,\xi)$ is {\em good}.

\begin{defi}\label{def:RR-U-xi}
Let $(\Ucal,\xi)$ be a good couple. For any $T$-vector bundle 
$E\to M$, the tensor product $\Thom_{\xi}(\Ucal)\otimes p^{*}E$ 
belongs to $\K_{T}(\T_{T}\Ucal)$ and we denote by  
    $$
    RR_{\Ucal}^{\xi}(M,E)\in R^{-\infty}(T)
    $$ 
its index.
\end{defi}

\begin{prop}\label{prop:U-xi}
Let $(\Ucal,\xi)$ be a good couple. 

$a)$\ If $\Ucal$ possess two $T$-invariant open subsets
$\Ucal^1,\Ucal^2$ such that 
$\overline{\Ucal^1}\cap \overline{\Ucal^2}\cap$ \break 
$\crpt(\parallel\Phi-\xi\parallel^2)=\emptyset$ and 
$(\Ucal^1\cup\Ucal^2)\cap\crpt(\parallel\Phi-\xi\parallel^2)=
\Ucal\cap\crpt(\parallel\Phi-\xi\parallel^2)$, then the couples $(\Ucal^1,\xi)$ 
and $(\Ucal^2,\xi)$ are good and 
$$
RR_{\Ucal}^{\xi}(M,-)=RR_{\Ucal^1}^{\xi}(M,-)+RR_{\Ucal^2}^{\xi}(M,-).
$$
In particular $RR_{\Ucal}^{\xi}(M,-)=RR_{\Ucal^1}^{\xi}(M,-)$ 
if $\Ucal^1$ is an open subset 
of $\Ucal$ such that $\Ucal^1\cap\crpt(\parallel\Phi-\xi\parallel^2)=
\Ucal\cap\crpt(\parallel\Phi-\xi\parallel^2)$.

$b)$\ If $\xi'\in\tgot^*$ is close enough to $\xi$, then $(\Ucal,\xi')$ 
is  good and 
$$
RR_{\Ucal}^{\xi}(M,-)=RR_{\Ucal}^{\xi'}(M,-).
$$
\end{prop}

{\em Proof.} The part $a)$ is a direct consequence of the excision property 
(see Proposition 4.1. in \cite{pep4}). Consider now the scalar product 
$(\Hcal-\xi^s_M,\Hcal-\xi_M)$ where $\xi^s=s\xi'+(1-s)\xi,\ s\in[0,1]$ 
and $(-,-)$ is a $T$-invariant Riemannian metric on $M$. We have 
$(\Hcal-\xi^s_M,\Hcal-\xi_M)=\|\Hcal-\xi_M\|^2+ s((\xi-\xi')_M,\Hcal-\xi_M)$
and then the following inequality holds on $M$
\begin{equation}
  \label{eq:maj-prime}
  (\Hcal-\xi^s_M,\Hcal-\xi_M)\geq \|\Hcal-\xi_M\|^2 
\Big( \|\Hcal-\xi_M\|-s\|(\xi-\xi')_M\|\Big). 
\end{equation}
Since $\partial\Ucal$ is compact we have the following inequalities on 
it: $\|\Hcal-\xi_M\|\geq c_1>0$ and $\|a_M\|\leq c_2\|a\|$ for
any $a\in\tgot$. So (\ref{eq:maj-prime}) implies the following 
inequality on $\partial\Ucal$:
$$
(\Hcal-\xi^s_M,\Hcal-\xi_M)\geq c_1( c_1-s\|\xi-\xi'\|)\quad {\rm for}\quad s\in[0,1]. 
$$
So if $\xi'$ is close enough to $\xi$, we have $\|\Hcal-\xi^s_M\|\geq
c_3>0$ on $\partial\Ucal$ for any $s\in[0,1]$. We have first prove 
that the couple $(\Ucal,\xi^s)$ 
is good  for any $s\in[0,1]$. We see then that the familly of transversally
elliptic symbols $\Thom_{\xi^s}(\Ucal),\, s\in[0,1]$ defines an homotopy  
between $\Thom_{\xi}(\Ucal)$ and 
$\Thom_{\xi'}(\Ucal)$. Hence $\Thom_{\xi}(\Ucal)=\Thom_{\xi'}(\Ucal)$ in 
$\K_{T}(\T_{T}\Ucal)$. $\Box$

\medskip

Part $a)$ of Proposition \ref{prop:U-xi} tells us that $RR^\xi_\Ucal(M,-)$ depends 
closely \break of the intersection $\Ucal\cap\crpt(\parallel\Phi-\xi\parallel^2)$. 
In particular $RR^\xi_\Ucal(M,-)=0$ when \break 
$\Ucal\cap\crpt(\parallel\Phi-\xi\parallel^2)=\emptyset$.  
Recall that   
\begin{equation}  \label{eq:decom-kirwan}
 \crpt(\parallel\Phi-\xi\parallel^2)=
 \bigcup_{\gamma\in\Bcal_\xi}M^\gamma\cap\Phi^{-1}(\gamma+\xi)
\end{equation}
where $\Bcal_\xi\subset\tgot^*$ is a finite set \cite{Kirwan84}. 

\begin{defi}\label{def:RR-xi}
 For any $\xi\in\tgot^*$ and $\gamma\in\Bcal_\xi$, we denote 
simply by 
$$
RR^\xi_\gamma(M,-):\K_{T}(M)\to R^{-\infty}(T)
$$ 
the map  $RR_{\Ucal}^{\xi}(M,-)$, where  $\Ucal$ is a $T$-invariant 
open neighborhood of $M^\gamma\cap\Phi^{-1}(\gamma+\xi)$ such that 
$\crpt(\parallel\Phi-\xi\parallel^2)\cap\overline{\Ucal}=
M^\gamma\cap\Phi^{-1}(\gamma+\xi)$.  
\end{defi}

Part $a)$ of Proposition \ref{prop:U-xi} insures that the maps 
$RR_\gamma^\xi(M,-)$ are well defined, and for any good couple 
$(\Ucal,\xi)$ we have
\begin{equation}
  \label{eq:RR-U-gamma}
   RR^\xi_\Ucal(M,-)=\sum_{\gamma\in\Bcal_\xi\cap\Phi(\Ucal)} RR_{\gamma}^{\xi}(M,-).
\end{equation}
If one takes $\Ucal=M$, we have 
$RR^\xi_\Ucal(M,-)=RR(M,-)=\sum_{\gamma\in\Bcal_\xi} RR_{\gamma}^{\xi}(M,-)$ 
(see \cite{pep4}[Section 4]).

\subsection{Periodic polynomial behaviour of the multiplicities}

We suppose here that the Hamiltonian $T$-manifold $(M,\Omega,\Phi)$ 
is prequantized by a $T$-complex line bundle $L$ satisfying (\ref{eq:kostant}) 
for a suitable invariant connection. In this section we will 
characterize the periodic polynomial behaviour of the multiplicities 
$\mm(\mu,k)$ with the help of the 
localized Riemann-Roch character $RR^\xi_0(M,-)$.

Let us introduce some vocabulary. We say that two generalized characters 
$\chi^\pm=\sum_{\mu\in\Lambda^*}a^\pm_\mu\,\Cbb_\mu$ {\em coincide} on
a region $D\subset\tgot^*$, if $a^+_\mu=a^{-}_\mu$ for every 
$\mu \in D\cap\lambda^*$. A generalized character 
$\chi=\sum_{\mu}a_\mu\, \Cbb_\mu$ is {\em supported} on a region 
$D\subset\tgot^*$ if $a_\mu=0$ for $\mu\notin D$. A weight 
$\mu\in\Lambda^*$ {\em occurs} in $\chi=\sum_{\mu}a_\mu\, \Cbb_\mu$ 
if $a_\mu\neq 0$.

For $\xi\in\tgot^*$, we define $r_\xi>0$ as the smallest 
non-zero critical value of the function $\parallel\Phi-\xi\parallel$, 
and we denote by $B(\xi,r_\xi)$ the open ball of center $\xi$ 
and radius $r_\xi$.

\begin{theo}[\cite{pep4}]\label{th-RR-xi}
  For any $\xi\in\tgot^*$, the generalized character $RR_0^\xi(M,L^{\otimes k})$ 
coincides with $RR(M,L^{\otimes k})$ on the open ball $k\cdot B(\xi,r_\xi)$.
\end{theo}

The arguments of \cite{pep4} for the proof of this Theorem will be needed another 
time, so we recall them. Let $\xi\in\tgot^*$. We start with the decomposition
\begin{equation}
  \label{eq:loc-gene}
 RR(M,L^{\otimes k})=\sum_{\gamma\in\Bcal_\xi} RR_{\gamma}^{\xi}(M,L^{\otimes k}).
\end{equation}
We recall now, for a non-zero $\gamma\in\Bcal_\xi$, the localization of 
the map $RR_{\gamma}^{\xi}$ on the fixed point set $M^\gamma$ \cite{pep4}

\medskip

Let $N$ be the normal bundle of $M^\gamma$ in $M$. The almost complex structure on 
$M$ induces an almost complex struture on $M^\gamma$ and a complex structure on the 
bundles $N$ and $N_\Cbb:=N\otimes \Cbb$. Following 
(\ref{eq:N-beta-plus}) we define the $\gamma$-polarized complex vector bundles
$N^{+,\gamma}$ and $(N_\Cbb)^{+,\gamma}$. 

The manifold $M^\gamma$ is a symplectic submanifold of $M$ equipped with an
induced Hamiltonian action of $T$: its moment map is  the 
restriction of $\Phi$ on $M^\gamma$. Following Definition \ref{def:RR-xi}, we have 
on $M^\gamma$ a localized Riemann-Roch character $RR_\gamma^\xi(M^\gamma,-)$.
On $M^\gamma$, the Hamiltonian vectors fields of the functions 
$\parallel\Phi-\xi\parallel^2$ and $\parallel\Phi-(\xi+\gamma)\parallel^2$ 
coincide, hence 
\begin{equation}
  \label{eq:xi-gamma}
  RR_\gamma^\xi(M^\gamma,-)=RR_0^{\xi+\gamma}(M^\gamma,-).
\end{equation}

We prove in \cite{pep4}[Theorem 5.8.] that 
\begin{equation}
  \label{eq:loc-RR-beta}
  RR_\gamma^\xi(M,E)=\sum_{k\in\Nbb}(-1)^l
RR_\gamma^\xi(M^\gamma,E\vert_{M^\gamma}\otimes \det(N^{+,\gamma})\otimes
S^k(N_\Cbb^{+,\gamma}))
\end{equation}
for every $T$-vector bundle $E$. Here $l$ is the  locally constant
fonction on $M^\gamma$ equal 
to the complex rank of $N^{+,\gamma}$.

\begin{prop}\cite{pep4}[Section 5]\label{def:lambda-N}
Let $\overline{N}$ be the $T$-vector bundle $N$ with the opposite 
complex structure on the fibers. The sum 
$(-1)^l\sum_{k\in\Nbb} \det(N^{+,\gamma})\otimes S^k(N_\Cbb^{+,\gamma})$ 
is an inverse of $\wedge_{\Cbb}^{\bullet}\overline{N}$  that we denote
$\left[\wedge_{\Cbb}^{\bullet}\overline{N}\right]^{-1}_{\gamma}$.
\end{prop}

If we use the notations of Proposition \ref{def:lambda-N} and 
(\ref{eq:xi-gamma}), the localization (\ref{eq:loc-RR-beta}) 
can be rewritten as
\begin{equation}
  \label{eq:loc-RR-beta-bis}
  RR_\gamma^\xi(M,E)=
RR_0^{\xi+\gamma}\left(M^\gamma,E\vert_{M^\gamma}\otimes
\left[\wedge_{\Cbb}^{\bullet}\overline{N}\right]^{-1}_{\gamma}\right).
\end{equation}

Let $i:T_\gamma\croc T$ be the inclusion of the subtorus generated by $\gamma$. 
For a $T$-vector bundle $F\to M^\gamma$, it is easy to show that 
a weight $\mu\in\Lambda^*$ occurs in $RR_\gamma^\xi(M^\gamma,F)$ only if 
$i^*(\mu)$ occurs as a weight for the $T_\gamma$-action on the fibers of $F$
(see Lemma 9.4. in \cite{pep4}). 
Since the $T_\gamma$ weights on the bundles $N_\Cbb^{+,\gamma}$ and 
$N^{+,\gamma}$ are polarized by $\gamma$, the localization 
(\ref{eq:loc-RR-beta}) gives the following 

\begin{prop}
 For a non-zero $\gamma\in\Bcal_\xi$, the generalized character 
$RR_\gamma^\xi(M,L^{\otimes k})$ is supported on  the
half space $\{a\in\tgot^*\,\vert\, (\gamma, a-k(\xi+\gamma))\geq 0\}$. 
\end{prop}

Since the condition $(\gamma, a-k(\xi+\gamma))\geq 0$ implies that 
$\parallel a-k\xi\parallel\geq k \parallel\gamma\parallel\geq k r_\xi$, the last 
proposition shows that every weights of the open ball $k\cdot B(\xi,r_\xi)$ 
does not occurs in $RR_\gamma^\xi(M,L^{\otimes k})$. This last remark together 
with (\ref{eq:loc-gene}) prove Theorem \ref{th-RR-xi}. 

\medskip

For the localized Riemann-Roch character $RR_0^\xi(M,-)$ we have 
the following Lemma which is very similar to Lemma \ref{lem:P-xi-0}. 

\begin{lem}\label{lem:RR-xi-prime} Let $\cgot\subset\tgot^*$ be a 
connected component of regular values of $\Phi$.
  For every $\xi,\xi'\in\cgot$, we have $RR_0^\xi(M,-)=RR_0^{\xi'}(M,-)$.
\end{lem}

{\em Proof.} We have to show that the map $\xi\mapsto RR_0^\xi(M,-)$ is 
locally constant on $\cgot$. Let $\xi\in\cgot$ and take an open neigborhood
$\Ucal$ of $\Phi^{-1}(\xi)$ small enough such that the stabilizer 
$T_m=\{t\in T\,\vert\, t\cdot m=m\}$ is finite for every $m\in\overline{\Ucal}$. 
We see then that $\overline{\Ucal}\cap\crpt(\parallel \Phi-\xi'\parallel^2)=
\Phi^{-1}(\xi')$ and $\partial\Ucal\cap\crpt(\parallel \Phi-\xi'\parallel^2)=
\emptyset$ if $\xi'$ is close enough to $\xi$: hence  
$RR_0^{\xi'}(M,-)=RR_\Ucal^{\xi'}(M,-)$ for $\xi'$ close enough to $\xi$. 
Part $b)$ of Proposition \ref{prop:U-xi} finishes the proof. $\Box$

\medskip

When $\xi$ is a regular value of $\Phi$, the localized Riemann-Roch character 
$RR_0^\xi(M,-)$ as been computed in \cite{pep4} as follows. Let 
$RR(\Mcal_\xi,-)$ be the Riemann-Roch map defined on the orbifold
$\Mcal_\xi=\Phi^{-1}(\xi)/T$ by means of an almost complex structure compatible 
with the induced symplectic structure. For every $T$-vector bundle
$E\to M$ we define the following familly of orbifold vector bundles  
over $\Mcal_\xi$:
\begin{equation} \label{eq:ecal-mu}
  \Ecal_{\xi,\mu}:=\Big(E\vert_{\Phi^{-1}(\xi)}\otimes\Cbb_{-\mu}\Big)/T,
  \quad \mu\in\Lambda^*.
\end{equation}
For every $T$-vector bundle $E$ on $M$, we proved in \cite{pep4}[Section 6.2.] the 
following equality in $R^{-\infty}(T)$
\begin{equation} \label{eq:RR-xi-O}
 RR_0^\xi(M,E)=\sum_{\mu\in \Lambda^*}RR(\Mcal_\xi,\Ecal _{\xi,\mu})\,\Cbb_\mu.
\end{equation}
This decomposition was first obtained by Vergne \cite{Vergne96} when $T$ is the
circle group and when $M$ is Spin. The number 
$RR(\Mcal_\xi, \Ecal_{\xi,\mu})\in \Zbb$ is then equal to the $T$-invariant 
part of the index $RR_0^\xi(M,E)\otimes\Cbb_{-\mu}$. 

\begin{rem}\label{lem:discrete}
  Let $t\to t^\lambda$ be a character of $T$. Suppose that a subgroup 
$H\subset T$ acts trivially on $M$ and with the character $t\in H\to t^\lambda$ on the 
the fibers of the $T$-vector bundle $E$. Then $H$ acts with the character 
$t\in H\to t^{\lambda-\mu}$ on $RR_0^\xi(M,E)\otimes\Cbb_{-\mu}$, and then 
$RR(\Mcal_\xi, \Ecal_{\xi,\mu})\neq 0$ only if $t^{\lambda-\mu}=1$ for every $t\in H$. 
So the sum in (\ref{eq:RR-xi-O}) can be restricted to $\lambda + \Lambda^*_H$, 
where $\Lambda^*_H$ is the sub-lattice of $\Lambda^*$ formed by the 
element $\alpha\in\Lambda^*$ satisfying $t^\alpha=1,\ \forall\ t\in H$.

This remark applies also 
on the usual character $RR(M,E)=\sum_{\mu\in\Lambda^*}m_\mu\Cbb_\mu$.
The mutiplicity $m_\mu\in\Zbb$ is equal to the (virtual) dimension of 
the $T$-invariant part of $RR(M,E)\otimes\Cbb_{-\mu}$. With the same 
hypothesis than above we see that $m_\mu\neq 0$ only if 
$\mu\in \lambda + \Lambda^*_H$.
\end{rem}

Let $\Gamma_M$ be the generic stabilizer for the action of $T$ on
$M$. Consider a weight $\alpha_{o}$ such that $\Gamma_M$ acts 
on the fibers of $L$ with the character $t\mapsto t^{\alpha_{o}}$. 
We define the sublattice $\Xi(M,L)\subset \Lambda^*\times \Zbb$ by 
\begin{equation}
  \label{eq:xi-lattice}
 \Xi(M,L):=\{(\mu,k)\in \Lambda^*\times \Zbb\ \vert\ k\alpha_{o}-\mu\in
\Lambda^*_{\Gamma_M}\}.
\end{equation}

We know then that  $\mm(\mu,k)=0$ if $(\mu,k)\notin \Xi(M,L)$.

\begin{prop}\label{prop:m-poly}
  Let $\cgot$ be a connected component of regular values of $\Phi$ and 
let $\cone(\cgot)$ be the corresponding cone in $\tgot^*\times\Rbb^{> 0}$ 
(see (\ref{eq:cone-c})). Let $\xi\in\cgot$. For any $(\mu,k)\in\cone(\cgot)
\cap\Xi(M,L)$ we have $\mm(\mu,k)=RR(\Mcal_\xi, \Lcal^k_{\xi,\mu})$ where 
\begin{equation}
  \label{eq:line-k-xi-mu}
 \Lcal^k_{\xi,\mu}=
(L^{\otimes k}\vert_{\Phi^{-1}(\xi)}\otimes\Cbb_{-\mu})/T. 
\end{equation} 
\end{prop}

\medskip

{\em Proof.} Let $(\mu,k)\in\cone(\cgot)$ and let $\xi'=\frac{\mu}{k}\in
\cgot$. We known from Theorem \ref{th-RR-xi} that 
the generalized character $RR_0^{\xi'}(M,L^{\otimes k})$ 
coincides with $RR(M,L^{\otimes k})$ on the open ball $k\cdot B(\xi',r_{\xi'})=
B(\mu,k r_{\xi'})$. So $\mm(\mu,k)$ is equal to the $\mu$-multiplicity in 
$RR_0^{\xi'}(M,L^{\otimes k})$. Take now any $\xi\in\cgot$. We know after
Lemma \ref{lem:RR-xi-prime} that $RR_0^\xi(M,-)=RR_0^{\xi'}(M,-)$ and 
(\ref{eq:RR-xi-O}) shows that the $\mu$-multiplicity in 
$RR_0^{\xi}(M,L^{\otimes k})$ is equal to 
$RR(\Mcal_\xi, \Lcal^k_{\xi,\mu})$. 
$\Box$

\begin{defi}
Take $\xi\in\cgot$. The map $\mm_{\cgot}:\Lambda^*\times\Zbb\to \Zbb$ is 
defined by the equation 
  \begin{equation}
    \label{eq:def:m-c}
    \mm_{\cgot}(\mu,k)=RR(\Mcal_\xi, \Lcal^k_{\xi,\mu}),
  \end{equation}
where $\Lcal^k_{\xi,\mu}$ is the orbifold line bundle defined by 
(\ref{eq:line-k-xi-mu}). In other words, the map $\mm_{\cgot}$ is defined by the 
following equality in $R^{-\infty}(T)$
$$
\sum_{\mu\in\Lambda^*}\mm_{\cgot}(\mu,k)\ \Cbb_\mu=RR_0^\xi(M,L^{\otimes k}).
$$
for all $k\in\Zbb$. After remark \ref{lem:discrete}, we know that 
$\mm_{\cgot}$ is supported on the 
sublattice $\Xi(M,L)$ defined in (\ref{eq:xi-lattice}).
\end{defi}


We will now exploit the Riemann-Roch for orbifold due to Atiyah-Kawasaki 
\cite{Atiyah74,Kawasaki79} to 
show that the map $\mm_{\cgot}$ is a periodic polynomial.


\subsection{Riemann-Roch theorem on orbifolds}\label{subsec:RR-theo}
First we recall how is defined the Riemann-Roch character 
$RR(\Mcal_\xi,\Ecal_\xi)$ when 
$\xi$ is a regular value of $\Phi$, and $\Ecal_\xi= E_{\vert\Phi^{-1}(\xi)}/T$
is the reduction of a complex $T$-vector bundle $E$ over $M$. The number 
$RR(\Mcal_\xi,\Ecal_\xi)\in\Zbb$ is defined has the $T$-invariant part of 
the index of a transversally 
elliptic operator $D_E$ on $\Phi^{-1}(\xi)$. Since the index of $D_E$ depend 
only of the class of its symbol $\sigma(D_E)$ in $K_T(\T_T\Phi^{-1}(\xi))$, 
it is enough to define the transversally elliptic symbol $\sigma(D_E)$. Since 
the action of $T$ on $\Phi^{-1}(\xi)$ is locally free, $V:=\T_T\Phi^{-1}(\xi)$ is 
a vector bundle. It carries a canonical symplectic structure on the fibers and 
we choose any compatible complex structure making $V$ into a Hermitian vector bundle. 
At $(m,v)\in\T\Phi^{-1}(\xi)$, the map $\sigma(D_E)(m,v)$ is the Clifford action 
$$
\Clif_{m}(v_1)\otimes Id_{E_m}\ :\ (\wedge_{\Cbb}^{even} V_m)\otimes E_m
\longrightarrow (\wedge_{\Cbb}^{odd} V_m)\otimes E_m.
$$
where $v_1\in V_m$ is the $V$-component of the vector
$v\in\T_m\Phi^{-1}(\xi)$. We explain
now the formula of Atiyah-Kawasaki for $RR(\Mcal_\xi,\Ecal_\xi)$ 
when $\xi\in\Phi(M)$ is a regular 
value of $\Phi$ \cite{Atiyah74,Kawasaki79}.

\medskip

Let $\Fbb$ be the collection of the finite subgroup of $T$ which are stabilizer 
of points in $M$. Consider the orbit type stratification of $\Phi^{-1}(\xi)$ 
and denote by $\Sbb_\xi$ the set of its orbit type strata. Each statum 
$S$ is a connected component of the smooth submanifold
\begin{equation}
  \label{eq:stab-H}
  \Phi^{-1}(\xi)_{H_S}:= \{m\in\Phi^{-1}(\xi)\,\vert\, \stab_T(m)=H_S\}.
\end{equation}
for a unique $H_S\in\Fbb$. The orbifold 
$\Mcal_\xi$ decomposes as a disjoint union $\cup_{S\in\Sbb_\xi}S/T$ of 
smooth components, and each quotient $\overline{S}/T$ is a  suborbifold 
of $\Mcal_\xi$. The generic stabilizer $\Gamma_M$ of $T$ on $M$ is 
also the generic stabilizer of $T$ on the fiber\footnote{Since a neighborhood of 
$\Phi^{-1}(\xi)$ in $M$ is $T$-equivariantly diffeomorphic to 
$\Phi^{-1}(\xi)\times\tgot^*$.} 
$\Phi^{-1}(\xi)$, and is associated to an open and dense stratum $S_{max}$.

Suppose that $E\to M$ is an Hermitian $T$-vector bundle. On each suborbifold 
$\overline{S}/T$, we get the orbifold complex vector
bundle 
\begin{equation}
  \label{eq:bundle-E-H}
  \Ecal_S:=E_{\vert\overline{S}}/T.
\end{equation} 
We define twisted characteristic classes $\chern^-(\Ecal_S)$ and 
$D^-(\Ecal_S)$ by 
\begin{equation}
  \label{eq:chern}
  \chern^\gamma(\Ecal_S):=\tr
\left(\gamma^{\Ecal_S}\cdot e^{\frac{i}{2\pi}R(\Ecal_S)}\right),\quad \gamma\in H, 
\end{equation}
and 
\begin{equation}
  \label{eq:D-class}
  D^\gamma(\Ecal_S):=\det 
\left(1-(\gamma^{\Ecal_S})^{-1}\cdot e^{-\frac{i}{2\pi}R(\Ecal_S)}\right),
\quad \gamma\in H.
\end{equation}
Here $R(\Ecal_S)\in\Acal^2(\overline{S}/T,\End(\Ecal_S))$ 
is the curvature of an horizontal 
Hermitian connection on $E_{\vert\overline{S}}$, and $\gamma\mapsto
\gamma^{\Ecal_S}$ is the linear action of $H_S$ 
on the fibers of $E_{\vert\overline{S}}$. 

Let $N_S$ be the normal bundle of $\overline{S}$ in $\Phi^{-1}(\xi)$. The 
symplectic struture on $M$ induces a symplectic form $\Omega_S$ on each 
suborbifold $\overline{S}/T$, and a symplectic structure on the fibers 
of the bundle $N_S$. Choose a compatible almost complex structure on 
$\overline{S}/T$, and a compatible complex structure on the fibers 
of $N_S$ making the tangent bundle of $\overline{S}/T$ and 
$\Ncal_S:=N_S/T$ into Hermitian vector bundle. Consider a Hermitian connexion on 
$\T(\overline{S}/T)$, with curvature $R(\overline{S}/T)$, and let 
\begin{equation}
  \label{eq:todd-form}
  \todd(\overline{S}/T)=\det\left(\frac{(i/2\pi)R(\overline{S}/T)}
{1-e^{-(i/2\pi)R(\overline{S}/T)}}\right)
\end{equation}
be the corresponding Todd forms. Like in (\ref{eq:D-class}), we associate to 
the complex orbifold vector bundle $\Ncal_S$, the twisted form 
$D^-(\Ncal_S)$ which is a map form $H_S$ to $\Acal^{even}(\overline{S}/T)$.
The $0$-degree part of $D^\gamma(\Ncal_S)$ is equal to 
$\det(1-(\gamma^{\Ncal_S})^{-1})$, hence $D^\gamma(\Ncal_S)$ is invertible in 
$\Acal^{even}(\overline{S}/T)$ when $\gamma$ belongs to 
\begin{equation}
  \label{eq:H-S-o}
  H_S^o=\{\gamma\in H_S\ \vert\ \det(1-(\gamma^{\Ncal_S})^{-1})\neq 0\}.
\end{equation}
Note that $H_S^o$ corresponds to the set of $\gamma\in H_S$ for which 
$\overline{S}$ is a connected component of $(\Phi^{-1}(\xi))^\gamma$.

\begin{theo}[Atiyah-Kawasaki] \label{th:RR-AK}
The number $RR(\Mcal_\xi,\Ecal_\xi)\in\Zbb$ is given by the 
formula
\begin{equation}
  \label{eq:RR-AK}
 RR(\Mcal_\xi,\Ecal_\xi) = \sum_{S\in\Sbb_\xi}\frac{1}
{\vert H_S\vert}\sum_{\gamma\in H_S^o}\int_{\overline{S}/T}
\frac{\todd(\overline{S}/T)\chern^\gamma(\Ecal_S)}
{D^\gamma(\Ncal_S)}.
\end{equation}
\end{theo}

 We exploit now Theorem \ref{th:RR-AK} to show that the map 
$\mm_{\cgot}:\Lambda^*\times\Zbb\to \Zbb$ which is defined by (\ref{eq:def:m-c}) 
is periodic polynomial. We need the classical computation of the first Chern class 
of the line bundle 
\begin{equation}
  \label{eq:L-S-mu}
  \Lcal^k_{S,\mu}=(L^{\otimes k}\otimes\Cbb_{-\mu})_{\vert\overline{S}}/T.
\end{equation}
The curvature form $\omega_\xi\in H^2(\Mcal_\xi)\otimes\tgot$ of the 
principal $T$-bundle $\Phi^{-1}(\xi)\to\Mcal_\xi$ restricts to a curvature 
form $\omega_S\in H^2(\overline{S}/T)\otimes\tgot$ on each strata.

\begin{lem}\label{lem:c-1-L}
The first Chern class of the line bundle $\Lcal^k_{S,\mu}$ is given by 
$$
c_1(\Lcal^k_{S,\mu})=\frac{1}{2\pi}
\Big(k\Omega_S -\langle k\xi-\mu,\omega_S\rangle\Big).
$$
\end{lem}

For a strata $S$, we consider $\alpha_S\in\Lambda^*$ such that 
$\gamma\in H_S\mapsto \gamma^{\alpha_S}$ corresponds to the action 
of $H_S$ on the fibers of $L_{\vert\overline{S}}$.
Finally we have the decomposition 
\begin{equation}
  \label{eq:m-c-decompose}
  \mm_{\cgot}(\mu,k)=\sum_{S\in\Sbb_\xi} P_S(\mu,k),
\end{equation} 
where
\begin{equation}
  \label{eq:Pol-S}
P_S(\mu,k)=\frac{1}{\vert H_S\vert}
\sum_{\gamma\in H_S^o}\gamma^{k\alpha_S-\mu}\int_{\overline{S}/T}
\frac{\todd(\overline{S}/T)}{D^\gamma(\Ncal_S)}\ 
e^{\frac{1}{2\pi}\Big(k\Omega_S -\langle k\xi-\mu,\omega_S\rangle\Big)}.  
\end{equation}
When $S$ is the principal open dense stratum $S_{max}$ the map $P_S$ is  
\begin{equation}
  \label{eq:Pol-S-max}
P_{max}(\mu,k)=\frac{\sum_{\gamma\in \Gamma_M}
\gamma^{k\alpha_{o}-\mu}}{\vert \Gamma_M\vert}
\int_{\Mcal_\xi}\todd(\Mcal_\xi)
e^{\frac{1}{2\pi}(k\Omega_\xi -\langle k\xi-\mu,\omega_\xi\rangle)}.  
\end{equation}
The term $\frac{\sum_{\gamma\in \Gamma_M}
\gamma^{k\alpha_{o}-\mu}}{\vert \Gamma_M\vert}$ 
is equal to $1$ when $(\mu,k)$ belongs to the lattice $\Xi(M,L)$ 
(see (\ref{eq:xi-lattice})), 
and is equal to $0$ in the other cases. 
From (\ref {eq:Pol-S}) we see that $P_S$ is a periodic polynomial of 
degree less than $\frac{dim(\overline{S}/T)}{2}$, and for $S=S_{max}$ 
we have on $\Xi(M,L)$
\begin{equation}
  \label{eq:Pol-S-bis}
P_{max}(\mu,k)=\frac{1}{(2\pi)^l}\int_{\Mcal_\xi}
\frac{(k\Omega_\xi -\langle k\xi-\mu,\omega_\xi\rangle)^l}{l !} +O(l-1)  
\end{equation}
 where $l=\frac{dim\Mcal_\xi}{2}$ and $O(l-1)$ denotes a polynomial of degree
less than $l-1$. If we use the polynomial $\hd_{\cgot}$ defined in Section 
\ref{DH-mesures} we can conclude our computations with the following 
\begin{prop}\label{prop:m-c-poly}
The map $\mm_{\cgot}$ is a periodic polynomial of degree $l=\frac{dim\Mcal_\xi}{2}$ 
supported on $\Xi(M,L)$. For $(\mu,k)\in\Xi(M,L)$ we have 
$$
\mm_{\cgot}(\mu,k)=\vert \Gamma_M\vert\frac{k^l}{(2\pi)^l}{\hd_{\cgot}}
(\frac{\mu}{k}) \ +\ O(l-1), 
$$
where $O(l-1)$ means a periodic polynomial of degree less than $l-1$.
\end{prop}

If $\cgot_\pm$ are two adjacent connected components of regular values of $\Phi$, 
we know after Theorem \ref{DH-theo-principal} that $\hd_{\cgot_+}\neq\hd_{\cgot_-}$ 
so we can conclude from Proposition \ref{prop:m-c-poly} that 
\begin{equation}
  \label{eq:mm-plus-moins}
 \mm_{\cgot_+}\neq\mm_{\cgot_-}. 
\end{equation}


\subsection{Jump formulas for the $\mm_{\cgot}$}\label{subsec:jump-formula}

Let $\cgot_+$ and $\cgot_-$ be two adjacent connected component of 
regular values of $\Phi$ separated by an hyperplane $\p$. The aim of 
this section is to compute the periodic polynomial  
$\mm_{\cgot_+}-\mm_{\cgot_-}$. 

Let $\beta\in\Lambda$ be the unimodular orthogonal vector 
to the hyperplane $\p\subset\tgot^*$ which is pointing out of $\cgot_-$. 
Let $T_{\p}\subset T$ be the subtorus of dimension $1$, with Lie algebra 
$\tgot_{\p}:=\{X\in\tgot\,\vert\, \langle \xi-\xi',X\rangle=0,\ 
\forall \xi,\xi'\in\p\}$. We make the choice of a decomposition  
$T=T_\p\times T/T_\p$, where $T/T_\p$ denotes a subtorus de $T$.

Let $\cgot'$ be the relative interior of 
$\overline{\cgot_+}\cap\overline{\cgot_-}$ in $\p$. Following Definition 
\ref{def:M-delta-quantum} we denote by $M_{\cgot'}$ the union of 
the connected component $Z$  of the fixed point set 
$M^{T_{\p}}$ for which $\Phi(Z)$ contains $\cgot'$. The symplectic
submanifold $M_{\cgot'}$ carries an Hamiltonian action of $T/T_{\p}$ 
with moment map $\Phi_{\cgot'}$ equal to the restriction of $\Phi$ on 
$M_{\cgot'}$. 

We consider two points $\xi_\pm\in \cgot_\pm$ such that $\xi=
\frac{1}{2}(\xi^+ + \xi^-)\in\cgot'$. We suppose furthermore 
that $\xi^+ -\xi^-$ is orthogonal to $\p$. Using the identification 
$\tgot^*\simeq\tgot$ given by the scalar product the vector 
$\gamma=\frac{1}{2}(\xi_+-\xi_-)$, seen  as a vector of $\tgot_\p$, 
belongs to $\Rbb^{>0}\beta$. We noticed in Section
\ref{subsec:Jump-dh} that for all 
$m\in\Phi^{-1}(\xi)$ the stabilizer $\tgot_m$ is either 
equal to $\tgot_\p$ or to $\{0\}$. Then $\xi$ is 
a regular value of $\Phi_{\cgot'}$ and there exists an open 
$T$-invariant neighborhood $\Ucal$ of $\Phi^{-1}(\xi)$ in $M$ such that 
for all $m\in\overline{\Ucal}$ either 
$\tgot_m:=\{0\}$, or $\tgot_m=\tgot_\p$ and $\Phi(m)\in\p$. 

One see easily that the couple $(\Ucal,\xi)$ is good and part $b)$ of 
Proposition \ref{prop:U-xi} tells us that
\begin{equation}
  \label{eq:RR-U-pm}
 RR^\xi_\Ucal(M,-)=RR^{\xi_-}_\Ucal(M,-)=RR^{\xi_+}_\Ucal(M,-) 
 \end{equation}
when $\xi_\pm$ are close enough to $\xi$. Since 
$\Ucal\cap\crpt(\parallel\Phi-\xi\parallel^2)=\Phi^{-1}(\xi)$ we have
\break 
$RR^\xi_\Ucal(M,-)=RR^\xi_0(M,-)$. If $\xi_\pm$ are close enough to $\xi$
 we have 
 \begin{equation}
   \label{eq:U-cap-xi-pm}
 \Ucal\cap\crpt(\parallel\Phi-\xi_\pm\parallel^2)=
\Phi^{-1}(\xi_\pm)\bigcup M^{\gamma}\cap\Phi^{-1}(\xi).  
 \end{equation}
The former decomposition is due to (\ref{eq:decom-kirwan}) and to 
the fact that the stabiliser of $\tgot$ on $\Ucal$ are either equal to 
$\tgot_\p$ or to $\{0\}$. Notice that $\xi_-
+\gamma=\xi_+ +\gamma=\xi$. The decomposition (\ref{eq:U-cap-xi-pm})
gives 
\begin{equation}
  \label{eq:RR-U-mp}
RR_{\Ucal}^{\xi_\pm}(M,-)=RR_{0}^{\xi_\pm}(M,-)+ RR_{\mp\gamma}^{\xi_\pm}(M,-),  
\end{equation}
where $RR_{\gamma}^{\xi_-}(M,-)$ (resp. $RR_{-\gamma}^{\xi_+}(M,-)$) is the 
Riemann-Roch character localized on $M^\gamma\cap\Phi^{-1}(\xi)$ 
by the vectors field $\Hcal- (\xi_-)_M$ (resp. $\Hcal- (\xi_+)_M$). Now 
(\ref{eq:RR-U-pm}) and (\ref{eq:RR-U-mp}) prove the following 

\begin{prop}\label{RR-0-diff}
If $\xi_\pm$ are close enough to $\p$, we have 
$$
RR_0^{\xi_+}(M,-)-RR_0^{\xi_-}(M,-)=
RR_{\gamma}^{\xi_-}(M,-)-RR_{-\gamma}^{\xi_+}(M,-).
$$
\end{prop}

We know from Proposition \ref{prop:m-poly} that $\mm_{\cgot_\pm}(\mu,k)$ is 
equal to the $\mu$-mutiplicity of $RR_0^{\xi_\pm}(M,L^{\otimes k})$. Hence
$\mm_{\cgot_+}(\mu,k)-\mm_{\cgot_-}(\mu,k)$ is equal to the $\mu$-mutiplicity 
of \break $RR_{\gamma}^{\xi_-}(M,L^{\otimes k})-RR_{-\gamma}^{\xi_+}(M,L^{\otimes k})$.

Let $N_{\p}$ be the normal bundle of $M^{T_\p}$ in $M$, and let 
$\left[\wedge_{\Cbb}^{\bullet}\overline{N_{\p}}\right]^{-1}_{\pm\beta}$ 
be the polarized inverses of $\wedge_{\Cbb}^{\bullet}\overline{N_{\p}}$ 
(see Proposition \ref{def:lambda-N}). Since $\xi=\xi_+ -\gamma=
\xi_- +\gamma$ and $\gamma\in \Rbb^{>0}\beta$, 
the localization (\ref{eq:loc-RR-beta-bis}) gives 
\begin{equation}\label{eq:loc-xi-1}
RR_{\gamma}^{\xi_-}(M,E)= \sum_{Z}
RR_0^{\xi}\left(Z,E\vert_{Z}\otimes
\left[\wedge_{\Cbb}^{\bullet}\overline{N_Z}\right]^{-1}_{\beta}\right)
\end{equation}
and
\begin{equation}
  \label{eq:loc-xi-2}
  RR_{-\gamma}^{\xi_+}(M,E)= \sum_{Z}
RR_0^{\xi}\left(Z,E\vert_{Z}\otimes
\left[\wedge_{\Cbb}^{\bullet}\overline{N_Z}\right]^{-1}_{-\beta}\right).
\end{equation}

In (\ref{eq:loc-xi-1}) and (\ref{eq:loc-xi-2}) the sum  are 
taken over the connected components $Z$ of $M_{\cgot'}\subset M^{T_\p}$ and 
we denote $N_Z$ the restriction 
of the bundle $N_{\p}$ to $Z$. Let us make few remarks concerning 
the maps $RR_0^{\xi}(Z,-):\K_T(Z)\to R^{-\infty}(T)$. Since $T_\p$ acts trivially on 
$Z$, the decomposition  $T=T/T_\p\times T_\p$ induces a canonical isomorphism 
$\K_T(Z)\simeq\K_{T/T_\p}(Z)\otimes R(T_\p)$: i.e. every $T$-equivariant 
vector bundle $E\to Z$ decomposes as
\begin{equation}
  \label{eq:t-p-decompose}
 E= \sum_{\alpha\in\Lambda^*_{\tgot_\p}}E^\alpha\otimes\Cbb_\alpha.
\end{equation}
Here $\Lambda^*_{\tgot_\p}=\Lambda^*\cap\tgot_\p^*$ is the set of weights for the 
subtorus $T_\p$, each $E^\alpha$ is a $T/T_\p$-equivariant vector bundle on $Z$ and 
$\Cbb_\alpha$ denotes the one dimensional $T_\p$-representation associated to 
$\alpha\in\Lambda^*_{\tgot_\p}$.

For every $T$-equivariant vector bundle $E\to Z$, the character 
$RR_0^{\xi}(Z,E)$ is equal to the $T$-equivariant index of 
the $T$-transversally elliptic symbol $\Thom_{\xi}(\Vcal)\otimes
p^*(E)$, where $\Vcal$ is a small neighborhood of $\Phi^{-1}(\xi)\cap
Z$ in $Z$ (see Definition \ref{def:RR-U-xi}).  Since the $T_\p$ action 
is trivial on $Z$ the symbol $\Thom_{\xi}(\Vcal)$ is also $T/T_\p$ 
transversally elliptic and the action of $T_\p$ is trivial on it. 
We have then 
\begin{equation}
  \label{eq:rr-xi-o}
RR_0^{\xi}(Z,E)=\sum_{\alpha\in\Lambda^*_{\tgot_\p}}
RR_0^{\xi}(Z,E^\alpha)\otimes\Cbb_\alpha
\end{equation}
where $RR_0^{\xi}(Z,E^\alpha)$ belongs to $R^{-\infty}(T/T_\p)$. 
Since $\xi$ is a regular value 
of $\Phi_{\cgot'}$, the character $RR_0^{\xi}(Z,E^\alpha)$ is computed by Theorem 
\ref{eq:RR-xi-O} applied to the Hamiltonian $T/T_\p$-manifold $Z$. 
For every $T$-vector bundle $E\to M$ we define the familly 
$\Ecal^{\mu_1}_{\xi,\mu_2}$ of orbifold 
vector bundles over the reduced space $\Zcal_\xi=
Z\cap \Phi^{-1}(\xi)/(T/T_\p)$ by 
\begin{equation}
  \label{eq:E-xi-mu}
  \Ecal^{\mu_1}_{\xi,\mu_2}:=
\left(E^{\mu_1}\otimes\Cbb_{-\mu_2}\right)\vert_{\Phi^{-1}(\xi)\cap Z}/(T/T_\p).
\end{equation}
Here $\mu_1\in\Lambda^*_{\tgot_\p}$ and  $\mu_2\in\Lambda^*_{\tgot/\tgot_\p}=
\Lambda^*\cap(\tgot/\tgot_\p)^*$. Finally 
Theorem \ref{eq:RR-xi-O} and (\ref{eq:rr-xi-o}) give the following 
\begin{eqnarray}
  \label{eq:indice-t-p}
  RR_0^\xi(Z,E)&=&
\sum_{\mu_1\in\Lambda^*_{\tgot_\p}}\sum_{\mu_2\in\Lambda^*_{\tgot/\tgot_\p}}
RR(\Zcal_\xi, \Ecal^{\mu_1}_{\xi,\mu_2})\ \otimes 
\underbrace{\Cbb_{\mu_1}}_{\in \, R(T_\p)} \otimes
\underbrace{\Cbb_{\mu_2}}_{\in \, R(T/T_\p)} 
\nonumber \\
               &=&\sum_{\mu\in\Lambda^*}
RR(\Zcal_\xi, \Ecal^{\mu_1}_{\xi,\mu_2})\ \Cbb_{\mu} \quad .
\end{eqnarray}

In (\ref{eq:indice-t-p}) we write $\mu\in\Lambda^*$ as a sum of 
$\mu_1\in\Lambda^*_{\tgot_\p}$ with
$\mu_2\in\Lambda^*_{\tgot/\tgot_\p}$ 
so that $\Cbb_{\mu}\in R(T)$ is equal to the tensor product 
$\Cbb_{\mu_1}\otimes\Cbb_{\mu_1}$.
 
\bigskip

Now we finish the computation of the periodic polynomial 
$\mm_{\cgot_+}(\mu,k)-\mm_{\cgot_-}(\mu,k)$. The hyperplane $\p$ is 
defined by the equation $\frac{\langle\xi,\beta\rangle}{2\pi} = 
r_\p,\, \xi\in\tgot^*$, for some
$r_\p\in\Zbb$. We have the decomposition 
\begin{equation}
  \label{eq:diff-decompose}
 \left[\wedge_{\Cbb}^{\bullet}\overline{N_Z}\right]^{-1}_{\beta}-
\left[\wedge_{\Cbb}^{\bullet}\overline{N_Z}\right]^{-1}_{-\beta}=
\sum_{\alpha\in\Lambda^*_{\tgot_\p}} (-1)^{n_Z(\alpha)}\ 
S_Z^\alpha\otimes\Cbb_\alpha
\end{equation}
where $S_Z^\alpha\to Z$ is a complex $T/T_\p$-vector bundle, and the integers 
$n_Z(\alpha)$ are defined by the the following relations
\begin{eqnarray}
  \label{eq:n-Z-alpha}
 n_Z(\alpha)&=& {\rm rk}_\Cbb( N_Z^{+,\beta})\quad  if\quad  
\langle\alpha,\beta\rangle\geq 0\nonumber\\
 n_Z(\alpha)&=& {\rm rk}_\Cbb( N_Z^{+,-\beta}) +1 \quad if \quad 
\langle\alpha,\beta\rangle< 0.
\end{eqnarray}

Let $\beta^*$ be the vector of $\tgot_\p^*$ which is defined by the
relation $\langle\beta^*,\beta\rangle=2\pi$, so that 
$\Lambda^*_{\tgot_\p}=\Zbb\beta^*$. The $T_\p$-weight on 
$\det(N_Z^{+,\pm\gamma})$ is $\pm s_Z^\pm\beta^*$ where 
$s_Z^\pm\in\Nbb$ is the absolute value of the trace of 
$\frac{1}{2\pi}\Lcal_\beta$ on $N_Z^{+,\pm\gamma}$. Since the 
$T_\p$-weights on $S^k(N_Z^{+,\gamma})$ (resp. 
$S^k(N_Z^{+,-\gamma})$) are of the form $p\beta^*$ with 
$p\geq 0$ (resp. $p\leq 0$), we see that 
\begin{equation}
  \label{eq:S-nul}
  S_Z^{\alpha}=0 \quad {\rm if }\quad -s_Z^- 
< \frac{\langle\alpha,\beta\rangle}{2\pi}< s_Z^+.
\end{equation}

The line bundles $L^{\otimes k}\vert_{M_{\cgot'}}, k\in \Zbb$ 
can be considered either as 
$T$-vector bundles or as $T/T_\p$-vector vector bundles: we denote them 
respectively by $L^{\otimes k}\vert_T$ and by $L^{\otimes k}\vert_{T/T_\p}$. 
The $T_\p$-weight on $L^{\otimes k}\vert_{M_{\cgot'}}$ is equal to 
$kr_\p\beta^*$, hence we have 
\begin{equation}
  \label{eq:L-k}
  L^{\otimes k}\vert_T=L^{\otimes k}\vert_{T/T_\p}\otimes \Cbb_{kr_\p\beta^*} \ .
\end{equation}

For every $\mu\in\Lambda^*$ and $Z\subset M_{\cgot'}$ we define the 
orbifold vector bundle\footnote{More precisely 
$\Sbb^k_{Z,\mu}$ is $\pm 1$ times a orbifold vector bundle over 
the reduced space $\Zcal_\xi$.} $\Sbb^k_{Z,\mu}$ 
on the reduced space $\Zcal_\xi=(Z\cap\Phi^{-1}(\xi))/(T/T_\p)$ by 
\begin{equation}
  \label{eq:bundle-xi-mu}
  \Sbb^k_{Z,\mu}:=(-1)^{n_Z(\mu_1-kr_\p\beta^*)}
\left(L^{\otimes k}\vert_{T/T_\p}\otimes 
S_Z^{\mu_1-kr_\p\beta^*}\otimes\Cbb_{-\mu_2}\right)
\vert_{Z\cap\Phi^{-1}(\xi)}/(T/T_\p),
\end{equation}
where $\mu=\mu_1+\mu_2$ with $\mu_1\in\Lambda^*_{\tgot_\p}$ and 
$\mu_2\in\Lambda^*_{\tgot/\tgot_\p}$. 

\medskip

\begin{theo}\label{th:th-principal-mult}
Let $\Sbb^k_{\xi,\mu}$ be the orbifold vector bundle on 
$\Mcal^\p_\xi=\cup_{Z\subset M_{\cgot'}}\Zcal_\xi$ 
which is equal to $\Sbb^k_{Z,\mu}$ over each connected component $\Zcal_\xi$. 
  For every $(\mu,k)\in\Lambda^*\times\Zbb$, we have 
  \begin{equation}
    \label{eq:mm-diff}
 \mm_{\cgot_+}(\mu,k)-\mm_{\cgot_-}(\mu,k)=RR(\Mcal^\p_\xi,\Sbb^k_{\xi,\mu}).   
  \end{equation}
In particular $\mm_{\cgot_+}(\mu,k)=\mm_{\cgot_-}(\mu,k)$ if
\begin{equation}
  \label{eq:m-egale}
 - s^- < \frac{\langle\mu,\beta\rangle}{2\pi}- k r_\p < s^+, 
\end{equation}
with $s^\pm=\inf_{Z\subset M_{\cgot'}} s^\pm_Z$.
\end{theo}

\medskip

{\em Proof.} After (\ref{eq:loc-xi-1}), (\ref{eq:loc-xi-2}), and Proposition 
\ref{RR-0-diff} 
we know that $\mm_{\cgot_+}(\mu,k)-\mm_{\cgot_-}(\mu,k)$ is equal to the 
$\mu$-mutiplicity in 
\begin{eqnarray}\label{eq:RR-T-delta}
 \lefteqn{ \sum_{Z\subset M_{\cgot'}} RR_{0}^{\xi}\left(Z,L^{\otimes k}\vert_Z\otimes
(\left[\wedge_{\Cbb}^{\bullet}\overline{N_Z}\right]^{-1}_{\beta}-
\left[\wedge_{\Cbb}^{\bullet}\overline{N_Z}\right]^{-1}_{-\beta})\right)}\nonumber\\
&=&\sum_{\alpha\in\Lambda^*_{\tgot_\p}}\sum_{Z\subset M_{\cgot'}}
(-1)^{n_Z(\alpha)}\underbrace{RR_{0}^{\xi}\left(Z,L^{\otimes k}\vert_{T/T_\p}
\otimes S_Z^\alpha\right)}_{\in\ R^{-\infty}(T/T_\p)}\otimes
\underbrace{\Cbb_{\alpha+kr_\p\beta^*}}_{\in \ R(T_\p)}\nonumber \\
&=&\sum_{\mu\in\Lambda^*} RR(\Mcal^\p_\xi,\Sbb^k_{\xi,\mu}) \ \Cbb_\mu 
\end{eqnarray}

Here (\ref{eq:RR-T-delta}) is a direct consequence of
(\ref{eq:indice-t-p}) modulo the
shift $\mu_1\mapsto\mu_1 - kn_\p\beta^*$ in 
$\Lambda^*(\tgot_\p)$. We get (\ref{eq:m-egale}) using the vanishing 
conditions (\ref{eq:S-nul}). $\Box$

\medskip

When condition (\ref{eq:m-egale}) is the optimal one? In other word, do we have 
$\mm_{\cgot_+}(\mu,k)\neq \mm_{\cgot_-}(\mu,k)$ for some $\mu\in\Lambda^*$ such 
that $\frac{\langle\mu,\beta\rangle}{2\pi}- k r_\p = s^+$ , and do we have 
$\mm_{\cgot_+}(\mu,k)\neq \mm_{\cgot_-}(\mu,k)$ for some $\mu\in\Lambda^*$ such 
that $\frac{\langle\mu,\beta\rangle}{2\pi}- k r_\p = s^-$ ?

\medskip

Until the end of this section we consider couples $(\mu,k)$ such that 
$\mu=$ \break $(k r_\p +s^+)\beta^* +\mu_2$ 
with $\mu_2\in\Lambda^*_{\tgot/\tgot_\p}$. For the orbifold vector bundle  
$\Sbb^k_{Z,\mu}$ we have $\Sbb^k_{Z,\mu}= (-1)^{{\rm rk}_\Cbb( N_Z^{+,\beta})} 
\widetilde{\Sbb}^k_{Z,\mu_2}$ where
\begin{equation}
  \label{eq:vector-S-maximal}
 \widetilde{\Sbb}^k_{Z,\mu_2}= 
\left(L^{\otimes k}\vert_{T/T_\p}\otimes 
\det( N_Z^{+,\beta})\otimes\Cbb_{-\mu_2}\right)
\vert_{Z\cap\Phi^{-1}(\xi)}/(T/T_\p)
\end{equation}
when  $s^+=s^+_Z$ and $\Sbb^k_{Z,\mu}=0$ when $s^+<s^+_Z$. Hence  
$\mm_{\cgot_+}(\mu,k)-\mm_{\cgot_-}(\mu,k)$ is equal to
\begin{equation}
  \label{eq:mm-pm-maximal}
  \sum_{Z,s^+_Z=s^+} (-1)^{{\rm rk}_\Cbb( N_Z^{+,\beta})}
RR(\Zcal_\xi,\widetilde{\Sbb}^k_{Z,\mu_2})
\end{equation}

Now we use the results of Section \ref{subsec:RR-theo} 
to compare the behaviour of the maps
\begin{eqnarray}
  \label{eq:RR-Z-mu-2}
  \Lambda^*_{\tgot/\tgot_\p}\times\Zbb&\longrightarrow &\Zbb,\\ 
(\mu_2,k)&\longmapsto& RR(\Zcal_\xi,\widetilde{\Sbb}^k_{Z,\mu_2})\nonumber
\end{eqnarray}
for the different components $Z$ satisfying $s^+_Z=s^+$. 
Let $\Gamma_Z\subset T/T_\p$ be the generic stabiliser of 
$T/T_\p$ on a component $Z$. Let $\alpha_Z,\delta_Z\in 
\Lambda^*_{\tgot/\tgot_\p}$ such that the action of 
$\Gamma_Z$ on the fibers of $L_{\vert Z}$ and $\det( N_Z^{+,\beta})$ 
are respectively $t\to t^{\alpha_Z}$ and 
$t\to t^{\delta_Z}$. After Remark \ref{lem:discrete} we know that 
the map (\ref{eq:RR-Z-mu-2}) is supported on the subset 
\begin{equation}
  \label{eq:xi-Z}
 \Xi_Z:=\{(\mu_2,k)\in\Lambda^*_{\tgot/\tgot_\p}\times\Zbb\ \vert\ 
t^{k\alpha_Z+\delta_Z+\mu_2}=1,\  \forall\ t\in\Gamma_Z\}.
\end{equation}
The only difference with  the computations done in Section 
\ref{subsec:RR-theo} is the presence of the line bundle 
$\det( N_Z^{+,\beta})$. But this do not change the global behaviour of 
the map (\ref{eq:RR-Z-mu-2}) on $\Xi_Z$: 
it is a periodic polynomial map of degree
$l_Z=\dim(\Zcal_\xi)/2$ and we have 
\begin{equation}
  \label{eq:pol-behaviour-RR-mu2}
RR(\Zcal_\xi,\widetilde{\Sbb}^k_{Z,\mu_2})=\frac{1}{(2\pi)^{l_Z}}\int_{\Zcal_\xi}
(k\Omega_{\Zcal_\xi} -\langle k\xi-\mu_2,\omega_{\Zcal_\xi}\rangle)^{l_Z} +O(l_Z-1)  
\end{equation}
for all $(\mu_2,k)\in\Xi_Z$. We see then that 
the sum (\ref{eq:mm-pm-maximal}) does not vanish for large values of $(\mu_2,k)$ 
when the number $(-1)^{{\rm rk}_\Cbb( N_Z^{+,\beta})}$ are equal for all the 
components $Z$ of maximal dimension. Finally we can conclude with

\begin{prop}\label{prop:inegalite-optimale}
Consider the set of components $Z\subset M_{\cgot'}$ for which  $s^+_Z$ is minimal. 
Among them consider the subset  $\Fcal$ where $\dim(Z)$ is maximal. If the 
integers ${\rm rk}_\Cbb( N_Z^{+,\beta})$, $Z\in\Fcal$ have the same parity, then 
the condition ``$\frac{\langle\mu,\beta\rangle}{2\pi}- k r_\p < s^+$'' is 
optimal in (\ref{eq:m-egale}). 

In the same way, consider the set of components 
$Z\subset M_{\cgot'}$ for which  $s^-_Z$ is minimal. 
Among them consider the subset  $\Fcal'$ where $\dim(Z)$ 
is maximal. If the \break integers ${\rm rk}_\Cbb( N_Z^{+,\beta})$, 
$Z\in\Fcal'$ have the same parity, then 
the condition \break ``$-s^-<\frac{\langle\mu,\beta\rangle}{2\pi}- k r_\p $'' is 
optimal in (\ref{eq:m-egale}). 
\end{prop}


\section{Multiplicities of group representations}\label{sec:cas-coadjoint}

Let $K$ be a semi-simple compact Lie group with Lie algebra $\kgot$, and let 
$T$ be a maximal torus in $K$ with Lie algebra $\tgot$. In this section 
we denote $(-,-)$ the scalar product on $\kgot$ induced by the 
Killing form, and we keep the same notation for the induced 
scalar products on $\tgot^*$ and on $\tgot$.


An element $\lambda\in \tgot^*$ belong to the real weight lattice 
$\Lambda^*\subset \tgot^*$ if $\imath\, \lambda$ is the differential of 
a character of $T$ that we denote $t\mapsto t^\lambda$:
if $t=\exp X$ then $t^{\lambda}:=e^{\imath\langle\lambda,X\rangle}$. Let 
$\Rgot\subset\Lambda^*$ be the set of roots for the action of 
$T$ on $\kgot\otimes\Cbb$, and let $\Lambda^*_R$ be the sub-lattice 
of $\Lambda^*$ generated by $\Rgot$. We choose 
a system of positive roots $\Rgot^+\subset\Rgot$, and we denote $\tgot^*_+$ the 
corresponding Weyl chamber.

The irreducible representations of $K$ are parametrized by the set $\Lambda^*_+=
\Lambda^*\cap\tgot^*_+$. For $\lambda\in\Lambda^*_+$ 
we denote by $V_\lambda$ the irreducible representation 
of $K$ with heighest weight $\lambda$. Here we are interested 
in the $T$-multiplicities in $V_\lambda\vert_T$. Let 
$\mm: \Lambda^*\times \Lambda^*_+\to\Nbb$ be the map defined by 
\begin{equation}
  \label{eq:T-mult}
 V_\lambda\vert_T=\sum_{\mu\in\Lambda^*}\mm(\mu,\lambda)\ \Cbb_\mu
\end{equation}
for every $\lambda\in\Lambda^*_+$. 

\begin{defi}
  \label{def:mm-lambda}
For every $\lambda\in\Lambda^*_+$, we denote $\mm^\lambda: 
\Lambda^*\times \Zbb^{>0}\to\Nbb$ the map defined by 
$\mm^\lambda(\mu,k)=\mm(\mu,k\lambda)$. So
$\mm^\lambda(\mu,k)$ is equal to the multiplicity 
of $\Cbb_\mu$ in $V_{k\lambda}\vert_T$.
\end{defi}

\subsection{Borel-Weil Theorem}

First we recall the realization of the $K$-representation $V_\lambda$ given by the 
Borel-Weil Theorem. The coadjoint orbit $K\cdot\lambda$ is equipped with the 
Kirillov-Kostant-Souriau symplectic form $\Omega$ which is defined by:
\begin{equation}
  \label{eq:KKS}
  \Omega(X_M,Y_M)_m=\langle m,[X,Y]\rangle,\quad {\rm for}\quad  
m\in K\cdot\lambda\quad {\rm and}
\quad X,Y\in\kgot.
\end{equation}
The action of $K$ on $K\cdot\lambda$ is Hamiltonian with moment map 
$K\cdot\lambda\croc\kgot^*$ equal to the inclusion. The action of 
$T$ on $K\cdot\lambda$ is also Hamiltonian with moment map 
$\Phi:K\cdot\lambda \to\tgot^*$ equal to the composition of the 
inclusion $K\cdot\lambda\croc\kgot^*$ with the projection map 
$\kgot^*\to\tgot^*$.

There exists a unique $K$-invariant complex structure 
on $K\cdot\lambda$ compatible with the symplectic form. In this situation 
the Kostant-Souriau prequantum line bundle over $K\cdot\lambda$ is 
$$
\Cbb_{[\lambda]}=K\times_{K_\lambda}\Cbb_\lambda.
$$
Here we use the canonical identification $K/K_\lambda\simeq K\cdot\lambda,\  
[k]\mapsto k\cdot\lambda$, where $K_\lambda$ is the stabilizer of
$\lambda$ in $K$. The line bundle $\Cbb_{[\lambda]}$ over 
the complex manifold $K\cdot\lambda$ 
carries a canonical holomorphic structure. If one work with the 
symplectic form $k\Omega$, for an integer $k\geq 1$, the 
corresponding Kostant-Souriau prequantum 
line bundle is $\Cbb_{[\lambda]}^{\otimes k}=
K\times_{K_\lambda}\Cbb_{k\lambda}=\Cbb_{[k\lambda]}$.

Let $\Hcal^q(K\cdot\lambda,\Cbb_{[\lambda]}^{\otimes k})$ 
be $q$th cohomology group of the sheaf of holomorphic 
section of $\Cbb_{[\lambda]}^{\otimes k}$ over  $K\cdot\lambda$. 
The Borel-Weil Theorem tells us that 
\begin{equation}
  \label{eq:BW-1}
 \Hcal^0(K\cdot\lambda,\Cbb_{[\lambda]}^{\otimes k})= V_{k\lambda}
\end{equation}
and 
\begin{equation}
  \label{eq:BW-2}
 \Hcal^q(K\cdot\lambda,\Cbb_{[\lambda]}^{\otimes k})= 0\quad {\rm for}\quad q\geq 1.
\end{equation}

If $RR^K(K\cdot\lambda,-):\K_K(K\cdot\lambda)\to R(K)$ is the $K$-Riemann-Roch 
character defined by the compatible complex structure, (\ref{eq:BW-1}) and 
(\ref{eq:BW-2}) give 
\begin{equation}
  \label{eq:RR-V-lambda}
  RR^K(K\cdot\lambda,\Cbb_{[\lambda]}^{\otimes k})=
V_{k\lambda} \quad {\rm in}\quad  R(K)
\end{equation}

Now if we denote by $RR(K\cdot\lambda,-):\K_T(K\cdot\lambda)\to R(T)$ 
the $T$-equivariant Riemann-Roch character, we have $V_{k\lambda}\vert_T=
RR(K\cdot\lambda,\Cbb_{[\lambda]}^{\otimes k})$.  
The multiplicity fonction $\mm^\lambda:
\Lambda^*_+\times \Nbb^*\to\Nbb$ is characterized 
by the relation 
\begin{equation}
  \label{eq:T-mult-bis}
 RR(K\cdot\lambda,\Cbb_{[\lambda]}^{\otimes k})= 
\sum_{\mu\in\Lambda^*}\mm^\lambda(\mu,k)\ \Cbb_\mu,\quad {\rm in}\ R(T) ,
\end{equation}
for $k\geq 1$.

The sub-lattice $\Lambda^*_R$ of $\Lambda^*$ generated by the 
roots is characterized by the (finite) center $Z(K)$ of $K$ as 
follows. For $\alpha\in\Lambda^*$ we have
\begin{equation}
  \label{eq:lambda-R}
\lambda\in\Lambda^*_R\ \Longleftrightarrow\  t^\lambda= 1,\  
\forall t \in Z(K),   
\end{equation}
and for $t\in T$ we have 
$t\in Z(K)\ \Longleftrightarrow\  t^\lambda= 1,\  \forall 
\lambda\in\Lambda^*_R$. The finite abelian group  
$\Lambda^*/\Lambda^*_R$ is then naturally identified  with the 
dual of $Z(K)$. We have the following well-known fact.

\begin{prop}
  \label{prop:support-R}
We have $\mm^\lambda(\mu,k)\neq 0$ only if $\mu-k\lambda\in\Lambda^*_R$.
\end{prop}

{\em Proof.} The center $Z(K)$ of $K$ acts trivially on $K\cdot\lambda$ and 
with the character $t\in Z(K)\mapsto t^{k\lambda}$ on 
the fibers of the line bundle $\Cbb_{[\lambda]}^{\otimes k}$. Since 
$\mm^\lambda(\mu,k)$ is equal to the dimension of the $T$-invariant subspace of
$RR(K\cdot\lambda,\Cbb_{[\lambda]}^{\otimes k})\otimes\Cbb_{-\mu}$, we have following
 Lemma \ref{lem:discrete} that $\mm^\lambda(\mu,k)\neq 0$ only if 
$t^{\mu-k\lambda}=1$, $\forall t\in Z(K)$. We conclude then with 
(\ref{eq:lambda-R}). $\Box$

In this section we are interested in the periodic polynomials
\begin{equation}
  \label{eq:multiplicite-c}
\mm_{\cgot}^\lambda : \Lambda^* \times \Zbb \longrightarrow  \Zbb.
\end{equation}
defined for every connected component $\cgot\subset\tgot^*$ 
of regular values of the moment map $\Phi:K\cdot\lambda \to\tgot^*$: 
the map $\mm_{\cgot}^\lambda$ coincide with $\mm^\lambda$ on the set 
$\{(\mu,k)\in\Lambda^* \times \Zbb^{>0}\ \vert\ \mu\in k\cgot\}$. 
Like in Proposition \ref{prop:support-R}, the formula given in 
Proposition \ref{prop:m-poly} for $\mm_{\cgot}^\lambda$ tells us 
that $\mm_{\cgot}^\lambda(\mu,k)\neq 0$ only if 
$\mu-k\lambda\in\Lambda^*_R$.

To apply Theorem \ref{th:th-principal-mult} to the periodic polynomials
$\mm_{\cgot}^\lambda$, we have to compute the critical values of the 
moment map $\Phi:K\cdot\lambda \to\tgot^*$.

\subsection{Critical points of $\Phi:K\cdot\lambda \to\tgot^*$}

Let $\{\alpha_1,\cdots,\alpha_{dim T}\}$ be the simple roots of the set 
$\Rgot_+$ of positive weights. The fundamental weights 
$\varpi_k, 1\leq k\leq \dim T$ are defined by the conditions 
\begin{equation}
  \label{eq:fundamental-w}
 2\frac{(\varpi_i,\alpha_j)}{\vert\alpha_j\vert^2}=\delta_{i,j}\quad {\rm for\ all} 
\quad 1\leq i,j\leq \dim T. 
\end{equation}
Recall that the fundamental weights generate the lattice $\Lambda^*_{alg}$ of 
{\em algebraic} integral element of $\tgot^*$. 
We have $\Lambda^*\subset \Lambda^*_{alg}$ and equality holds 
only if $K$ is simply-connected.

Let $W$ be the Weyl group of $(K,T)$. We will see  
\begin{equation}
  \label{eq:Fcal}
 \Fcal= \{\sigma\cdot\varpi_i\ \vert\ \sigma\in W,\  1\leq i\leq \dim T\}.
\end{equation}
as a subset of $\tgot$ modulo the identification 
$\tgot\simeq\tgot^*$ given by the scalar product.  The singular 
points of $\Phi$ have the following nice description. This result
first appeared in Heckman's Thesis \cite{Heckman-PhD}.

\begin{prop}[\cite{Heckman-PhD,GLS96}]
  \label{prop:pt-critiques}
The critical points of $\Phi:K\cdot\lambda \to\tgot^*$ is the union 
of the fixed points set $(K\cdot\lambda)^\beta,\ \beta\in\Fcal$. For each 
$\beta\in\Fcal$ we have 
$$
(K\cdot\lambda)^\beta=\bigcup_{\sigma\in W} K^\beta\cdot \sigma\lambda.
$$
Here $K^\beta$ is the stabilizer of $\beta$ in $K$. 
\end{prop}

The fixed points of the action of $T$ on $K\cdot\lambda$ characterize the 
image of $\Phi$ completely: $\Phi(K\cdot\lambda)$ is the convex polytope 
\begin{equation}
  \label{eq:convex}
  {\conv}(W\cdot\lambda):=\ {\rm convex\ hull\ of\ } W\cdot\lambda.
\end{equation}
This result was first proved by Kostant \cite{Kostant73}. 
This is particular case of the convexity 
theorem of Atiyah, Guillemin and Sternberg 
\cite{Atiyah.82,Guillemin-Sternberg82.bis}. From Proposition 
\ref{prop:pt-critiques},
we know that the singular values of $\Phi:K\cdot\lambda \to\tgot^*$ are the 
convex polytopes 
\begin{equation}
  \label{eq:convex-singular}
  {\conv}(W^\beta\cdot\sigma\lambda),\quad \beta\in\Fcal,\ \sigma\in W/W^\beta.
\end{equation}
where $W^\beta$ is the stabilizer of $\beta$ in $W$, i.e. $W^\beta$ 
is the Weyl group of $(K^\beta,T)$. Each convex polytope 
${\conv}(W^\beta\cdot\sigma\lambda)$ lies in the hyperplane 
\begin{equation}
  \label{eq:delta-beta-sigma}
  \Delta_{\beta,\sigma}=\{\xi\in\tgot^*\ \vert\ 
(\xi-\sigma\lambda,\beta)=0\}.
\end{equation}

Note that $K^\beta\cdot \sigma\lambda$ coincide with $K^\beta\cdot \sigma'\lambda$ 
if and only if $\sigma\lambda\in W^\beta\sigma'\lambda$. Two polytopes 
${\conv}(W^\beta\cdot\sigma\lambda)$ and
${\conv}(W^\beta\cdot\sigma'\lambda)$ intersect if and only if 
$\langle\sigma\lambda,\beta\rangle=\langle\sigma'\lambda,\beta\rangle=0$. 

\begin{defi}\label{def:lambda-gene}
  An element $\lambda\in\Lambda^*_+$ is {\em generic} 
if for every fundamental root $\varpi_i$ 
we have 
\begin{equation}\label{eq:lambda-gene}
 (\sigma\lambda,\varpi_i)\neq(\sigma'\lambda,\varpi_i)
\end{equation}
each times that $\sigma\lambda\notin W^i\sigma'\lambda$ 
(here $W^i=\{\sigma\in W\,\vert\,\sigma\varpi_i=\varpi_i\}$). 
\end{defi}

This  condition of genericity imposes that the hyperplanes 
$\Delta_{\beta,\sigma}$ and 
$\Delta_{\beta,\sigma'}$ are distinct each times that the submanifolds 
$K^\beta\cdot \sigma\lambda$ and $K^\beta\cdot \sigma'\lambda$ are not equal.  

\begin{ex}
  \label{ex:su-4}
Consider the case of $\hbox{\rm SU}(4)$. Take the coadjoint orbit trough $\lambda=
(2,1,-1,-2)$, and $\sigma,\sigma'$ such that $\sigma\lambda=(2,-2,1,-1)$
and $\sigma'\lambda=(1,-1,2,-2)$. Take the fundamental weight 
$\varpi_2=\frac{1}{2}(1,1,-1,-1)$. In this case $\lambda$ is not 
``generic'' since $\sigma\lambda\notin W^i\sigma'\lambda$ but 
$(\sigma\lambda,\varpi_2)=(\sigma'\lambda,\varpi_2)=0$.
\end{ex}

\subsection{Main theorems}

Let $\cgot_+$ and $\cgot_-$ be two adjacent connected components of 
regular values of $\Phi:K\cdot\lambda\to\tgot^*$. The intersection 
$\overline{\cgot_+}\cap\overline{\cgot_+}$ is contained in an hyperplane 
orthogonal  to $\beta\in\Fcal$.

\begin{defi}
Let $\Acal(\cgot_+,\cgot_-)$ be the set of all $\sigma\in  W/W^\beta$ such 
that the convex 
polytope ${\conv}(W^\beta\cdot\sigma\lambda)$  
contains $\overline{\cgot_+}\cap\overline{\cgot_+}$. 
\end{defi}

The set 
$$
\bigcup_{\sigma\in\Acal(\cgot_+,\cgot_-)}K^\beta\cdot \sigma\lambda
$$ 
corresponds to the subset of the critical points of $\Phi$ that intersect 
$\Phi^{-1}(\xi)$ when $\xi\in\overline{\cgot_+}\cap\overline{\cgot_+}$.  

\begin{rem}
When $\lambda$ is a regular element of $\tgot^*$, all  polytopes 
${\conv}(W^\beta\cdot\sigma\lambda)$ are of codimension $1$. When $\lambda$ 
is ``generic'' (see Def. \ref{def:lambda-gene}), the set 
$\Acal(\cgot_+,\cgot_-)$ is reduced to one element.
\end{rem}

The multiplicity function $\mm^\lambda:\Lambda^*\times \Nbb^*\to\Nbb$ is 
invariant under the action of the Weyl group: 
$\mm^\lambda(\sigma\mu,k)=\mm^\lambda(\mu,k)$ for every 
$\sigma\in W$. The set of connected component of regular 
values of $\Phi$ is also invariant under the action of $W$. 

So for the rest of this section we  restrict our 
attention to case where $\cgot_+$ and $\cgot_-$ are 
separated by an hyperplane orthogonal to a fundamental weight 
$\beta=\varpi_i$: the vector $\beta=\varpi_i$ is pointing out 
of $\cgot_-$. Consider $\sigma\in \Acal(\cgot_+,\cgot_-)$ and let  
$K^i\cdot \sigma\lambda$ be the corresponding connected component 
of $(K\cdot\lambda)^\beta$ (here $K^i$ denote the stabilizer of 
$\varpi_i$ in $K$). The tangent space of $K\cdot\lambda$ at $\sigma\lambda$ 
is the following $K^{\sigma\lambda}$-module 
\begin{equation}
  \label{eq:tangent-s-lambda}
  \T_{\sigma\lambda}(K\cdot\lambda)=\sum_{(\alpha,\sigma\lambda)>0}\kgot_\alpha
\end{equation}
  where $\kgot_\alpha\subset\kgot\otimes\Cbb$ is the one-dimensional 
complex subspace associated to the weight $\alpha\in\Rgot$. In the same way, 
the tangent space of $K^i\cdot\sigma\lambda$ at $\sigma\lambda$ 
is the $K^i\cap K^{\sigma\lambda}$-module defined by 
\begin{equation}
  \label{eq:tangent-s-lambda-i}
  \T_{\sigma\lambda}(K^i\cdot\sigma\lambda)=
\sum_{\stackrel{(\alpha,\sigma\lambda)>0}{(\alpha,\varpi_i)= 0}}
\kgot_\alpha
\end{equation}
Finally the normal bundle  of $K^i\cdot\sigma\lambda\simeq 
K^i/(K^i\cap K^{\sigma\lambda})$ in $K\cdot\lambda$ is $\Ncal_{\sigma,i}=$ \break 
$K^i\times_{K^i\cap K^{\sigma\lambda}}N_{\sigma,i}$ where 
\begin{equation}
  \label{eq:normal-sigma}
  N_{\sigma,i}=\sum_{\stackrel{(\alpha,\sigma\lambda)>0}{(\alpha,\varpi_i)\neq 0}}
\kgot_\alpha
\end{equation}

For an element $\mu\in\tgot^*$, we have $\mu=\sum_{i=1}^{\dim T}[\mu]_k\, \alpha_k$ 
where the $[\mu]_k\in\Rbb$ are  defined by the relation $[\mu]_k= 
2\frac{(\varpi_k,\mu)}{\vert\alpha_k\vert^2}$. 
When $\mu\in\Lambda^*_R$ the coefficients
$[\mu]_k$ are all integers.

\begin{defi}\label{def:s-sigma-pm}
  For $\sigma\in \Acal(\cgot_+,\cgot_-)$ we define the positive integers 
$$
s^\pm_{\sigma,i}=\pm\sum_{\stackrel{(\alpha,\sigma\lambda)>0}{\pm(\alpha,\varpi_i)> 0}}
[\alpha]_i .
$$
Note that $s^+_{\sigma,i}+s^-_{\sigma,i}$ is larger than half of the codimension 
of $K^i\cdot\sigma\lambda$ in $K\cdot\lambda$.
\end{defi}

\begin{theo}\label{theo:mult-coadjoint}
Let $\cgot_+$ and $\cgot_-$ be two adjacent connected component of 
regular values of $\Phi:K\cdot\lambda\to\tgot^*$ separated by an hyperplane 
orthogonal to a fundamental weight $\varpi_i$: we denote $r_i$ the commum value 
$[\xi]_i$ for all $\xi$ in this hyperplane. Let
$\mm_{\cgot_\pm}^\lambda : \Lambda^* \times \Zbb \longrightarrow  \Zbb$ be the 
corresponding periodic polynomials: they are supported on the sub-lattice 
$\Xi_\lambda:=\{(\mu,k)\, \vert\, \mu\in k\lambda +\Lambda^*_R\}$.

For all $(\mu,k)\in\Xi_\lambda$, we have 
$\mm^\lambda_{\cgot_+}(\mu,k)=\mm_{\cgot_-}^\lambda(\mu,k)$ when
\begin{equation}
  \label{eq:m-egale-lambda}
 - s^-_i < [\mu]_i- k r_i< s^+_i. 
\end{equation}
Here the positive integer $s^\pm_i$ are defined by 
\begin{equation}
  \label{eq:n-pm-lambda}
  s^\pm_i=\inf_{\sigma\in \Acal(\cgot_+,\cgot_-)} s^\pm_{\sigma,i}.
\end{equation}  

When $\Acal(\cgot_+,\cgot_-)$ is reduced to one element $\sigma$, for example 
if $\lambda$ is ``generic'', the 
integer $s^+ +s^-$ is larger than half of the codimension 
of $K^i\cdot\sigma\lambda$ in $K\cdot\lambda$. 
\end{theo}

Another way to express the result of Theorem 
\ref{theo:mult-coadjoint} is to introduce 
like in \cite{Szenes-Vergne02} the convex polytope 
\begin{equation}
  \label{eq:carre-R}
  \Box(\cgot_+,\cgot_-)=\bigcap_{\sigma\in\Acal(\cgot_+,\cgot_-)}
\left(\sum_{(\alpha,\sigma\lambda)>0}[0,1[\,\alpha\right).
\end{equation}

Let $\p$ be the hyperplane which separates $\cgot_+$ and $\cgot_-$. Equation 
(\ref{eq:m-egale-lambda}) is equivalent to saying that 
\begin{equation}
  \label{eq:m-egale-carre}
 \mm^\lambda_{\cgot_+}(\mu,k)=\mm_{\cgot_-}^\lambda(\mu,k)\quad {\rm if}\quad 
\mu\in k\p +\Box(\cgot_+,\cgot_-).
\end{equation}

\bigskip

\begin{coro}
  Let $\cgot$ be a connected component of 
regular values of $\Phi$ which is bording a facet of the polytope 
$\Phi(K\cdot\lambda)$ orthogonal to the fundamental weight $\varpi_i$: 
the facet is ${\conv}(W^i\cdot\sigma\lambda)$ 
for a unique $\sigma\in W/W^i$. We suppose that $\varpi_i$ is pointing out of 
$\cgot$. We denote $r_i$ the commum value 
$[\xi]_i$ for all $\xi$ in the facet. For all $(\mu,k)\in\Xi_\lambda$, 
we have $\mm^\lambda_{\cgot}(\mu,k)=0$ when 
\begin{equation}
  \label{eq:m-egale-lambda-bord}
 - s^-_{\sigma,i} < [\mu]_i- k r_i< s^+_{\sigma,i}. 
\end{equation}
\end{coro}

\medskip

{\em Proof. } Theorem \ref{theo:mult-coadjoint} is a direct consequence of 
Theorem \ref{th:th-principal-mult}. The main difference between them is 
the decomposition of the lattice supporting the periodic 
polynomials. In the former we use 
the decomposition $\Lambda^*=\Lambda^*_{\tgot_\p}\oplus\Lambda^*_{\tgot/\tgot_\p}$ 
associated to the choice of a subtorus $T/T_\p$. Here, 
since $\mm^\lambda_{\cgot_\pm}$ is supported on $\lambda + \Lambda^*_R$, we use 
the decomposition $\Lambda^*_R=\Zbb\alpha_i\oplus\sum_{k\neq i}\Zbb\alpha_k$.

We start like after Proposition \ref{RR-0-diff}: $\mm^\lambda_{\cgot_+}(\mu,k)-
\mm^\lambda_{\cgot_-}(\mu,k)$ is equal to the $\mu$-mutiplicity in 
$\sum_{\sigma\in\Acal(\cgot_+,\cgot_-)}A^-_\sigma -A^+_\sigma$ where
\begin{equation}
  \label{eq:mu-lambda-pm}
 A^\pm_\sigma= 
RR_0^{\xi}\left(K^i\cdot\sigma\lambda,\Cbb_{[\lambda]}^{\otimes k}\otimes
\left[\wedge_{\Cbb}^{\bullet}\overline{\Ncal_{\sigma,i}}\right]^{-1}_{\mp\varpi_i}\right)
\end{equation}

Here $\xi$ belongs to the relative interior of 
$\overline{\cgot_+}\cap\overline{\cgot_+}$, the line bundle 
$\Cbb_{[\lambda]}^{\otimes k}$ is equal to 
$K^i\times_{K^i\cap K^{\sigma\lambda}}\Cbb_{k\sigma\lambda}$ and 
$\left[\wedge_{\Cbb}^{\bullet}
\overline{\Ncal_{\sigma,i}}\right]^{-1}_{\pm\varpi_i}$ 
corresponds to $(-1)^{{\rm rk}_\Cbb(N_{\sigma,i}^\pm)}$ times 
$$
K^i\times_{K^i\cap K^{\sigma\lambda}}\left(\det(N_{\sigma,i}^\pm)\otimes 
S^\bullet((N_{\sigma,i}\otimes\Cbb)^\pm)\right),
$$
with 
$$
N_{\sigma,i}^\pm=\sum_{\stackrel{(\alpha,\sigma\lambda)>0}{\pm(\alpha,\varpi_i)>0}}
\kgot_\alpha\ ,
$$
and 
$$
(N_{\sigma,i}\otimes\Cbb)^\pm=
\sum_{\stackrel{(\alpha,\sigma\lambda)\neq 0}{\pm(\alpha,\varpi_i)>0}}
\kgot_\alpha.
$$

Now we can apply Remark \ref{lem:discrete} with the subgroup $H\subset T$ equal to  
the center $Z(K^i)$ of $K^i$: an element $\gamma\in\Lambda^*$ belong to 
$\sum_{k\neq i}\Zbb\alpha_k$ if and only if $t^\gamma=1$ for all $t\in Z(K^i)$.

The group $Z(K^i)$ acts trivially on the manifolds $K^i\cdot\sigma\lambda$, 
and with the characters associated to the weights
$$
k\sigma\lambda+\sum_{\stackrel{(\alpha,\sigma\lambda)>0}
{(\alpha,\varpi_i)>0}}\alpha+\delta
\quad {\rm with}\quad (\delta,\varpi_i)\geq 0
$$
on the bundle $\Cbb_{[\lambda]}^{\otimes k}\otimes
\left[\wedge_{\Cbb}^{\bullet}\overline{\Ncal_{\sigma,i}}\right]^{-1}_{\varpi_i}$, and 
with the characters associated to the weights
$$
k\sigma\lambda+\sum_{\stackrel{(\alpha,\sigma\lambda)>0}
{(\alpha,\varpi_i)<0}}\alpha+\delta
\quad {\rm with}\quad (\delta,\varpi_i)\leq 0
$$
on the bundle $\Cbb_{[\lambda]}^{\otimes k}\otimes
\left[\wedge_{\Cbb}^{\bullet}\overline{\Ncal_{\sigma,i}}\right]^{-1}_{-\varpi_i}$. 
Now the 
$\mu$-multiplicity in $A^\pm_\sigma$ is not equal to $0$ only if 
\begin{equation}
  \label{eq:condition-rem-F}
 k\sigma\lambda+\sum_{\stackrel{(\alpha,\sigma\lambda)>0}{\pm(\alpha,\varpi_i)>0}}
\alpha+\delta
-\mu \in \sum_{k\neq i}\Zbb\alpha_k \quad {\rm with}\quad \pm(\delta,\varpi_i)\geq 0.
\end{equation}
Condition (\ref{eq:condition-rem-F}) implies that 
$[\mu]_i\geq k[\sigma\lambda]_i+s^+_{\sigma,i}$ or $[\mu]_i\leq k[\sigma\lambda]_i- 
s^-_{\sigma,i}$. Finally we have prove that  $\mm^\lambda_{\cgot_+}(\mu,k)=
\mm^\lambda_{\cgot_-}(\mu,k)$ if 
$$
-s^-_{\sigma,i}<[\mu]_i - k[\sigma\lambda]_i<s^+_{\sigma,i}
$$
for all $\sigma\in\Acal(\cgot_+,\cgot_-)$. $\Box$


\subsection{The case of $\hbox{\rm SU}(n)$}\label{subsec:su-n}

Let $T$ be the maximal torus of $\hbox{\rm SU}(n)$ consisting of the diagonal 
matrices. The dual $\tgot^*$ can be identified with the subspace 
$x_1 +\cdots +x_n=0$ of $\Rbb^n$. The roots are 
$\Rgot=\{e_i-e_j\, \vert 1\leq i\neq j\leq n\}$ and we will choose the positives 
ones to be $\Rgot^+=\{e_i-e_j\, \vert \, 1\leq i<j\leq n\}$. The simple roots 
are then $\alpha_i=e_i -e_{i+1}$, for $1\leq i\leq n-1$, and for these simple roots, 
the fundamental weights are 
\begin{equation}
  \label{eq:omega-su-n}
  \omega_k=\frac{1}{n}(\underbrace{n-k,n-k,\cdots,n-k}_{k\ {\rm times}}, 
\underbrace{-k,-k,\cdots,-k}_{n-k\ {\rm times}}),\quad 1\leq k\leq n-1.
\end{equation}

Consider now the coadjoint orbit $O_\lambda$ for $\lambda\in\tgot^*$. Let 
$\Phi:O_\lambda\to\tgot^*$ the moment map associated to the Hamiltonian action 
of $T$ on $O_\lambda$. The center of $\hbox{\rm SU}(n)$, that we denote $\Zbb_n$ 
corresponds to the set of matrices $zI$ with $z^n=1$. Recall the following 
well-known fact.

\begin{lem}
  Let $\xi$ be a regular value of $\Phi:O_\lambda\to\tgot^*$. Then for every 
$m\in \Phi^{-1}(\xi)$ the stabilizer subgroup $T_m:=\{t\in T\,\vert\, t\cdot m=m\}$ 
is equal to $\Zbb_n$.
\end{lem}

{\em Proof.} Since $\xi$ is a regular value, we know that $T_m$ is finite for every 
$m\in \Phi^{-1}(\xi)$. The dual of the Lie algebra $\sugot(n)$ decomposes as 
$\sugot(n)^*=\tgot^*\oplus\sum_{\alpha\in\Rgot^+}\sugot(n)^*_\alpha$ where 
$\sugot(n)^*_\alpha\simeq\Cbb_{-\alpha}$ as $T$-module. For $m\in\Phi^{-1}(\xi)$, we 
have $m=m_0 + \sum_{\alpha\in\Rgot^+}m_\alpha$ with $m_\alpha\in \sugot(n)^*_\alpha$, and 
then $T_m=\cap_{m_\alpha\neq 0}\ker(t\mapsto t^\alpha)$. So the 
lattice $\Lambda^*_m$ generated by the set 
$\{\alpha\in\Rgot^+\,\vert\, m_\alpha\neq 0\}$ 
is a subgroup of $\Lambda^*_R$ with $\Lambda^*_R/\Lambda^*_m$  finite. 
We have to show that $\Lambda^*_m=\Lambda^*_R$. For this purpose we introduce 
the following equivalence relation on $\{1,\ldots, n\}$: 
$$
i\sim j\Longleftrightarrow  e_i-e_j\in \Lambda^*_m.
$$
Suppose that $\{1,\ldots, n\}/\sim$ is not reduced to a point: 
let $C_1$ and $C_2$ be two distinct equivalent classes and let 
$\beta=(\beta_1,\ldots,\beta_n)$ be the element 
of $\tgot^*$ defined by: $\beta_i=\frac{1}{\vert C_1\vert}$ if 
$i\in C_1$, $\beta_i=\frac{-1}{\vert C_2\vert}$ if $i\in C_2$, 
and $\beta_i=0$ in the other cases. We see then that 
$(\beta,\alpha)=0$ for all $\alpha\in \Lambda^*_m$: it is 
in contradiction with the fact that $\Lambda^*_R/\Lambda^*_m$ 
is finite. We have proved that $e_i-e_j\in \Lambda^*_m$ 
for all $i,j\in\{1,\ldots, n\}$. $\Box$

\bigskip 

Suppose now that $\lambda$ is a positive weight, and let $\cgot$ 
a connected component of regular values of $\Phi:O_\lambda\to\tgot^*$. 
We know that  the corresponding periodic polynomial
$\mm_{\cgot}^\lambda : \Lambda^* \times \Zbb \longrightarrow  \Zbb$ 
is supported on the sub-lattice $\Xi_\lambda:=\{(\mu,k)\, \vert\, 
\mu\in k\lambda +\Lambda^*_R\}$.

\begin{coro}
  The map $\mm_{\cgot}^\lambda : \Xi_\lambda 
\longrightarrow  \Zbb$ is a polynomial of degree 
$\frac{(n-1)(n-2)}{2}-d_\lambda$, where $d_\lambda$ is the number 
of positive roots orthogonal to $\lambda$.
.
\end{coro}

{\em Proof.} Take $\xi\in\cgot$. Following Proposition \ref{prop:m-poly}, 
the periodic-polynomial $\mm^\lambda_{\cgot}$ is defined by 
$\mm^\lambda_{\cgot}(\mu,k)=RR((O_\lambda)_\xi, \Lcal^k_{\xi,\mu})$ 
for all $(\mu,k)\in\Xi_\lambda$. Here $(O_\lambda)_\xi=\Phi^{-1}(\xi)/T$ is a 
{\em smooth} manifold, and the line bundle 
$\Lcal^k_{\xi,\mu}=(L^{\otimes k}\vert_{\Phi^{-1}(\xi)}\otimes\Cbb_{-\mu})/T$
is also {\em smooth} since the center $\Zbb_n$ acts trivially on 
$L^{\otimes k}\vert_{\Phi^{-1}(\xi)}\otimes\Cbb_{-\mu}$. Now the 
Atiyah-Singer integral formula  for the Riemann-Roch number 
$RR((O_\lambda)_\xi, \Lcal^k_{\xi,\mu})$ shows that 
$\mm^\lambda_{\cgot}$ is a polynomial 
of degree $\frac{\dim (O_\lambda)_\xi}{2} =\frac{\dim O_\lambda}{2} -(n-1)=
\frac{(n-1)(n-2)}{2}-d_\lambda$. $\Box$ 

\bigskip

Now we rewrite Theorem \ref{theo:mult-coadjoint} for the group $\SU(n)$. 
Let $\lambda=(\lambda_1\geq\cdots\geq\lambda_n)$ be a  positive weight and let 
$\cgot_+$ and $\cgot_-$ be two adjacent connected components of 
regular values of $\Phi:O_\lambda\to\tgot^*$ separated by an hyperplane 
orthogonal to a fundamental weight $\varpi_i$: the vector $\varpi_i$ is pointing out 
of $\cgot_-$. Let $q(\xi)=(\varpi_i,\xi) -r_i$ be the defining equation 
of this hyperplane.

The conditions $(e_k-e_l,\sigma\lambda)>0$ and 
$(e_k-e_l,\varpi_i)>0$ are respectively  
equivalent to $\lambda_{\sigma(k)}>\lambda_{\sigma(l)}$ and $k\leq i<
l$. For $\SU(n)$, the number $[\alpha]_i$ is equal to $0,1$ 
or $-1$ for any roots $\alpha$ and any 
$i=1,\cdots,n-1$. Hence for every $\sigma\in\Acal(\cgot_+,\cgot_-)$, 
the integers $s_{\sigma,i}^-,s_{\sigma,i}^+\geq 0$ defined in Def. 
(\ref{def:s-sigma-pm}) are equal to 
\begin{equation}
  \label{eq:s-plus-su-n}
  s_{\sigma,i}^+={\rm rk}_\Cbb(N_{\sigma,i}^{+})=
\sharp\{k\leq i< l\quad {\rm such\ that}\quad 
\lambda_{\sigma(k)}>\lambda_{\sigma(l)} \}, 
\end{equation}
\begin{equation}
  \label{eq:s-moins-su-n}
  s_{\sigma,i}^-={\rm rk}_\Cbb( N_{\sigma,i}^{-})=
\sharp\{k\leq i< l\quad {\rm such\ that}\quad 
\lambda_{\sigma(k)}<\lambda_{\sigma(l)} \}, 
\end{equation}
and the sum $s_{\sigma,i}^+ + s_{\sigma,i}^-$ is equal to half of 
the codimension of $K^i\cdot\sigma\lambda$ in $K\cdot\lambda$, 
that is $s_{\sigma,i}^+ + s_{\sigma,i}^-=i(n-i)-
\dim(K^{\sigma\lambda}/K^i\cap K^{\sigma\lambda})/2$.

Now we precise the results of \cite{BGR03}.

\begin{theo}\label{th:SU-N} 
The polynomial $\mm_{\cgot_-}^\lambda-\mm_{\cgot_+}^\lambda: 
\Xi_\lambda\to\Zbb$ is divible by the linear factors 
$$
(q-s_i^-+1),(q-s_i^-+2),\ldots,q,\ldots,(q+s_i^+ -2),(q+s_i^+ -1), 
$$ 
where  $q$ is the defining equation of the hyperplane separating 
$\cgot_\pm$ and $s_i^\pm=$\break 
$\inf_{\sigma\in\Acal(\cgot_+,\cgot_-)}s_{\sigma,i}^\pm$.  
Moreover the linear factors $(q-s_i^-)$ and $(q-s_i^-)$ {\em do not divide} 
$\mm_{\cgot_-}^\lambda-\mm_{\cgot_+}^\lambda$. 
\end{theo}

{\em Proof.} The first part is the translation of 
Theorem \ref{theo:mult-coadjoint}. 
We have just to prove that the linear factors $(q-s_i^-)$ and $(q-s_i^-)$ 
do not divide $\mm_{\cgot_-}^\lambda-\mm_{\cgot_+}^\lambda$. This point is a  
direct application of Proposition \ref{prop:inegalite-optimale}. 
The only fact we use here is that 
${\rm rk}_\Cbb( N_{\sigma,i}^{\pm})=s_{\sigma,i}^\pm$. So the number 
${\rm rk}_\Cbb( N_{\sigma,i}^{\pm})$ is constant for all 
$\sigma\in\Acal(\cgot_+,\cgot_-)$ for which $s_{\sigma,i}^\pm=s_{i}^\pm$. $\Box$

\medskip

We rewrite now Theorem \ref{th:SU-N} in the particular case where 
$\Acal(\cgot_+,\cgot_-)$ contains just one element: it happens 
when  $\lambda$ is a ``generic'' positive weight 
(see Definition \ref {def:lambda-gene}), 
or when $\cgot_+$ does not intersect $\Phi(O_\lambda)$. Here a positive weight
$\lambda=(\lambda_1\geq\cdots\geq\lambda_n)$ is ``generic'' if for every couple 
of permutations $\sigma,\sigma'$ and any $k=1,\cdots,n-1$, we have
$$
\sum_{i=1}^k\lambda_{\sigma(i)}\neq\sum_{i=1}^k\lambda_{\sigma'(i)}
$$
when $(\lambda_{\sigma(1)},\cdots,\lambda_{\sigma(n)})\notin \Sgot_k\times\Sgot_{n-k}
(\lambda_{\sigma'(1)},\cdots,\lambda_{\sigma'(n)})$.

\begin{coro} Let $\lambda$ be a {\em regular} weight. Let 
$\cgot_+$ and $\cgot_-$ be two adjacent connected components of 
regular values of $\Phi:O_\lambda\to\tgot^*$ and suppose that 
$\Acal(\cgot_+,\cgot_-)$ contains just one element $\sigma$. Then 
the polynomial $\mm_{\cgot_-}^\lambda-\mm_{\cgot_+}^\lambda: 
\Xi_\lambda\to\Zbb$ is divible by the $i(n-i)$ linear factors 
$$
(q-s_i^-+1),(q-s_i^-+2),\ldots,q,\ldots,(q+s_i^+ -2),(q+s_i^+ -1), 
$$ 
where  $s_i^\pm=s_{\sigma,i}^\pm$ are defined by (\ref{eq:s-plus-su-n}) and 
(\ref{eq:s-moins-su-n}). Moreover the linear factors $(q-s_i^-)$ and 
$(q-s_i^+)$ {\em do not divide} 
$\mm_{\cgot_-}^\lambda-\mm_{\cgot_+}^\lambda$.   
\end{coro}


\section{Vector partition functions}\label{sec:cas-C-d}

Let $T$ be a torus with Lie algebra $\tgot$ and let $\Lambda^*\subset\tgot^*$ be the 
weight lattice. Let $R=\{\alpha_1,\ldots,\alpha_d\}$ be a subset of not 
necessarily distint elements of $\Lambda^*$ which lie enterely in an open halfspace 
of $\tgot^*$. We associate with the collection $R$ a function 
$$
N_R: \Lambda^*\longrightarrow \Nbb
$$
called the {\em vector partition function} associated to $R$. By definition, for 
a weight $\mu$, the value $N_R(\mu)$ is the number of solutions of the
equation 
\begin{equation}
  \label{eq:partition-function}
  \sum_{j=1}^d k_j\alpha_j=\mu,\quad k_j\in\Zbb^{\geq 0},\quad j=1,\ldots,d.
\end{equation}

 Let $C(R)\subset\tgot^*$ be the closed convex cone generated by the elements of 
$R$, and denote by $\Lambda^*_R\subset\Lambda^*$ the sublattice generated 
by $R$. Obviously, $N_R(\mu)$ vanishes if $\mu$ does not 
belong to $C(R)\cap\Lambda^*_R$.

Suppose now that $R$ generates the vector space $\tgot^*$. 
Following \cite{Szenes-Vergne02}, we will call a vector {\em singular} with 
respect to $R$ if it is in a cone $C(\nu)$ generated by a subset 
$\nu\in R$ of cardinality strictly less than $dim T$. The connected components 
of $\tgot^*-\{singular\  vectors\}$ are called {\em conic chambers}. 
The periodic polynomial behavior of $N_R$ on closures of conic
chambers of the cone $C(R)$ is proved in \cite{Sturmfels95}. We have
the following refinement due to Szenes and Vergne 
\cite{Szenes-Vergne02}. Let us introduce the convex polytope 
\begin{equation}
  \label{eq:box-phi}
  \Box(\Phi)=\sum_{j=1}^d[0,1]\alpha_j.
\end{equation}
We remark that $\cgot-\Box(\Phi)$ is a neighborhood of
$\overline{\cgot}$ for any conic chamber $\cgot$ of the cone $C(R)$. 
We have the following qualitative result. 

\begin{theo}[\cite{Szenes-Vergne02}]\label{prop:N-conic-chamber-SV}
  Let $\cgot$ be a conic chamber of the cone $C(R)$. There exists a 
periodic polynomial $P_\cgot$ on $\Lambda^*$ such that for each 
$\mu\in \cgot-\Box(\Phi)$, we have 
$$
N_R(\mu)=P_\cgot(\mu).
$$
\end{theo}

In Section \ref{proof:N-conic-chamber-SV} we will give another proof of
Theorem \ref{prop:N-conic-chamber-SV}. Let $\cgot_\pm\subset\tgot^*$ be 
two adjacent conic chambers. The aim of this Section 
is to give a formula for the periodic polynomial  
$P_{\cgot_+}-P_{\cgot_-}$.

Let $\p\subset\tgot^*$ be the hyperplane that separates $\cgot_+$ 
and $\cgot_-$. Let $\beta\in\tgot$ be such that 
$\p=\{\xi\in\tgot^* \ \vert\ \langle\xi,\beta\rangle=0\}$ and 
$\cgot_\pm\subset \{\xi\in\tgot^* \ \vert\ \pm\langle\xi,\beta\rangle>0\}$. 
Note that the vector space $\p$ is generated by $R\cap\p$. We will now 
polarize the elements of $R$ that are outside $\p$.  We define  
\begin{equation}
  \label{eq:R-delta-prime}
  R'=\{\epsilon_j\alpha_j\ \vert\ \langle\alpha_j,\beta\rangle\neq 0
\ {\rm and}\ \epsilon_j={\rm sign}\ 
\langle\alpha_j,\beta\rangle\},
\end{equation} 
\begin{equation}
  \label{eq:delta-pm}
  \delta^\pm=\sum_{\pm\langle\alpha_j,\beta\rangle>0}\alpha_j\ ,
\end{equation}
and 
\begin{equation}
  \label{eq:r-pm}
  r^\pm= \sharp\{j\ \vert \ \pm\langle\alpha_j,\beta\rangle>0\}.
\end{equation}

We now look at the vector space $\p$ equipped with the subset 
$R\cap\p\subset\Lambda^*\cap\p$ which lie enterely in 
an open halfspace: let $N_{R\cap\p}:\Lambda^*\cap\p\to\Nbb$ be 
the corresponding vector partition function. Let $\cgot'$ be 
the cone in $\p$ which is equal to the relative interior of 
$\overline{\cgot_+}\cap\overline{\cgot_-}$. It is easy to see 
that $\cgot'$ is a conic chamber in $\p$ with respect to $R\cap\p$. 
Following Proposition \ref{prop:N-conic-chamber-SV} there exists 
a periodic polynomial $P_{\cgot'}$ on $\Lambda^*\cap\p$ such 
that for each $\mu\in\overline{\cgot'}\cap\Lambda^*$, we have 
$$
N_{R\cap\p}(\gamma)=P_{\cgot'}(\gamma).
$$

Let  $N_{R'}:\Lambda^*\to\Nbb$ be the vector partition 
function associated to the polarized set of weight 
$R'$ (see (\ref{eq:R-delta-prime})). The main result of 
this Section is the following  

\begin{theo}\label{theo:cas-C-d}
The periodic polynomial $P_{\cgot_+}-P_{\cgot_-}:\Lambda^*\to \Zbb$
satisfies 
\begin{equation}
  \label{eq:P-pm-C-d}
P_{\cgot_+}(\mu)-P_{\cgot_-}(\mu)=\sum_{\gamma\in\Lambda^*\cap\p}
D(\mu-\gamma)P_{\cgot'}(\gamma),\quad \mu\in\Lambda^*,
\end{equation}
where $D:\Lambda^*\to \Zbb$ is defined by 
$$
D(\mu)=(-1)^{r^-}N_{R'}(\mu+\delta^-)-(-1)^{r^+}N_{R'}(-\mu-\delta^+).
$$ 
\end{theo}

\bigskip

The proof of Theorem \ref{theo:cas-C-d} will be given in Section 
\ref{proof-theo-cas-C-d}. 

\medskip

\begin{coro}
$P_{\cgot_+}(\mu)=P_{\cgot_-}(\mu)$ for all the weights 
$\mu\in\Lambda^*$ satisfying the condition 
$$
- \langle\delta^+,\beta\rangle<\langle\mu,\beta\rangle
< -\langle\delta^-,\beta\rangle.
$$
The former ineqalities are optimal since 
$$
\left(P_{\cgot_+}-P_{\cgot_-}\right)(-\delta^- +\gamma)=
(-1)^{r^-}P_{\cgot'}(\gamma)
$$
and 
$$
\left(P_{\cgot_+}-P_{\cgot_-}\right)(-\delta^+ +\gamma)=
(-1)^{1+r^+}P_{\cgot'}(\gamma)
$$
for all $\gamma\in\Lambda^*\cap\p$.
\end{coro}

\medskip

{\em Proof.} In (\ref{eq:P-pm-C-d}), the term 
$D(\mu-\gamma)P_{\cgot'}(\gamma)$ does not
vanish only if $\mu-\gamma\in -\delta^- + C(R')$ or 
$-(\mu-\gamma)\in \delta^+ + C(R')$ for some $\gamma\in C(R\cap\p)$. 
These two conditions impose respectively that 
$\langle\mu,\beta\rangle\geq -\langle\delta^-,\beta\rangle$ and 
$\langle\mu,\beta\rangle\leq -\langle\delta^+,\beta\rangle$. If one 
take $\mu=-\delta^- +\gamma$ with $\gamma\in\Lambda^*\cap\p$, 
(\ref{eq:P-pm-C-d}) becomes $(P_{\cgot_+}-P_{\cgot_-})(-\delta^- +\gamma)=
\sum_{\gamma'\in\Lambda^*\cap\p}D(-\delta^- +\gamma-\gamma')
P_{\cgot'}(\gamma')$ with 
$$
D(-\delta^- +\gamma-\gamma')=
(-1)^{r^-}N_{R'}(\gamma  -\gamma')-(-1)^{r^+}N_{R'}(\delta^- -\delta^+ 
-\gamma+\gamma').
$$
Since the cone $C(R')$ intersects $\p$ only at $\{0\}$, 
$N_{R'}(\gamma-\gamma')=0$ if $\gamma\neq\gamma'$. Since 
$\langle\delta^- -\delta^+,\beta\rangle< 0$ we always 
have $N_{R'}(\delta^- -\delta^+ -\gamma+\gamma')=0$. We get  
finally that $(P_{\cgot_+}-P_{\cgot_-})(-\delta^- +\gamma)=
(-1)^{r^-}P_{\cgot'}(\gamma)$. One can show in the same way that 
$(P_{\cgot_+} -P_{\cgot_-})(-\delta^+ +\gamma)=
-(-1)^{r^+}P_{\cgot'}(\gamma)$. $\Box$


\subsection{Quantization of $\Cbb^d$}

We consider the complex vector space $\Cbb^d$ equipped with the canonical 
symplectic form $\Omega=\frac{i}{2}\sum_{i=1}^d dz_j\wedge d\overline{z}_j$. 
The standard complex struture $J$ on $\Cbb^d$ is compatible with $\Omega$. 
Let $T$ be a torus, let $\alpha_j\in\tgot^*$, $j=1,\ldots,d$ be weights 
of $T$, and let $T$ acts on $\Cbb^d$ as
\begin{equation}
  \label{eq:action-T-C}
  t\cdot(z_1,\ldots,z_d)=(t^{-\alpha_1}z_1,\ldots,t^{-\alpha_d}z_d).
\end{equation}
The action of $T$ preserve the symplectic form $\Omega$ and the moment map 
associated with this action is 
\begin{equation}
  \label{eq:moment-C}
  \Phi(z)= \frac{1}{2}\sum_{i=1}^d \vert z_j\vert^2 \alpha_j.
\end{equation}

The pre-quantization data $(L,\langle,\rangle,\nabla)$ 
on the Hamiltonian $T$-manifold $(\Cbb^d,\Omega,\Phi)$ is a 
trivial line bundle $L$ with a trivial action of $T$ equipped 
with the Hermitian structure $\langle s,s'\rangle_z= 
e^{\frac{-\vert z\vert^2}{2}}s\overline{s'}$ and the 
Hermitian connexion $\nabla=d-\beta$ where $\beta=
\frac{1}{2}\sum_{i=1}^d \overline{z}_j dz_j$.

The quantization of the Hamiltonian $T$-manifold  $(\Cbb^d,\Omega)$, 
that we denote $\Qcal^T(\Cbb^d)$, is the {\em Bargman space} of 
entire holomorphic functions on $\Cbb^d$ which are $\Lcal^2$ 
integrable with respect to the Gaussian measure 
$e^{\frac{-\vert z\vert^2}{2}}\Omega^d$.

We suppose now that the set of weights
$R=\{\alpha_1,\ldots,\alpha_d\}$ is {\em polarized} by $\eta\in\tgot$, 
which means that $\langle\alpha_j,\eta\rangle >0$ for all $j$. 
The $T$-representation $\Qcal^T(\Cbb^d)$ is then admissible 
and we have the following equality in $R^{-\infty}(T)$:
\begin{equation}
  \label{eq:Q-T-C}
 \Qcal^T(\Cbb^d)=\sum_{\mu\in\Lambda^*} N_R(\mu)\, \Cbb_\mu,
\end{equation}
where $N_R: \Lambda^*\to \Nbb$ is the vector partition function 
associated to $R$. In other words, the generalized character 
of $\Qcal^T(\Cbb^d)$ coincides with the generalized character 
of the symmetric algebra $S^\bullet(\overline{\Cbb^d})$, where
$\overline{\Cbb^d}$ means $\Cbb^d$ with the opposite complex structure. 

\medskip

For the remaining part of Section \ref{sec:cas-C-d}, we assume that
the set of weights $R=\{\alpha_1,\ldots,\alpha_d\}$ is {\em
  polarized}, and {\em generate} the vector space $\tgot^*$. The first
assumption is equivalent to the fact that the moment map 
$\Phi:\Cbb^d\to\tgot^*$ is {\em proper}, and the second assumption 
is equivalent to the fact that the generic stabiliser of $T$ on
$\Cbb^d$ is {\em finite}. Notice that the vectors of $\tgot^*$ which 
are {\em singular} with respect to $R$ correspond to the {\em singular 
values} of $\Phi$. 

In the next section we will show that $\Qcal^T(\Cbb^d)$, viewed as an
element of $R^{-\infty}(T)$, can 
be realized as the index of transversally elliptic symbols on
$\Cbb^d$. After we will apply the techniques developped in Section 
\ref{sec:quantum-version}. The main difference here is
that we work with a {\em non-compact} manifold.


\subsection{Transversally elliptic symbols on $\Cbb^d$}\label{subsec:C-d}
Let $p:\T\Cbb^d\to\Cbb^d$ be the canonical projection. The Thom symbol 
$$
\Thom(\Cbb^d)\in\Gamma\left(\T\Cbb^d,\hom(p^*(\wedge^{\rm even}_\Cbb\T\Cbb^d), 
p^*(\wedge^{\rm odd}_\Cbb\T\Cbb^d))\right)
$$
is defined as follows. At $(z,v)\in\T\Cbb^d$, the Thom symbol $\Thom(\Cbb^d)$ is 
equal to the Clifford map 
$$
Cl(v):\wedge^{\rm even}_\Cbb\Cbb^d\longrightarrow \wedge^{\rm odd}_\Cbb\Cbb^d
$$
which is defined by $Cl(v)w=v\wedge w -c(v)w$. 
Here $c(v):\wedge^{\bullet}_\Cbb\Cbb^d\mapsto
\wedge^{\bullet-1}_\Cbb\Cbb^d$ denotes the contraction 
map relatively to the standard Hermitian structure on $\Cbb^d$. 
Obviously the symbol $\Thom(\Cbb^d)$ is not elliptic since its 
characteristic set is equal to the zero section in $\T\Cbb^d$ 
(hence is not compact).

Now we deform the symbol $\Thom(\Cbb^d)$ in order to obtain 
transversally elliptic symbols. Since $\Cbb^d$ can be realized 
as an open subset of a compact $T$-manifold 
we have a well defined index map
$$
\indice^T_{\Cbb^d}:\K_T(\T_T\Cbb^d)\longrightarrow R^{-\infty}(T).
$$ 

\begin{defi}\label{def:Thom-deforme}
  For any $\eta\in\tgot$, we define the symbol $\Thom^\eta(\Cbb^d)$ by 
$$
\Thom^\eta(\Cbb^d)(z,v)=\Thom(\Cbb^d)(z,v-\eta_{\Cbb^d}(z)), 
\quad (z,v)\in\T\Cbb^d,
$$
where $\eta_{\Cbb^d}$ is the vectors field on $\Cbb^d$ generated by $\eta$.
\end{defi}

The symbols $\Thom^\eta(\Cbb^d)$ were studied in \cite{pep4}. 
It is easy to see that  $\Thom^\eta(\Cbb^d)$ is tranversally 
elliptic if and only if the vector subspace 
$(\Cbb^d)^\eta$ is reduced to $\{0\}$, i.e. if 
$\langle\alpha_j,\eta\rangle\neq 0$ for all $j=1,\ldots,d$. 
We prove in Proposition 5.4. of \cite{pep4} that 
\begin{equation}
  \label{eq:indice-thom-eta}
 \indice^T_{\Cbb^d}\left(\Thom^\eta(\Cbb^d)\right)=
S^\bullet(\overline{\Cbb^d})\quad {\rm in}
\quad R^{-\infty}(T),
\end{equation}
when $\langle\alpha_j,\eta\rangle> 0$ for all $j=1,\ldots,d$.

In order to compute the multiplicities $N_R(\mu)$ of 
$\Qcal^T(\Cbb^d)$ we introduce the following tranversally 
elliptic symbols. Take a scalar product $b(\cdot,\cdot)$ on 
$\tgot^*$, and denote by $\xi\mapsto\xi^b,\tgot^*\simeq\tgot$ 
the induced isomorphism. For each $\xi\in\tgot^*$, the 
Hamiltonian vectors field of the function 
$\frac{-1}{2}\|\Phi-\xi\|^2_b$ is the vectors field 
$$
z\mapsto \Big( (\Phi(z)-\xi)^b\Big)_{\Cbb^d}(z).
$$
that we denote $\Hcal^b-\xi^b_{\Cbb^d}$. 

\begin{defi}\label{def:Thom-deforme-phi}
  For any $\xi\in\tgot^*$, and any scalar product $b(\cdot,\cdot)$ on 
$\tgot^*$, we define the symbol $\Thom_{\xi,b}(\Cbb^d)$ by 
$$
\Thom_{\xi,b}(\Cbb^d)(z,v)=\Thom(\Cbb^d)(z,v-(\Hcal^b-\xi^b_{\Cbb^d})(z)), 
\quad (z,v)\in\T\Cbb^d.
$$
\end{defi}

Let  $\Char(\Thom_{\xi,b}(\Cbb^d))\subset\T\Cbb^d$ be the characterictic set of 
$\Thom_{\xi,b}(\Cbb^d)$. We know that
$\Char(\Thom_{\xi,b}(\Cbb^d))\cap\T_T\Cbb^d$ 
is equal to the critical set $\crpt(\|\Phi-\xi\|^2_b)$ of the function 
$\|\Phi-\xi\|^2_b:\Cbb^d\to\Rbb$ (see Section \ref{ssection.thom}). 
A straightforward computation gives that $z\in\crpt(\|\Phi-\xi\|^2_b)$ if and only if 
\begin{equation}
  \label{eq:critique-C-d}
  b(\Phi(z)-\xi,\alpha_j)\, z_j=0 \quad {\rm for\ all}
\quad j=1,\ldots,d.
\end{equation}
The former relations implies in particular that 
$b(\Phi(z)-\xi,\Phi(z))=\frac{1}{2}\sum_j b(\Phi(z)-\xi,\alpha_j)\, 
\vert z_j\vert^2=0$. Hence $\|\Phi(z)\|^2_b=b(\Phi(z),\xi)$ which 
implies 
\begin{equation}
  \label{eq:phi-less-xi}
  \|\Phi(z)\|_b\leq \|\xi\|_b.
\end{equation}
Take now $\eta\in\tgot$ such that $\langle\alpha_j,\eta\rangle> 0$ for all $j$, 
and let $\eta_b\in\tgot^*$ such that $(\eta_b)^b=\eta$. We have then 
\begin{equation}
  \label{eq:C-eta}
C_{\eta} \|z\|^2\leq \langle\Phi(z),\eta\rangle=b(\Phi(z),\eta_b)\leq 
\|\Phi(z)\|_b\|\eta_b\|_b  
\end{equation}
where $C_\eta=\frac{1}{2}\inf_j \langle\alpha_j,\eta\rangle$, and 
$z\mapsto\|z\|^2$ is the usual hermitian form on $\Cbb^d$. 
With (\ref{eq:critique-C-d}) and (\ref{eq:C-eta}) 
we get the following 

\begin{lem}\label{lem:crit-C-d}
  The critical set $\crpt(\|\Phi-\xi\|^2_b)\subset\Cbb^d$ is contained in the ball
of radius 
$$
\frac{\|\xi\|_b \,\|\eta_b\|_b}{C_\eta},
$$ 
where $\eta\in\tgot$ is such that 
$C_\eta=\frac{1}{2}\inf_j \langle\alpha_j,\eta\rangle >0$.
\end{lem}

We have then proved that the symbols $\Thom_{\xi,b}(\Cbb^d)$ are transversally  
elliptic.

\begin{prop}\label{prop:indice-thom-xi}
The class of the transversally elliptic symbol $\Thom_{\xi,b}(\Cbb^d)$  
in $\K_T(\T_T\Cbb^d)$ does not depend of the data $\xi,b$, and is
equal to the class defined by $\Thom^\eta(\Cbb^d)$ where $\eta\in\tgot$
is choosen so that $\langle\alpha_j,\eta\rangle >0$ for all $j$. 
\end{prop}

{\em Proof}. After Lemma \ref{lem:crit-C-d}, we know that for any scalar product 
$b(\cdot,\cdot)$ on $\tgot^*$, the characteristic set of 
$\Thom_{0,b}(\Cbb^d)$ intersects  $\T_T\Cbb^d$ at $\{0\}$. If $b_0$ and $b_1$ 
are two scalar products on $\tgot^*$ we consider the 
family $b_t=t b_1 +(1-t)b_0,\ 0\leq t\leq 1$, of scalar products on $\tgot^*$. 
Hence $\Thom_{0,b_t}(\Cbb^d)$, $t\in[0,1]$, defines an 
homotopy of transversally elliptic symbols. We have proved that 
$\Thom_{0,b_0}(\Cbb^d)=\Thom_{0,b_1}(\Cbb^d)$ in $\K_T(\T_T\Cbb^d)$ 
for any $\xi\in\tgot^*$.

Fix now the scalar product $b$ and an element $\xi\in\tgot^*$. For any $t\in[0,1]$ the 
characteristic set of $\Thom_{t\xi,b}(\Cbb^d)$ intersects  $\T_T\Cbb^d$ in  
 the ball of radius 
$$
\frac{\|\xi\|_{b} \,\|\eta_b\|_{b}}{C_\eta}.
$$
Hence $\Thom_{t\xi,b}(\Cbb^d)$, $t\in[0,1]$, defines an homotopy of 
transversally elliptic symbols: $\Thom_{\xi,b}(\Cbb^d)=\Thom_{0,b}(\Cbb^d)$ 
in $\K_T(\T_T\Cbb^d)$. We have proved that the class of 
the transversally elliptic symbol $\Thom_{\xi,b}(\Cbb^d)$ in 
$\K_T(\T_T\Cbb^d)$ does not depend of the data $\xi,b$. 

Since the weights $\alpha_j$ lie enterely in an open halfspace 
of $\tgot^*$, there exists a scalar product $b_+(\cdot,\cdot)$ 
on $\tgot^*$ for which we have 
$$
b_+(\alpha_i,\alpha_j)>0
$$
for all $i,j=1,\ldots,d$. Let $\Hcal^{b_+}$ be the Hamiltonian vectors field of 
the function $\frac{-1}{2}\|\Phi\|^2_{b_+}$, and let $\eta_{\Cbb^d}$ be the 
vectors field on $\Cbb^d$ generated by  $\eta\in\tgot$ such 
that $\langle\alpha_j,\eta\rangle >0$ for all $j$. 
A straightforward computation gives that 
\begin{equation}
  \label{eq:H-eta-plus}
  (\Hcal^{b_+}(z),\eta_{\Cbb^d}(z))>0
\end{equation}
for all non zero $z\in\Cbb^d$. Consider now the 
following familly of symbols on $\Cbb^d$ 
$$
\sigma_t(z,v)=\Thom(\Cbb^d)(z,v-(t\Hcal^{b_+}+(1-t)\eta_{\Cbb^d})(z)), 
\quad (z,v)\in\T\Cbb^d.
$$
so that $\sigma_0=\Thom^\eta(\Cbb^d)$ and $\sigma_1=\Thom_{0,b_+}(\Cbb^d)$.
The inequality (\ref{eq:H-eta-plus}) shows that 
$\Char(\sigma_t)\cap\T_T\Cbb^d=\{0\}$ for all $t\in[0,1]$. 
Hence $\sigma_t$, $t\in[0,1]$, defines an homotopy of 
transversally elliptic symbols: $\Thom^\eta(\Cbb^d)=
\Thom_{0,b_+}(\Cbb^d)$ in $\K_T(\T_T\Cbb^d)$. $\Box$

For the remaining part of this paper we fix a scalar product on 
$\tgot^*$, and we consider the family of transversally elliptic 
symbols $\Thom_{\xi}(\Cbb^d)$, $\xi\in\tgot^*$ 
(to simplify, we do not mention the scalar product in the
notation). Proposition \ref{prop:indice-thom-xi} and
(\ref{eq:indice-thom-eta}) imply the following 

\begin{prop}\label{prop:N-R-indice}
  For every $\xi\in\tgot^*$, $\Qcal^T(\Cbb^d)$ is equal to 
the generalized character 
$\indice^T_{\Cbb^d}\left(\Thom_{\xi}(\Cbb^d)\right)$. 
\end{prop}

Now we apply the techniques developped in Section
\ref{sec:quantum-version} in order to compute the multiplicities of 
$\indice^T_{\Cbb^d}\left(\Thom_{\xi}(\Cbb^d)\right)$.


\subsection{Localization in a non-compact setting} \label{subsec:loc-non-compact}

Like in Section \ref{ssection.thom} we start with the 

\begin{defi}\label{def.thom.loc-C-d}
For any $\xi\in\tgot^*$ and any $T$-invariant {\em relatively compact} 
open subset $\Ucal\subset \Cbb^d$ 
we define the symbol $\Thom_{\xi}(\Ucal)$ by the relation
$$
\Thom_{\xi}(\Ucal)(z,v):=\Thom(\Cbb^d)(z,v-(\Hcal-\xi_{\Cbb^d})(z))
\quad (z,v)\in\T \Ucal.
$$
\end{defi}

The symbol $\Thom_{\xi}(\Ucal)$ is transversally elliptic when 
$\crpt(\parallel\Phi-\xi\parallel^2)\cap\partial\Ucal=\emptyset$ 
(the couple $(\Ucal,\xi)$ is called {\em good}) and we denote by  
    $$
    RR_{\Ucal}^{\xi}(\Cbb^d)\in R^{-\infty}(T)
    $$ 
its index. Proposition \ref{prop:U-xi} is still valid here. In
particular, for a good couple $(\Ucal,\xi)$, we have
$RR_{\Ucal}^{\xi'}(\Cbb^d)=RR_{\Ucal}^{\xi}(\Cbb^d)$ if $\xi'$ is close
enough to $\xi$. Consider now the decomposition
$$
\crpt(\parallel\Phi-\xi\parallel^2)=
\bigcup_{\gamma\in\Bcal_\xi}(\Cbb^d)^\gamma\cap\Phi^{-1}(\gamma+\xi).
$$
Here $\Bcal_\xi\subset\tgot^*$ is finite set since $\Cbb^d$ has a
finite number of stabilizer.  Since $0\in(\Cbb^d)^\gamma$
and $z\mapsto \langle\Phi(z),\gamma\rangle$ is constant on
$(\Cbb^d)^\gamma$, we have
\begin{equation}
  \label{eq:ps-gamma-xi}
  (\gamma +\xi,\gamma)=0
\end{equation}
for all $\gamma\in\Bcal_\xi$.

\begin{defi}\label{def:RR-xi-C}
 For any $\xi\in\tgot^*$ and $\gamma\in\Bcal_\xi$, we denote 
simply by 
$$
RR^\xi_\gamma(\Cbb^d)\in R^{-\infty}(T)
$$ 
the generalized character $RR_{\Ucal}^{\xi}(\Cbb^d)$, where  $\Ucal$
is a $T$-invariant relatively compact open neighborhood of 
$(\Cbb^d)^\gamma\cap\Phi^{-1}(\gamma+\xi)$ such that 
$\crpt(\parallel\Phi-\xi\parallel^2)\cap\overline{\Ucal}=
(\Cbb^d)^\gamma\cap\Phi^{-1}(\gamma+\xi)$.  
\end{defi}

Since $RR_{\Cbb^d}^{\xi}(\Cbb^d)$ is equal to $\Qcal^T(\Cbb^d)$ (see
Proposition \ref{prop:N-R-indice}), part $a)$ of Proposition
\ref{prop:U-xi} insures that we have the decomposition 
$$
   \Qcal^T(\Cbb^d)=\sum_{\gamma\in\Bcal_\xi} RR^\xi_\gamma(\Cbb^d).
$$

Let $\cgot\subset\tgot^*$ be a conic chamber of the cone $C(R)$, and
take $\xi$ in $\cgot$. Then $\xi$ is a regular value of the 
moment map $\Phi:\Cbb^d\to\tgot^*$ defined in (\ref{eq:moment-C}). 
Let $\Omega_\xi$ be the symplectic structure on the orbifold 
$(\Cbb^d)_\xi=\Phi^{-1}(\xi)/T$ 
that is induced from $\Omega$. The orbifold $(\Cbb^d)_\xi$ is also 
equipped with a {\em complex structure} $J_\xi$ that is induced 
from the standard complex structure on $\Cbb^d$, in such a way that  
the orbifold $((\Cbb^d)_\xi,\Omega_\xi,J_\xi))$ is a K\"ahler
orbifold. If $\xi$ belongs to the lattice $\Lambda^*$, the reduced 
space $(\Cbb^d)_\xi$ is the K\"ahler toric variety corresponding 
to the polytope $\{s\in(\Rbb^{\geq 0})^d\ \vert\ \sum s_j\alpha_j=\xi\}$ 
of $\Rbb^d$. For every $\mu\in\Lambda$ we consider the holomorphic 
orbifold line bundle
$$
\Lcal_{\xi,\mu}=(\Phi^{-1}(\xi)\times\Cbb_{-\mu})/T
$$
on $(\Cbb^d)_\xi$. 

\begin{defi}\label{def:Pol-C-d}
The periodic polynomial $P_\cgot:\Lambda^*\to\Zbb$ associated to 
the conic chamber $\cgot$ is given by  
\begin{equation}
  \label{eq:P-c-vector}
  P_\cgot(\mu)=RR((\Cbb^d)_\xi,\Lcal_{\xi,\mu}),
\end{equation}
where the right hand side is the Riemann-Roch number associated to the 
holomorphic orbifold line bundle $\Lcal_{\xi,\mu}$.
\end{defi}

Another way to define the periodic polynomial $P_\cgot$ is to consider
the generalized character $RR^\xi_0(\Cbb^d)$ for $\xi\in\cgot$: here
$\gamma=0$ parametrizes the component 
$\Phi^{-1}(\xi)$ of $\crpt(\parallel\Phi-\xi\parallel^2)$. Following 
(\ref{eq:RR-xi-O}) we have 
\begin{equation}
  \label{eq:pol-P-c}
 RR^\xi_0(\Cbb^d)=\sum_{\mu\in \Lambda^*}P_\cgot(\mu)\,\Cbb_\mu\quad 
{\rm in}\quad R^{-\infty}(T).
\end{equation} 
After Lemma \ref{lem:RR-xi-prime}, we know that
$RR^{\xi'}_0(\Cbb^d)=RR^\xi_0(\Cbb^d)$ when $\xi,\xi'$ are two elements
of $\cgot$: hence the polynomial $P_\cgot$ does not depend of the choice of 
$\xi$ in $\cgot$.

\subsection{Proof of Theorem \ref{prop:N-conic-chamber-SV}}
\label{proof:N-conic-chamber-SV}

Consider a weight $\mu\in(\cgot-\Box(\Phi))\cap\Lambda^*$ 
of the form $\mu=\xi'-\sum_j
t_j\alpha_j$ with ${\xi'}\in\cgot$ and $t_j\in[0,1]$. We start with the
decomposition 
$$
\Qcal^T(\Cbb^d)=
\sum_{\gamma\in\Bcal_{\xi'}} RR^{\xi'}_\gamma(\Cbb^d).
$$ 
Since $N_R(\mu)$ and $P_\cgot(\mu)$ are respectively the multiplicity of 
$\Cbb_\mu$ in $\Qcal^T(\Cbb^d)$ and in $RR^{\xi'}_0(\Cbb^d)$, 
the proof will be complete if we show that the multiplicity of 
$\Cbb_\mu$ in $RR^{\xi'}_\gamma(\Cbb^d)$ is equal to
zero when $\gamma\neq 0$. 

Consider a {\em non-zero} element $\gamma$ in $\Bcal_{\xi'}$. For the
character $RR^{\xi'}_\gamma(\Cbb^d)$ the localization 
(\ref{eq:loc-RR-beta-bis}) gives  
\begin{equation}
  \label{eq:loc-RR-gamma}
  RR^{\xi'}_\gamma(\Cbb^d)=
RR_0^{{\xi'}+\gamma}((\Cbb^d)^\gamma)\otimes
\left[\wedge_{\Cbb}^{\bullet}\overline{N}\right]^{-1}_{\gamma},
\end{equation}
where $N=\sum_{(\alpha_j,\gamma)\neq
  0}\Cbb_{-\alpha_j}$ corresponds to the normal bundle 
of $(\Cbb^d)^\gamma$ in $\Cbb^d$. The inverse 
$\left[\wedge_{\Cbb}^{\bullet}\overline{N}\right]^{-1}_{\gamma}$ is
equal to $(-1)^l\Cbb_{\delta(\gamma)}\otimes S^\bullet(N_\Cbb^{+,\gamma})$ 
where 
$$
\delta(\gamma)=-\sum_{(\alpha_j,\gamma)<0}\alpha_j. 
$$
Since $\gamma$ acts trivially on 
$(\Cbb^d)^\gamma$ all the weights $\mu'\in\Lambda^*$ that appear in
$RR_0^{{\xi'}+\gamma}((\Cbb^d)^\gamma)$ satisfy
$(\mu',\gamma)=0$. Since the weights of 
$N_\Cbb^{+,\gamma}$ are polarized by $\gamma$, we see from
(\ref{eq:loc-RR-gamma}) that all the weights $\mu'\in\Lambda^*$ that 
appear in $RR^{\xi'}_\gamma(\Cbb^d)$ must satisfy
\begin{equation}
  \label{eq:inegalite-mu-delta}
  (\mu',\gamma)\geq (\delta(\gamma),\gamma).
\end{equation}
Consider now the weight $\mu={\xi'}-\sum_j t_j\alpha_j$. Since
$\xi'\in\cgot$, the equality (\ref{eq:ps-gamma-xi}) implies  
$({\xi'},\gamma)<0$  and then
$$
(\mu,\gamma)=\underbrace{({\xi'},\gamma)}_{<0}+
\underbrace{\sum_{(\alpha_j,\gamma)>0}-
t_j(\alpha_j,\gamma)}_{\leq
0}-\sum_{(\alpha_j,\gamma)<0}t_j(\alpha_j,\gamma)<
-\sum_{(\alpha_j,\gamma)<0}(\alpha_j,\gamma).
$$
So we have proved that 
$(\mu,\gamma) < (\delta(\gamma),\gamma)$,
hence the multiplicity of $\Cbb_\mu$ in $RR^{\xi'}_\gamma(\Cbb^d)$ 
is equal to zero.  $\Box$

\subsection{Proof of Theorem \ref{theo:cas-C-d}}\label{proof-theo-cas-C-d}

Let $\cgot_\pm$ be two adjacent conic chambers and 
let  $\p\subset\tgot^*$ be the hyperplane 
that separates $\cgot_+$ and $\cgot_-$. Let $\beta\in\tgot$ be such that 
$\p=\{\xi\in\tgot^* \ \vert\ \langle\xi,\beta\rangle=0\}$ and 
$\cgot_\pm\subset \{\xi\in\tgot^* \ \vert\ \pm\langle\xi,\beta\rangle>0\}$. 

We consider two points $\xi_\pm\in \cgot_\pm$ such that $\xi=
\frac{1}{2}(\xi^+ + \xi^-)\in\p$  belongs to the relative interior 
$\cgot'$ of $\overline{\cgot_+}\cap\overline{\cgot_-}$. 
We suppose also that the orthogonal projection of $\xi_\pm$ on 
$\p$ are equal to $\xi$. We know that
$P_{\cgot_+}(\mu)-P_{\cgot_-}(\mu)$ 
is equal to the $\mu$-mutiplicity of $RR_0^{\xi_+}(\Cbb^d)-
RR_0^{\xi_-}(\Cbb^d)$. Proposition \ref{RR-0-diff} tells us that
$$
RR_0^{\xi_+}(\Cbb^d)- RR_0^{\xi_-}(\Cbb^d)=
RR_{\gamma}^{\xi_-}(\Cbb^d)-RR_{-\gamma}^{\xi_+}(\Cbb^d),
$$
where $\gamma\in\Rbb^{>0}\beta$ is such that $\xi_-+\gamma=\xi_+ -
\gamma=\xi$. The
localization (\ref{eq:loc-RR-beta-bis}) gives then  
\begin{equation}
  \label{eq:P-pm-rhs}
  RR_{\gamma}^{\xi_-}(\Cbb^d)-RR_{-\gamma}^{\xi_+}(\Cbb^d)=
RR_{0}^{\xi}((\Cbb^d)^\beta)\otimes
\left(\left[\wedge_{\Cbb}^{\bullet}\overline{N}\right]^{-1}_{\beta}-
\left[\wedge_{\Cbb}^{\bullet}\overline{N}\right]^{-1}_{-\beta}\right).
\end{equation}

The element $\xi$ belongs to the relative interior $\cgot'$ of 
$\overline{\cgot_+}\cap\overline{\cgot_-}$ which is a conic chamber
$\cgot'$ in $\p$ with respect to $R\cap\p$. Let $P_{\cgot'}:
\Lambda^*\cap\p\to\Zbb$ be the periodic polynomial map which coincides 
with the vector partition function $N_{R\cap\p}$ on 
$\overline{\cgot'}\cap\Lambda^*$. If we work with the vector space 
$(\Cbb^d)^\beta$ equipped with the hamiltonian action of $T/T_\p$, 
(\ref{eq:RR-xi-O}) gives the following equality in 
$R^{-\infty}(T/T_\p)$
\begin{equation}
  \label{eq:RR-xi-p-0}
  RR_{0}^{\xi}((\Cbb^d)^\beta)
=\sum_{\gamma\in\Lambda^*\cap\p}P_{\cgot'}(\gamma)\Cbb_\gamma. 
\end{equation}
A straightforward computation gives
\begin{equation}
  \label{eq:lambda-N-plus}
  \left[\wedge_{\Cbb}^{\bullet}\overline{N}\right]^{-1}_{\beta}=
(-1)^{r^-}\sum_{\mu\in\Lambda^*} N_{R'}(\mu+\delta^-)\Cbb_{\mu}
\end{equation}
and 
\begin{equation}
  \label{eq:lambda-N-moins}
 \left[\wedge_{\Cbb}^{\bullet}\overline{N}\right]^{-1}_{-\beta}=
(-1)^{r^+}\sum_{\mu\in\Lambda^*} N_{-R'}(\mu+\delta^+)\Cbb_{\mu}, 
\end{equation}
where $r^\pm,\delta^\pm,R'$ are defined in 
(\ref{eq:R-delta-prime}), (\ref{eq:delta-pm}) and (\ref{eq:r-pm}). 
Since $N_{-R'}(\mu)=N_{R'}(-\mu)$, the equations (\ref{eq:RR-xi-p-0}), 
(\ref{eq:lambda-N-plus}) and (\ref{eq:lambda-N-moins}) show that the
RHS of (\ref{eq:P-pm-rhs}) is equal to 
$$
\sum_{\mu\in\Lambda^*}\sum_{\gamma\in\Lambda^*\cap\p}
D(\mu)P_{\cgot'}(\gamma)\Cbb_{\mu+\gamma}
$$
with
$D(\mu)=(-1)^{r^-}N_{R'}(\mu+\delta^-)-(-1)^{r^+}N_{R'}(-\mu-\delta^+)$.
Finally we have proved that $P_{\cgot_+}(\mu)-P_{\cgot_-}(\mu)=
\sum_{\gamma\in\Lambda^*\cap\p}D(\mu-\gamma)P_{\cgot'}(\gamma)$. $\Box$

\end{document}